\documentclass[12pt]{amsart}
 
\usepackage[dvips, bookmarks=false, hyperfigures=false,
setpagesize=false]{hyperref}

\usepackage[mathscr]{eucal}
\usepackage{t-angles}

\DeclareMathAlphabet{\mathpzc}{OT1}{pzc}{m}{it}

\advance\textheight\topskip
\usepackage{mathabx}
\usepackage{mathptmx}

\usepackage[matrix,arrow,curve,dvips]{xy} 
\UsePSspecials{dvips}

\allowdisplaybreaks

\numberwithin{equation}{section}
\numberwithin{figure}{section}
\makeatletter
\@addtoreset{equation}{section}
\@addtoreset{figure}{section}
\@addtoreset{subsubsection}{section}

\def\@secnumfont{\bfseries}
\def\subsubsection{\@startsection{subsubsection}{3}%
  \z@{.5\linespacing\@plus.7\linespacing}{-.5em}%
  {\normalfont\bfseries}}
\def\paragraph{\@startsection{paragraph}{4}%
  \z@\z@{-\fontdimen2\font}%
  \normalfont\bfseries}
\def\subparagraph{\@startsection{subparagraph}{5}%
  \z@\z@{-\fontdimen2\font}%
  \normalfont\bfseries}
\makeatother

\usepackage[OT2,T1]{fontenc}
\newcommand{\mifody}{%
  \renewcommand\rmdefault{wncyr}%
  \renewcommand\sfdefault{wncyss}%
  \renewcommand\encodingdefault{OT2}%
  \normalfont
  \selectfont}
\usepackage{times}

\newcommand{\UresSL}[1]{\overline{\mathscr{U}}_{\q} s\ell(#1)}
\newcommand{\theaut}{\alpha}

\newcommand{\VertexII}[2]{V^{#2}_{#1,\,0}}
\newcommand{\VertexIV}[4]{V^{\{#2,\,#4\}}_{#1,\,#3,\,0}}
\newcommand{\VertexVI}[6]{V^{\{#2,\,#4,\,#6\}}_{#1,\,#3,\,#5,\,0}}
\newcommand{\bF}{\bar{F}}

\newcommand{\scr}{F}
\newcommand{\thev}[1]{\scr(#1)}

\newcommand{\lS}{\mathscr{L}}

\newcommand{\modO}{\mathscr{O}}
\newcommand{\modP}{\mathscr{P}}
\newcommand{\modV}{\mathscr{V}}
\newcommand{\modW}{\mathscr{W}}
\newcommand{\modX}{\mathscr{X}}

\newcommand{\lvoa}{\mathcal{L}}
\newcommand{\modVF}{\mathcal{V}}
\newcommand{\cF}{\mathcal{F}}

\newcommand{\Nich}{\mathfrak{B}}

\newcommand{\cK}{\mathbb{X}}
\newcommand{\confS}{\cK}

\newcommand{\FYDBraid}{\mathcal{B}}
\newcommand{\Fbraid}{\mathcal{B}}

\newcommand{\FromLeftiiD}{\,\;\begin{tangles}{l}
    \hstr{50}\vstr{50}\lu
    \object{\raisebox{3.3pt}{\tiny$\bullet$}}
  \end{tangles}\;\;}

\newcommand{\FromRightiiD}{\;\;\begin{tangles}{l}
    \hstr{50}\vstr{50}\object{\raisebox{3.3pt}{\tiny$\bullet$}}\ru
  \end{tangles}\;\,}

\newcommand{\leftact}{\kern1pt{\rightharpoonup}\kern1pt}
\newcommand{\rightact}{\kern1pt{\leftharpoonup}\kern1pt}
\newcommand{\leftreg}{\kern1pt{\rightharpoonup}\kern1pt}
\newcommand{\rightreg}{\kern1pt{\leftharpoonup}\kern1pt}
\newcommand{\rightregX}{\kern1pt{\leftbarharpoon}\kern1pt}
\newcommand{\leftregii}{\kern2pt{\looparrowright}\kern1pt}
\newcommand{\rightregii}{\kern1pt{\looparrowleft}\kern1pt}

\newcommand{\leftregU}{\kern1pt{\rightharpoondown}\kern1pt}
\newcommand{\rightregU}{\kern1pt{\leftharpoondown}\kern1pt}
\newcommand{\rightregUX}{\kern1pt{\barleftharpoon}\kern1pt}
\newcommand{\leftregUii}{\kern2pt{\looparrowdownright}\kern1pt}
\newcommand{\rightregUii}{\kern1pt{\looparrowdownleft}\kern2pt}

\newcommand{\cU}{\mathscr{U}}
\newcommand{\eval}[2]{\langle#1,\,#2\rangle\,}

\newcommand{\oZ}{\mathbb{Z}}
\newcommand{\mint}[2]{\idotsint\limits_{#1}#2}

\newcommand{\q}{\mathfrak{q}}

\newcommand{\aint}[1]{\langle#1\rangle}
\newcommand{\qInt}[1]{\langle#1\rangle}

\newcommand{\abin}[2]{\mathchoice%
  {\abinomm{#1}{#2}}{\abinommm{#1}{#2}}%
  {\abinommm{#1}{#2}}{\abinommm{#1}{#2}}}
\newcommand{\abinomm}[2]{\mbox{\footnotesize$\displaystyle
    \genfrac{\langle}{\rangle}{0pt}{}{#1}{#2}$}}
\newcommand{\abinommm}[2]{\genfrac{\langle}{\rangle}{0pt}{}{#1}{#2}}

\newcommand{\Aint}[1]{\langle#1\rangle}
\newcommand{\Afac}[1]{\langle#1\rangle!\,}
\newcommand{\Abin}[2]{\mathchoice%
  {\Abinm{#1}{#2}}{\Abinmm{#1}{#2}}%
  {\Abinmm{#1}{#2}}{\Abinmm{#1}{#2}}}
\newcommand{\Abinm}[2]{\mbox{\footnotesize$\displaystyle
    \genfrac{\langle}{\rangle}{0pt}{}{#1}{#2}$}}
\newcommand{\Abinmm}[2]{\genfrac{\langle}{\rangle}{0pt}{}{#1}{#2}}

\newcommand{\alphaminus}{-\sqrt{\kern-1pt\frac{2}{p}}}

\newcommand{\Smash}{\mathbin{\hash}}

\newcommand{\set}[1]{\mathbf{#1}}

\newcommand{\cV}{\mathscr{V}}

\newcommand{\bup}[1]{^{{\scriptscriptstyle\!\underline{\,#1_{\vphantom{.}}\,}\!}}}

\newcommand{\zero}{_{_{(0)}}}
\newcommand{\mone}{_{_{(-1)}}}
\newcommand{\bzero}{_{_{\,\underline{\,0\,}}}}
\newcommand{\bmone}{_{_{\underline{\!{-1}\!}}}}

\newcommand{\KK}{\mathcal{K}\kern-5.7pt\raisebox{-3.9pt}{\footnotesize\textit{2}}\,}
\newcommand{\DD}{\mathscr{D}}
\newcommand{\cA}{\mathscr{A}}

\newcommand{\cross}{\textstyle\!\!{\times}\!\!}
\newcommand{\punct}{\textstyle{\circ}}

\newcommand{\dotact}{\mathbin{\pmb{.}}}

\newcommand{\cR}{\mathscr{R}}

\newcommand{\oC}{\mathbb{C}}

\newcommand{\id}{\mathrm{id}}
\newcommand{\tensor}{\otimes}

\newcommand{\cH}{\mathscr{H}}

\newcommand{\ccirc}{\mathbin{\raisebox{1pt}{\,$\scriptscriptstyle\circ$\,}}}

\newcommand{\catC}{\category{C}}
\newcommand{\A}{\raisebox{.5pt}{\large$\mathpzc{S}\kern-1pt$}}
\newcommand{\medA}{\mathpzc{S}}
\newcommand{\hA}{\mbox{\large$\mathpzc{s}\kern-.8pt$}}

\newcommand{\YDname}{\mathcal{Y\kern-3ptD}}
\newcommand{\HHyd}{{}\mbox{\small${}^{H}_{H}$}\YDname}
\newcommand{\Cepx}{\mathscr{C}}

\newcommand{\Shuffle}{\mathop{\text{\mifody\sf Sh}}\nolimits}

\newcommand{\Shift}[2]{#2^{\uparrow#1}}
\newcommand{\xShift}[2]{#2}
\newcommand{\shift}{\uparrow}

\newcommand{\Bbin}[2]{{\Shuffle^{#1}_{#2}}}
\newcommand{\Bfac}[1]{\mathfrak{S}_{#1}}
\newcommand{\Bint}[1]{\{#1\}}

\newcommand{\Yspace}{(\!(Y)\!)}
\newcommand{\Zspace}{(\!(Z)\!)}

\newcommand{\YDspace}{(\!(Y|}
\newcommand{\ZDspace}{(\!(Z|}
\newcommand{\UDspace}{(\!(U|}

\newcommand{\fobject}[1]{\object{\scriptstyle #1}}

\newcommand{\FromLeft}[1]{\mathcal{L}_{#1}}
\newcommand{\FromLeftii}[1]{\,\widetilde{\!\mathcal{L}\!}\,_{#1}}
\newcommand{\FromRight}[1]{\mathcal{R}_{#1}}
\newcommand{\FromRightii}[1]{\widetilde{\mathcal{R}}_{#1}}

\newcommand{\adjoint}{%
  \mathchoice{\mathbin{\blacktriangleright}}%
  {\mathbin{\mbox{\small${\blacktriangleright}$}}}%
  {\mathbin{{\blacktriangleright}}}%
  {\mathbin{{\blacktriangleright}}}}

\newcommand{\Bbraid}{\pmb{\Psi}}
\newcommand{\Adja}{\mathop{\mathrm{Ad}}\nolimits}
\newcommand{\Adj}{\mathop{\mathrm{Ad}}\nolimits}

\newcommand{\category}{\mathcal}

\newcommand{\deltaL}{\delta_{_{\text{L}}}}
\newcommand{\deltaR}{\delta_{_{\text{R}}}}

\newcommand{\ffrac}[2]{\raisebox{.5pt}{\mbox{\footnotesize$\displaystyle\frac{#1}{#2}$}}}

\newcommand{\oN}{\mathbb{N}}
\newcommand{\End}{\mathop{\mathrm{End}}\nolimits}

\newcommand{\cX}{\mathscr{X}}
\newcommand{\cY}{\mathscr{Y}}
\newcommand{\cZ}{\mathscr{Z}}
\newcommand{\mapX}{\mathop{\text{\mifody\sf Kh}}}
\newcommand{\mapI}{\mathop{\text{\mifody\sf I}}}

\newcommand{\Left}{{_{\text{L}}}}
\newcommand{\Right}{{_{\text{R}}}}

\newcommand{\leftii}{%
  \mathchoice%
  {\mathbin{\mbox{\footnotesize${\vartriangleright}$}}}%
  {\mathbin{\scriptstyle\vartriangleright}}%
  {\mathbin{{\vartriangleright}}}%
  {\mathbin{{\vartriangleright}}}}

\newcommand{\rightii}{\mathbin{\mbox{\footnotesize${\vartriangleleft}$}}}

\newcommand{\bref}[1]{\textbf{\ref{#1}}}

\swapnumbers

\newtheorem{lemma}[subsubsection]{Lemma}

\newtheorem{prop}[subsubsection]{Proposition}

\theoremstyle{definition}

\newtheorem{Rem}[subsection]{Remark}

\begin{document}

\title[The Nichols algebra of screenings
  ]{The Nichols algebra of screenings}

\author[Semikhatov]{A.M.~Semikhatov}

\author[Tipunin]{I.Yu.~Tipunin}

\address{Lebedev Physics Institute\hfill\mbox{}\linebreak
  \texttt{ams@sci.lebedev.ru}, \ \texttt{tipunin@pli.ru}}

\begin{abstract}
  Two related constructions are associated with screening operators in
  models of two-dimensional conformal field theory.  One is a local
  system constructed in terms of the braided vector space $X$ spanned
  by the screening species in a given CFT{} model and the space of
  vertex operators~$Y$ and the other is the Nichols algebra $\Nich(X)$
  and the category of its Yetter--Drinfeld modules, which we propose
  as an algebraic counterpart, in a ``braided'' version of the
  Kazhdan--Lusztig duality, of the representation category of
  vertex-operator algebras realized in logarithmic CFT{} models.
\end{abstract}

\maketitle
\thispagestyle{empty}

\begin{footnotesize}\addtolength{\baselineskip}{-4pt}
  \tableofcontents
\end{footnotesize}

\section{Introduction}
We discuss how braided categories\,---\,more specifically,
representation categories of braided Hopf algebras\,---\,can be
associated with nonsemisimple (logarithmic) models of conformal field
theory (CFT).  The subject of braided categories related to CFT{} has
been profoundly investigated in~\cite{[FRS],[FFRS],[FS]}, but
considering the problem at a somewhat different angle may not be
altogether unworthy.  The proposal in this paper is to take the
algebra of screening operators, regard it as a braided Hopf algebra,
and then construct a category of its Yetter--Drinfeld modules from
other CFT{} data, the vertex operators in a given model.  The algebras
generated by screening operators are in fact Nichols
algebras~\cite{[Nich],[AG],[Andr-remarks],[U-Nich],[Heck-Weyl],
  [Heck-2],[AHS],[HS-LMS],[Heck-class]}---universal braided Hopf
algebra quotients of tensor algebras that have received considerable
attention recently, originally motivated by Andruskiewitsch and
Schneider's classification program (see
\cite{[AS-pointed],[AS-onthe],[ARS]} and the references therein).  We
here extend their use to a bridge between CFT{} and braided
categories.

\subsubsection*{Nichols algebras}
Under the name of ``bialgebras of type one,'' Nichols algebras (more
precisely, their bosonizations) originally appeared in \cite{[Nich]}.
They have several definitions, whose equivalence is due
to~\cite{[Sch-borel]} and~\cite{[AG]} (where they still feature under
a more indigenous name).  In addition to the papers cited above, they
also appeared
in~\cite{[Lu-intro],[KT],[Lu-intro],[J],[KS],[Sch-borel],[Rosso-inv],
  [B]} (see~\cite{[Ag-all]} and the references therein for recent
progress).\footnote{An important technicality, which we nevertheless
  tend to ignore, is that there is a distinction between quantum
  symmetric algebras~\cite{[Rosso-inv]} and Nichols algebras proper;
  the latter are selected by the condition that the braiding is
  \textit{rigid} (which in particular guarantees that the duals $X^*$
  are objects in the same braided category with the $X$);
  see~\cite{[A-about],[G-free]}.}

The Nichols algebra $\Nich(X)$ of a braided linear space $X$ is a
graded braided Hopf algebra,
$\Nich(X)=\bigoplus_{n\geq0}\Nich(X)^{(n)}$ (a vector-space direct
sum), such that $\Nich(X)^{(0)}$ is merely the ground field and
$\Nich(X)^{(1)}=X$, and this last space has two properties: it
coincides with the space of \textit{all primitive elements}
$P(X)=\{x\in\Nich(X) \mid\Delta x=x\tensor 1 + 1\tensor x\}$ and it
\textit{generates all of $\Nich(X)$} as an algebra.

The Nichols algebras occurred independently in~\cite{[Wor]}, for the
purposes of constructing a quantum differential calculus and in a form
suggestive of their role of ``fully braided generalizations'' of
symmetric algebras, viz.
\begin{equation*}
  \Nich(X)
  = k\oplus X\oplus\bigoplus_{n\ge 2} X^{\otimes n}/\ker\Bfac{n},
\end{equation*}
where $\Bfac{n}$ is the total braided symmetrizer.  This gives another
characterization of Nichols algebras.

Finally, the Nichols algebra of a linear braided space $X$ can be
defined in terms of a duality pairing
$\eval{{\cdot}\,}{{\cdot}}:T(X^*)\tensor T(X)\to k$.  It follows that
$\Nich(X) = T(X)/I(X)$, where $I(X)$ is the kernel of the
pairing~\cite{[AG]}.

\subsubsection*{Motivation}
Nichols algebras are a core element in the construction of the
nilpotent part of deformed enveloping algebras
$\mathscr{U}^{+}_q(\mathfrak{g})$~\cite{[Lu-intro]} (also
see~\cite{[Rosso-CR],[Rosso-inv]}).  This last development has a
well-known ``physical'' reformulation, stating that quantum Serre
relations are satisfied by screening operators in
CFT~\cite{[Rosso-VO]} (also see the references therein,
\cite{[BMP-fock]} in particular).  In a conformal model, the action of
screenings on vertex operators does not yield elements of the
vertex-operator algebra of that model in general, but generates a
larger structure, whose properties have not been studied much just
because it is not a vertex-operator algebra and, in addition, carries
a subtle dependence on the choice of contours involved in the
definition of screenings.  Objects of this sort, sometimes called
screened (or dressed) vertex operators, are indispensable, however, in
constructing correlation functions in CFT~\cite{[DF]}.  We use them in
what follows to construct the categories of Hopf bimodules and
Yetter--Drinfeld modules over Nichols algebras of screenings.

This construction of (co)modules over an algebra of screenings is
closely related to a more geometric part of this paper---a local
system (or, a flat bundle) on the configuration space $\cK_{m,n}$ of
$n$ points on an $m$-punctured complex plane.  The points are to be
thought of as positions of screenings, and the punctures as the
positions of vertex operators.  The latter are fixed, and the former,
in physical terms, are ``integrated over'' along some contours.  The
idea of screening operators being represented by iterated contour
is refined by a ``space--time'' decomposition of the complex plane,
with the integrations running in the ``space'' direction, with which
we associate a stratification of the $\cK_{m,n}$ configuration space
of screenings.  The stratification readily leads to the local system
mentioned above, and the homology complex with coefficients in this
local system gives rise to a shuffle algebra, in which the Nichols
algebra $\Nich(X)$ is a subalgebra.  The local system becomes an
abundant source of Yetter--Drinfeld
$\Nich(X)$-modules.\enlargethispage{\baselineskip}

\subsubsection*{Outline of the construction}
In most of the paper, we work with much more general braided spaces
and categories than those observed\pagebreak[3] in CFT{} models, and
only specialize to diagonal braiding (and then, actually, to a very
particular case) in the Sec.~\ref{sec:1p}.

\begin{figure}[htb]
  \centering
  \begin{equation*}
    \xymatrix@R=4pt@C=80pt{
      \ar@{--}|(.2){\cross}|(.7){\cross}[0,2]&&\\
      \ar@*{[|(1.6)]}@{-}|(.1){\cross}|(.4){\punct}[0,2]&&\\
      \ar@*{[|(1.6)]}@{-}|(.6){\punct}[0,2]&&\\
      \ar@{--}|{\cross}[0,2]&&\\
      \ar@{--}|(.3){\cross}|(.7){\cross}[0,2]&&\\
      \ar@*{[|(1.6)]}@{-}|(.25){\punct}|(.45){\cross}|(.65){\cross}|(.8){\cross}[0,2]&&
      }
  \end{equation*}
  \caption{\small``Fixed'' parallel lines are drawn through the
    punctures ($\punct$).  Additional lines can move in between so as
    to remain parallel and not collide with or pass through one
    another or the fixed lines.  All lines are populated with
    screenings~(${\times}$).}
  \label{fig:intro}
\end{figure}
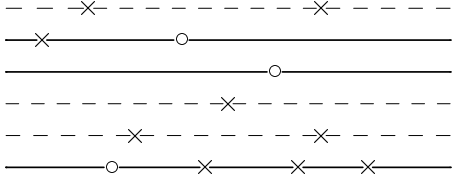

Screening operators, in the way they act on vertex operators, are
given by contour integrals in one form or another, and are in this
sense ``nonlocal.''  We abstract the contour integral representation
to a family parallel lines---``spatial directions on the complex
plane''---drawn through each puncture (solid lines in
Fig.~\ref{fig:intro}).
We also draw arbitrarily many additional parallel lines, not running
through the punctures, and populate them with screenings.  Once two
screenings (shown by crosses in~Fig.~\ref{fig:intro}) are on the same
line, the notion of precedence is clearly defined, and hence so is the
notion of the \textit{transposition} of any two screenings via
braiding.  The braiding is inherited from a given CFT{} model as a map
\begin{equation*}
  \Psi:X\tensor X\to X\tensor X
\end{equation*}
that satisfies the braiding equation, where $X$ is the linear span of
the different screening species.  Somewhat more abstractly---without
direct reference to CFT{}---we just assume that a braided linear
space~$X$ is associated with every cross in a picture such as in
Fig.~\ref{fig:intro}.

Placing $n$ crosses on the system of lines on an $m$-punctured complex
plane, as in Fig.~\ref{fig:intro},
defines a stratification of the configuration space of
crosses/screenings $\cK_{m,n}$.  A stratum---``cell''---is defined by
any particular distribution of $n$ crosses over such a system of
lines, with the movable lines (from $0$ to $n$ in number) allowed to
move so as to remain parallel and not collide with or pass through one
another or the fixed lines, and the crosses allowed to slide along the
lines so as to not collide with one another or with the punctures and
not pass through the punctures.

With each puncture, we associate another braided space $Y$ (in CFT{},
a linear span of vertex operators), and let $\Psi$ denote the braiding
in all cases $X\tensor Y\to Y\tensor X$, etc.  Notably, an essential
part of our construction does \textit{not} require the braiding
$Y\tensor Y\to Y\tensor Y$.

Given $X$ and $Y$, we define a \textit{local system} (the space of
sections of a flat bundle) over the configuration space $\confS_m$ of
(noncoincident and indistinguishable) points on the complex plane with
$m$ punctures.  For a fixed number $n$ of points, we take all cells,
i.e., all possibilities to distribute $n$ crosses over $m+\ell$ lines,
$0\leq\ell\leq n$.  A braided linear space linearly isomorphic
to~$X^{\tensor n}\tensor Y^{\tensor m}$ is then associated with each
cell.  The local system is defined by specifying how these spaces are
identified under the ``restriction map'' associated with taking the
boundary of each cell.

\begin{figure}[tbp]
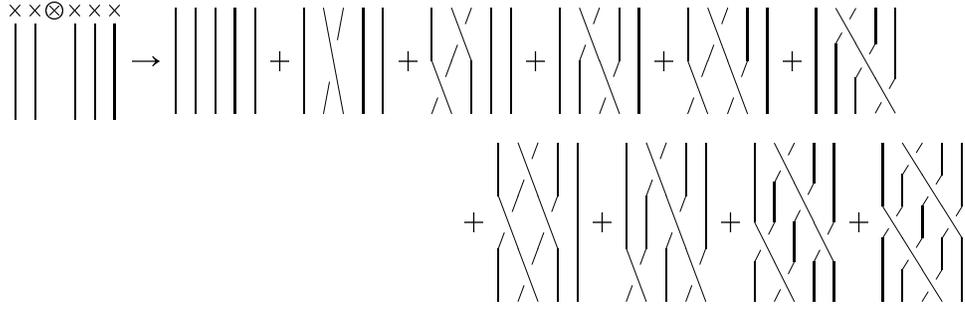

  \centering
  \begin{multline*}
    \begin{tangles}{l}
      \hstr{75}\fobject{{\times}}\step\fobject{{\times}}\step
      \fobject{\otimes}\step\fobject{{\times}}\step
      \fobject{{\times}}\step\fobject{{\times}}\\[-2pt]
      \hstr{75}\vstr{180}\id\step\id\step[2]\id\step\id\step\id
    \end{tangles}\
    \to\
    \begin{tangles}{l}
      \hstr{75}\vstr{200}\id\step\id\step\id\step\id\step\id
    \end{tangles} \
    +  \ \begin{tangles}{l}
      \hstr{75}\vstr{200}\id\step[1]\hx\step[1]\id\step[1]\id\\
    \end{tangles}\ +  \ \begin{tangles}{l}
      \hstr{75}\vstr{100}\id\step[1]\hx\step[1]\id\step[1]\id\\
      \hstr{75}\vstr{100}\hx\step[1]\id\step[1]\id\step[1]\id\\
    \end{tangles}\  +  \ \begin{tangles}{l}
      \hstr{75}\vstr{100}\id\step[1]\hx\step[1]\id\step[1]\id\\
      \hstr{75}\vstr{100}\id\step[1]\id\step[1]\hx\step[1]\id\\
    \end{tangles}\  +  \ \begin{tangles}{l}
      \hstr{75}\vstr{100}\id\step[1]\hx\step[1]\id\step[1]\id\\
      \hstr{75}\vstr{100}\hx\step[1]\hx\step[1]\id\\
    \end{tangles}\  +  \ \begin{tangles}{l}
      \hstr{75}\vstr{66}\id\step[1]\hx\step[1]\id\step[1]\id\\
      \hstr{75}\vstr{66}\id\step[1]\id\step[1]\hx\step[1]\id\\
      \hstr{75}\vstr{66}\id\step[1]\id\step[1]\id\step[1]\hx\\
    \end{tangles}
    \\[4pt]
    +  \ \begin{tangles}{l}
      \hstr{75}\vstr{100}\id\step[1]\hx\step[1]\id\step[1]\id\\
      \hstr{75}\vstr{100}\hx\step[1]\hx\step[1]\id\\
      \hstr{75}\vstr{100}\id\step[1]\hx\step[1]\id\step[1]\id\\
    \end{tangles}\  +  \ \begin{tangles}{l}
      \hstr{75}\vstr{100}\id\step[1]\hx\step[1]\id\step[1]\id\\
      \hstr{75}\vstr{100}\id\step[1]\id\step[1]\hx\step[1]\id\\
      \hstr{75}\vstr{100}\hx\step[1]\id\step[1]\hx\\
    \end{tangles} \ +  \ \begin{tangles}{l}
      \hstr{75}\vstr{75}\id\step[1]\hx\step[1]\id\step[1]\id\\
      \hstr{75}\vstr{75}\id\step[1]\id\step[1]\hx\step[1]\id\\
      \hstr{75}\vstr{75}\hx\step[1]\id\step[1]\hx\\
      \hstr{75}\vstr{75}\id\step[1]\hx\step[1]\id\step[1]\id\\
    \end{tangles}\  +  \ \begin{tangles}{l}
      \hstr{75}\vstr{60}\id\step[1]\hx\step[1]\id\step[1]\id\\
      \hstr{75}\vstr{60}\id\step[1]\id\step[1]\hx\step[1]\id\\
      \hstr{75}\vstr{60}\hx\step[1]\id\step[1]\hx\\
      \hstr{75}\vstr{60}\id\step[1]\hx\step[1]\id\step[1]\id\\
      \hstr{75}\vstr{60}\id\step[1]\id\step[1]\hx\step[1]\id
    \end{tangles}\qquad\qquad
  \end{multline*}
  \caption{\small The $X^{\otimes r}\tensor X^{\otimes s}\to
    X^{\otimes(r+s)}$ shuffle product (with $r=2$ and $s=3$ in the
    picture) involves $\binom{r+s}{r}$ terms, each of which represents
    a braid group element whose action arranges the crosses ``fed in''
    to the top of each strand (this and other braid diagrams are to be
    read from top down).  Each term shows precisely how the $2+3$
    crosses are sent into their ``target'' positions.}
  \label{fig:shuffle-intro}
\end{figure}

Taking the boundary means that a movable line in Fig.~\ref{fig:intro}
merges with another (movable or fixed) line, and the crosses carried
by the two lines are ``collectivized.''  For example, two movable
lines
\begin{align*}
  \xymatrix@1@C=80pt{ \ar@{--}|(.33){\cross}|(.66){\cross}[r]& }\qquad
  &\text{\small($r$ crosses)}
  \\*[-4pt]
  \xymatrix@1@C=100pt{
    \ar@{--}|(.2){\cross}|(.5){\cross}|(.7){\cross}[r]& } \quad
  &\text{\small($s$ crosses)}
  \\*[-6pt]
  \intertext{merge into} \xymatrix@1@C=100pt{
    \ar@{--}|(.16){\cross}|(.33){\cross}|(.50){\cross}|(.67){\cross}|(.85){\cross}[r]&
  } \quad &\text{\small($r+s$ crosses)}
\end{align*}
The restriction maps are naturally constructed in terms of the
braiding, and, expectedly, are given by ``quantum''
shuffles~\cite{[Rosso-inv]} (see Appendix~\ref{app:shuffles}). \ In a
merger of two movable lines, for example, the crosses are
collectivized by $(r,s)$-shuffle permutations lifted to the braid
group algebra.  This is illustrated
in Fig.~\ref{fig:shuffle-intro}.

This ``merger'' operation is associative and is part of a braided
bialgebra structure~\cite{[Rosso-inv]}.  The coproduct in this
``braided Hopf algebra of crosses'' is given by splitting each line
into two and, for a line with $j$ crosses, taking the sum over all the
$j+1$ possibilities of distributing the $j$ indistinguishable crosses
between two lines,
\begin{multline}\label{coproduct-intro}
  \xymatrix@R=6pt@C=80pt{
    \ar@{--}|(.25){\cross}|(.5){\cross}|(.75){\cross}[r]&
  }\longrightarrow
  \xymatrix@R=6pt@C=80pt{
    \ar@{==}|(.25){\cross}|(.5){\cross}|(.75){\cross}[r]&
  }\longrightarrow
  \\*
  \xymatrix@R=2pt@C=60pt{
    \ar@{--}|(.25){\cross}|(.5){\cross}[r]|(.75){\cross}&\\
    \ar@{--}[r]&
  }+
  \xymatrix@R=2pt@C=60pt{
    \ar@{--}|(.25){\cross}|(.75){\cross}[r]&\\
    \ar@{--}|(.5){\cross}[r]&
  }
  +\xymatrix@R=2pt@C=60pt{
    \ar@{--}|(.5){\cross}[r]&\\
    \ar@{--}|(.25){\cross}|(.75){\cross}[r]&
  }
  +\xymatrix@R=2pt@C=60pt{
    \ar@{--}[r]&\\
    \ar@{--}|(.25){\cross}|(.5){\cross}|(.75){\cross}[r]&
  }
\end{multline}

\smallskip

\noindent
(the ``deconcatenation'' coproduct).  The product and coproduct define
a braided bialgebra~\cite{[Rosso-inv]}.  \textit{The antipode given by
  orientation reversal} makes it into a braided Hopf algebra,
$\cH(X)$.  Those $r$-cross elements
\begin{equation}\label{r-crosses}
  \xymatrix@C=90pt{
      \ar@{--}|(.2){\cross}|(.5){\cross}|(.75){\cross}[r]&
    }
\end{equation}
that can be obtained from a single-cross one by iterated
multiplication span the Nichols algebra~$\Nich(X)\subset\cH(X)$.

The fixed lines (those with a puncture each), furthermore, become
elements of $\cH(X)$ (and $\Nich(X)$) bi(co)modules.  The left and
right \textit{actions} of elements of $\cH(X)$ (movable lines, such
as~\eqref{r-crosses}) on fixed lines (such as
$\xymatrix@1@C=24pt{\ar@*{[|(1.6)]}@{-}|(.15){\cross}
  |(.4){\cross}|(.65){\punct}|(.90){\cross}[0,2]&&}$) amount to
distributing the ``new'' crosses over the entire fixed line, again by
quantum shuffles, as is illustrated in Fig.~\ref{fig:action-intro}
(where the visualization conventions are somewhat different from those
in Fig.~\ref{fig:shuffle-intro}).  The left and right
\textit{coactions} are by deconcatenation of the crosses from the
half-line on the respective side of the puncture.  The resulting
bimodule bicomodule actually turns out to be a Hopf bimodule
(``bi-Hopf module,'' ``bicoviariant bimodule'') over~$\cH(X)$.
\begin{figure}[tbh]
  \centering
  \begin{multline*}
    \xymatrix@R=6pt@C=30pt{
      \ar@{--}|(.5){\cross}[0,2]&&
    }
    \dotact
    \xymatrix@R=6pt@C=30pt{
      \ar@*{[|(1.6)]}@{-}|(.3){\punct}|(.65){\cross}[0,2]&&
    }
    =
    \\[6pt]
    \xymatrix@R=6pt@C=30pt{
      \ar@*{[|(1.6)]}@{-}|(.1){\cross}|(.3){\punct}|(.70){\cross}[0,2]&&
    } + 
    \xymatrix@R=6pt@C=30pt{
      \ar@*{[|(1.6)]}@{-}|(.3){\punct}|(.5){\cross}|(.75){\cross}[0,2]
      \ar@/^10pt/[r]&&
    } + 
    \xymatrix@R=6pt@C=30pt{
      \ar@*{[|(1.6)]}@{-}|(.3){\punct}|(.6){\cross}|(.95){\cross}[0,2]
      \ar@/^12pt/[rr]&&
    }
  \end{multline*}
  \caption{\small The left action $\dotact$ of ``the Hopf algebra of
    crosses'' on its module.  The arrows indicate the braiding used to
    put the ``new'' cross into the different positions.
    (Alternatively, the braiding can be represented in terms of
    diagrams of the type of those in Fig.~\ref{fig:shuffle-intro}.)
    Under the right action, the ``new'' crosses populate the fixed
    line starting from the right.}
  \label{fig:action-intro}
\end{figure}

\begin{figure}[t]
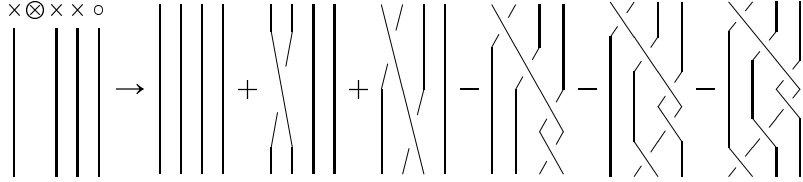

  \centering
  \begin{equation*}
    \begin{tangles}{l}
      \hstr{80}\fobject{{\times}}\step\fobject{\otimes}\step
      \fobject{{\times}}\step\fobject{{\times}}\step\fobject{{\circ}}\\
      \hstr{80}\vstr{280}\id\step[2]\id\step\id\step\id
    \end{tangles}\
    \to\
    \begin{tangles}{l}
      \\
      \hstr{80}\vstr{320}\id\step\id\step\id\step\id
    \end{tangles}
    \;+
    \
    \begin{tangles}{l}
      \hstr{80}\vstr{50}\id\step\id\step\id\step\id\\
      \hstr{80}\vstr{220}\hx\step\id\step\id\\
      \hstr{80}\vstr{50}\id\step\id\step\id\step\id\\
    \end{tangles}
    \;+
    \;
    \begin{tangles}{l}
      \hstr{80}\vstr{160}\hx\step\id\step\id\\
      \hstr{80}\vstr{160}\id\step\hx\step\id
    \end{tangles}
    \;-
    \;
    \begin{tangles}{l}\vstr{90}
      \hstr{90}\vstr{80}\hx\step\id\step\id\\
      \hstr{90}\vstr{80}\id\step\hx\step\id\\
      \hstr{90}\vstr{80}\id\step\id\step\hx\\
      \hstr{90}\vstr{80}\id\step\id\step\hx
    \end{tangles}
    \;-
    \;
    \begin{tangles}{l}
      \hstr{90}\vstr{66}\hx\step\id\step\id\\
      \hstr{90}\vstr{66}\id\step\hx\step\id\\
      \hstr{90}\vstr{66}\id\step\id\step\hx\\
      \hstr{90}\vstr{66}\id\step\id\step\hx\\
      \hstr{90}\vstr{66}\id\step\hx\step\id
    \end{tangles}
    \;-
    \;
    \begin{tangles}{l}
      \hstr{90}\vstr{55}\hx\step\id\step\id\\
      \hstr{90}\vstr{55}\id\step\hx\step\id\\
      \hstr{90}\vstr{55}\id\step\id\step\hx\\
      \hstr{90}\vstr{55}\id\step\id\step\hx\\
      \hstr{90}\vstr{55}\id\step\hx\step\id\\
      \hstr{90}\vstr{55}\hx\step\id\step\id\\
    \end{tangles}
  \end{equation*}
  \caption[Adjoint action $(r)\protect\adjoint(s;Y)$ for $r=1$ and
  $s=2$]{\small An example of the left adjoint action.  A single
    ``new'' cross populates a fixed line with two crosses as is
    described by the braid diagrams.  The cross arrives to each of the
    three possible positions in two ways (one with the plus and the
    other with the minus sign in front).  That the cross never stays
    to the right of the puncture is a manifestation of the fact that
    the space of right coinvariants in a Hopf bimodule is invariant
    under the left adjoint action.}
  \label{fig:adj-braids}
\end{figure}

From Hopf bimodules, we pass to the equivalent
category~\cite{[Besp-Dr-(Bi)]} of Yetter--Drinfeld modules.  The
(left--left) Yetter--Drinfeld modules are given by right coinvariants
in Hopf bimodules, which in our case means taking fixed lines with no
crosses to the right of the puncture,
\begin{equation}\label{2-module}
  \xymatrix@R=6pt@C=30pt{
      \ar@*{[|(1.6)]}@{-}|(.3){\cross}|(.6){\cross}|(.85){\punct}[0,2]&&
    }
\end{equation}
They carry the left ``adjoint'' action of the braided Hopf algebra, as
is illustrated in Fig.~\ref{fig:adj-braids}.

We next allow certain degenerations by considering
$\Nich(X)$-(co)modules that are (subspaces in) $\bigoplus_{i,j\geq
  0}X^{\otimes i}\tensor Y\tensor X^{\otimes j}\tensor Y$, with two
(or more) punctures on the same line:
\begin{equation*}
  \xymatrix@R=6pt@C=50pt{
      \ar@*{[|(1.6)]}@{-}|(.15){\cross}
      |(.3){\punct}|(.45){\cross}|(.65){\cross}|(.90){\punct}[0,2]&&
    }
\end{equation*}
Because the punctures are associated with positions of vertex
operators in CFT, such modules are referred to as two-vertex (or
multivertex) modules.

The multivertex modules carry the following $\cH(X)$ (and $\Nich(X)$)
coaction and action: the deconcatenation coaction \textit{up to the
  first puncture} and the \textit{``cumulative'' left adjoint action},
meaning that all punctures except the rightmost one are viewed on
equal footing with the crosses (an exact definition is given
in~\bref{multi-YD}).  Remarkably, the multivertex modules are
Yetter--Drinfeld modules.  We also evaluate the effect of exchanging
the two punctures in the two-vertex modules, i.e., the
``fusion-product braiding.''  This operation is important, in
particular, because its square commutes with ``everything'' and is
therefore related to the center of the category.\footnote{In CFT{}
  applications, the center is generally expected to be isomorphic to
  the space of torus amplitudes in the corresponding logarithmic CFT{}
  (isomorphic as both commutative associative algebras and $SL(2,\oZ)$
  representations; cf.~\cite{[FGST]} and~\cite{[FGST3],[FGSTq]}).}

When it comes to CFT{}-related applications, the braiding $\Psi$
reduces to a diagonal one.\footnote{Diagonal braiding has been the
  subject of considerable activity in the study of Nichols
  algebras~\cite{[Kh],[GH-lyndon],[Heck-Weyl],[Heck-class],[Heck-1+1],
    [Heck-2],[AA-gen],[Ag-all]}.}  It is convenient to first
specialize to the case where $\Psi$ is the braiding in some $\HHyd$,
the category of Yetter--Drinfeld modules over an ordinary Hopf
algebra~$H$.  A braided Hopf algebra whose braiding is rigid can be
realized as a Hopf algebra in the category of Yetter--Drinfeld modules
over an ordinary Hopf algebra $H$~\cite{[T-survey]},
although $H$ is by far not unique (but the Nichols algebra depends
only on the braiding, not on the Yetter--Drinfeld structure; Nichols
algebras have indeed been extensively studied in the Yetter--Drinfeld
setting).

In actual CFT{} applications in this paper, we restrict ourself to the
simplest case of a one-dimensional space $X$, with the braiding whose
associated total symmetrizer vanishes on $p$-fold tensor products.
This is the setting of the celebrated $(p,1)$ logarithmic CFT{}
models~\cite{[Kausch],[Gaberdiel-K],
  [FHST],[FGST],[FGST2],[S-q],[AM-3],[NT]}, which, their already long
history notwithstanding, are currently being studied from various
standpoints.  (In addition to the works just cited, we also refer the
reader to~\cite{[G-alg],[CF],[GR1],[RS],[Adamovic-M-12],[GR2],
  [GRW],[W],[AM-Z],[HLZ]} and the references therein for the
various aspects of logarithmic conformal field theories.)  We show how
the general constructions defined in terms of ``quantum shuffles''
specialize to the particular case of diagonal braiding associated with
$(p,1)$ logarithmic CFT models.  The general construction of
Yetter--Drinfeld $\Nich(X)$-modules specializes to the corresponding
Nichols algebra.  For the purposes of comparison with logarithmic
CFT{}, the modules thus constructed are then to be viewed \textit{not}
as objects in the braided category of Yetter--Drinfeld
$\Nich(X)$-modules but as ``just'' $\Nich(X)$ module comodules either
without the braiding at all or with the monodromy (``squared
braiding'') structure (the morphisms are then respectively those of
the module comodule structure or of the module comodule plus monodromy
structure).  Thus modified, the Yetter--Drinfeld $\Nich(X)$ modules
then conjecturally provide an equivalence with representations of the
triplet $W$ algebra in the $(p,1)$ model (also viewed as either an
Abelian category or a category with a monodromy structure; we note
that monodromy is intrinsically defined in~CFT{}).

\subsubsection*{Remarks on the known logarithmic Kazhdan--Lusztig
  duality}
Ordinary (nonbraided) Hopf algebras (``factorizable ribbon, although
not quasitriangular, quantum groups'') have previously been proposed
to capture (to different degree) the important features of
representation categories of vertex-operator algebras in logarithmic
CFT{}
models~\cite{[FGST],[FGST2],[FGST3],[FGSTq],[S-q],[NT],[BFGT],[BGT]}.
We expect the contact with this earlier development to be achieved in
terms of \textit{bosonization}~\cite{[Majid-bos]}, or in fact double
bosonization~\cite{[Majid-double],[Sommerh-deformed]}, of braided Hopf
algebras.  Bosonization, which turns ``braided statistics'' into
``bosonic statistics'' where objects are transposed at no expense, has
been known under the name of Radford's biproduct
since~\cite{[Radford-bos]}, before the actual invention of braided
Hopf algebras~\cite{[Mj-braided]}.  From a Hopf algebra
$\cH\in{}\HHyd$, the Radford biproduct produces an ordinary Hopf
algebra based on the smash product $\cH\Smash H$.

The quantum groups featuring in the known instances of logarithmic
Kazhdan--Lusztig duality have been arrived at using a Drinfeld double
of an ordinary Hopf algebra generated by screening(s) and a grading
operator.  The procedure, despite its reasonable success in capturing
the properties of relevant CFT{} models, was nevertheless somewhat
ad~hoc.  It would be quite interesting to see how the previous
proposal can be derived from our braided setting in this paper by a
version of double bosonization for the particular, rather special,
quantum groups.  It is worth noting here that for an ordinary Hopf
algebra $H$, the Drinfeld double $\DD(H)$ answers the question of
describing $H$-module comodules that satisfy the Yetter--Drinfeld
axiom as objects in a monoidal category of \textit{modules} over
\textsl{something}.  The answer is that $\textsl{something}=\DD(H)$
(the category of $\DD(H)$ representation is in fact equivalent to the
category of Yetter--Drinfeld $H$-modules~\cite{[Wor],[Sch-H-YD]}).  In
the braided case, an attempted Drinfeld double construction ``tangles
up''\label{doubleproblem} and does not work unless major simplifying
assumptions are made.\footnote{A ``nonobvious'' construction of a
  ``braided Drinfeld double'' was given in~\cite{[BV]}.}  In the
general case, one is left with just Yetter--Drinfeld modules, which
take over the role of representations of the quantum groups occurring
in the ``nonbraided'' logarithmic Kazhdan--Lusztig duality.

\subsubsection*{Remark on diagrams} 
There are \textit{two types of diagrams} in this paper: (i)~the
standard diagrams used in the theory of braided
(Hopf\;$|$bi$|$\dots\!)algebras (see, e.g., \cite{[Majid-book]}), and
(ii)~braid diagrams (i.e., those \textit{representing elements of the
  braid group algebra}).  These are similar, but not identical in
meaning.  The similarity is that
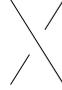
\begin{equation*}
  \begin{tangles}{l}
    \vstr{160}\x\\
  \end{tangles}
\end{equation*}

\noindent
represents an instance of braiding in either case. In the
braided-Hopf-algebra language, the lines can represent two braided
linear spaces, whereas in terms of the braid group algebra, these are
just the two strands of the simplest nontrivial braid.  The difference
is that, for example, the coproduct on a braided Hopf algebra $\cH$ is
represented in braided-algebra language as
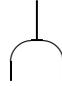
\begin{equation*}
  \begin{tanglec}
    \vstr{150}\cd
  \end{tanglec}
\end{equation*}

\smallskip

\noindent
In what follows, we deal with a \textit{specific} $\cH$, the one
constructed on the tensor algebra of some braided linear space
$X$. The braid diagrams are then those where each strand corresponds a
copy of this $X$.  In terms of braid-group diagrams, the above
coaction diagram translates into a \textit{sum} of diagrams of the
form
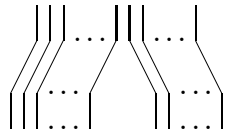
\begin{equation*}
  \begin{tangles}{l}
    \vstr{67}\step[2]\id\step[.5]\id\step[.5]\id\step
    \object{\smash{\dots}}\step\id
    \step[.5]\id\step[.5]\id\step\object{\smash{\dots}}\step\id\\[-4pt]
    \step[1]\dd\step[-.5]\dd\step[-.5]\dd\object{\smash{\dots}}\step\dd
    \step[.5]\d\step[-.5]\d\step\object{\smash{\dots}}\d\\
    \vstr{67}\step\id\step[.5]\id\step[.5]\id\step
    \object{\smash{\dots}}\step\id
    \step[2]
    \step[.5]\id\step[.5]\id\step\object{\smash{\dots}}\step\id
  \end{tangles}
\end{equation*}

\medskip

\noindent
(
that the two bunches are assigned to two different copies of $\cH$ is,
strictly speaking, a ``decoration'' of this trivial braid diagram).

A characteristic example of braided-algebra diagrams is provided
by~\eqref{yd-axiom} and~\eqref{adja} (page~\pageref{yd-axiom}), and
that of braid diagrams, by Fig.~\ref{fig:adj-braids}
(page~\pageref{fig:adj-braids}).  That figure in fact shows
how~\eqref{adja} specializes to the braid group algebra.

\section{Screening configurations as a
  stratification}\label{sec:start-start}
Let $\oC_{w_1,\dots,w_m}$ be the complex plane punctured at $m$ points
$w_1$, \dots, $w_m$.  We consider the configuration space of (an
arbitrary number of) indistinguishable, pairwise noncoincident points
(``screenings'') on $\oC_{w_1,\dots,w_m}$.  The space of such
configurations is
\begin{gather*}
 \cK(w_1,\dots,w_m)=\bigsqcup_{n>0}  \cK_n(w_1,\dots,w_m),\\
 \cK_n(w_1,\dots,w_m)=
 \Bigl(\underbrace{\oC_{w_1,\dots,w_m}\times
   \oC_{w_1,\dots,w_m}
   \times\dots\times\oC_{w_1,\dots,w_m}}_n\setminus\,\Delta\Bigr)
 \Bigr/\mathbb{S}_n,
\end{gather*}
where $\setminus\,\Delta$ means the removal of all subdiagonals (the
points are noncoincident), and $\bigl/\,\mathbb{S}_n$ is the quotient
by the action of the symmetric group on $n$ elements (the points are
indistinguishable).  We call the points ``crosses'' in accordance with
the chosen way of their pictorial representation.  In contrast to the
punctures, which are fixed, the crosses are \textit{movable} because
neither their positions on any line, nor the positions of the movable
lines are fixed.  

We fix two braided linear spaces $X$ and $Y$ and associate a copy of
$Y$ with each puncture and a copy of $X$ with each cross.

The aim of this section is to construct a local system
$\lS_m=\lS(w_1,\dots,w_m)$ whose associated flat bundle, with the
fiber linearly isomorphic to $X^{\otimes n}\otimes Y^{\otimes m}$,
generalizes the bundle of conformal blocks, with a
Knizhnik--Zamolodchikov-like flat connection.  We begin with
introducing a stratification of $\cK_n(w_1,\dots,w_m)$ based on a
``space--time'' decomposition of the complex plane.  We write
$\cK_{m,n}=\cK_n(w_1,\dots,w_m)$ for brevity.

\subsection{Local system}\label{sec:LS}
For a fixed $m$, we represent each $\cK_{m,n}$ as
\begin{equation*}
  \cK_{m,n} = \bF^{2n}
  \supset\bF^{2n-1}
  \supset\dots
  \supset\bF^{n},
\end{equation*}
where the $\bF^k$ are dimension-$k$ strata, $\bF^k$ being closed in
$\bF^{k+1}$, each consisting of several connected components (also
called strata):
\begin{equation*}
  \bF^k=\bigcup_{j\in I_k} \bF^k_j,\quad k=n,\dots, 2n,
\end{equation*}
where $I_k$ are some finite sets (for example,
$|I_{2n}|=\binom{n+m}{n}$).

We actually construct \textit{open} contractible subspaces
$F^{k}_{j_k}$ whose closure in $F^{k+1}$ is the corresponding
connected stratum, $\overline{F^{k}_{j}}=\bF^{k}_{j}$, and the
boundary of each $F^{k}_{j}$ is contained in 
$\bigcup_{i}F^{k-1}_{i}$.
Then
\begin{equation*}
  \cK_{m,n} = \bigcup_{n\leq k\leq 2n} \bigcup_{j_k\in I_k}F^{k}_{j_k},
\end{equation*}
We call the $F^{k}_{j_k}$ ``cells.''
 
Construction of the $F^{k}_{j}$ is described in~\bref{cells:1}.  For a
given $F^{k}_{a}$, the $F^{k-1}_{b}$ such that the codimension-$1$
part of $\partial F^{k}_{a}$ is contained in $\bigcup_{b\in I_{k-1}}
F^{k-1}_{b}$ are described in~\bref{sec:incidence}.

\subsubsection{Constructing the cells}\label{cells:1}
All cells $F^{k}_{j}$ are obtained by simply placing $n$ crosses on a
family of parallel lines (``spatial sections'') and then
hierarchically restricting the possible configurations of the lines
(see Fig.~\ref{fig:intro}).
\begin{enumerate}\addtocounter{enumi}{-1}
\item Select, once and for all, a straight line on
  $\oC_{w_1,\dots,w_m}$ and consider the family of parallel lines,
  assuming that by the choice of the first line, none of the lines
  passes through more than one puncture.  An orientation must also be
  selected, the same on all lines.  In addition, the lines must be
  ordered inside the family (``the global time'').  \footnote{A
    possible way to do this is to take the lines to be level lines of
    the function $f(z)=\alpha\,\mathrm{Re}(z) + \beta\,\mathrm{Im}(z)$
    with real $\alpha$ and $\beta$; some other choices of the function
    suggest interesting generalizations.}

\item Draw a line from the family passing through each puncture and
  label these $m$ lines ``monotonically'' in accordance with the order
  chosen in the family.  These $m$ lines are called \textit{fixed
    lines} in what follows.

\item Draw $r$, $0\leq r\leq n$, lines from the family not passing
  through any puncture; these lines are not fixed but can move so as
  to remain parallel and to not collide with or pass through one
  another or any of the fixed lines.  Place at least one cross on each
  of these movable lines.
\end{enumerate}

A cell is by definition a piece of $\cK_{m,n}$ consisting of all
configurations where the crosses are placed on fixed and movable lines
and are allowed to slide along the lines so as to \textit{not collide
  with one another and not pass through the punctures}, and the
movable lines move as just described.  The (real) dimension of a cell
constructed this way is $k = n + r$.

\medskip

\noindent
\textit{\thesubsubsection.1.\ \ Example.}  \ Two examples of cells in
$\cK_{2,7}$ are represented by the configurations
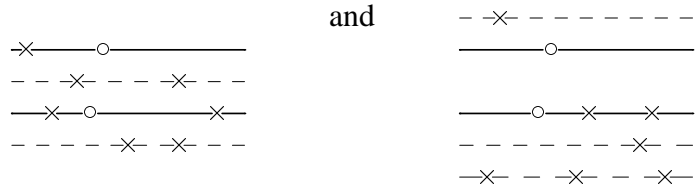
\begin{equation*}
 \xymatrix@R=6pt@C=40pt@R=4pt{
   \mbox{}&\\
   \ar@*{[|(1.6)]}@{-}|(.1){\cross}|(.4){\punct}[0,2]&&\\
   \ar@{--}|(.3){\cross}|(.7){\cross}[0,2]&&\\
   \ar@*{[|(1.6)]}@{-}|(.2){\cross}|(.35){\punct}|(.85){\cross}[0,2]&&\\
   \ar@{--}|{\cross}|(.7){\cross}[0,2]&&
 }
 \qquad\text{and}\qquad
 \xymatrix@R=6pt@C=40pt@R=4pt{
   \ar@{--}|(.2){\cross}[0,2]&&\\
   \ar@*{[|(1.6)]}@{-}|(.4){\punct}[0,2]&&\\
   \\
   \ar@*{[|(1.6)]}@{-}|(.35){\punct}|(.55){\cross}|(.8){\cross}[0,2]&&\\
   \ar@{--}|(.75){\cross}[0,2]&&\\
   \ar@{--}|(.15){\cross}|(.5){\cross}|(.85){\cross}[0,2]&&
 }
\end{equation*}

\medskip

\noindent
\textit{\thesubsubsection.2.\ \ Example} \ On a \textit{nonpunctured}
plane, all cells are enumerated following the simple pattern (with
dotted lines separating different cells)
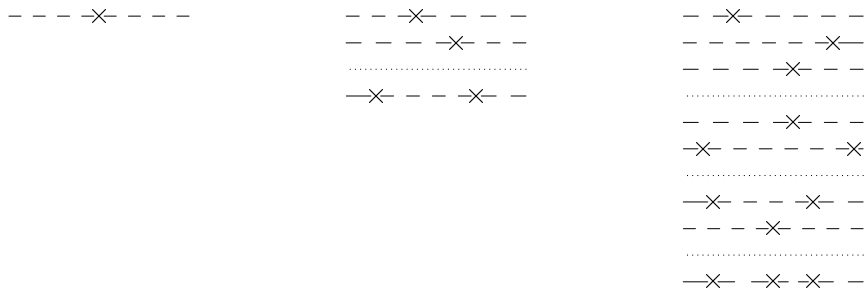
\begin{xalignat*}{3}
 &{\xymatrix@R=6pt@C=30pt@R=2pt{
     \ar@{--}|(.5){\cross}[0,2]&&
   }}
 &&{\xymatrix@R=6pt@C=30pt@R=2pt{
     \ar@{--}|(.4){\cross}[0,2]&&\\
     \ar@{--}|(.6){\cross}[0,2]&&\\
     \ar@{.}[0,2]&&\\
     \ar@{--}|(.2){\cross}|(.7){\cross}[0,2]&&\\
   }}
 &&{\xymatrix@R=6pt@C=30pt@R=2pt{
     \ar@{--}|(.3){\cross}[0,2]&&\\
     \ar@{--}|(.8){\cross}[0,2]&&\\
     \ar@{--}|(.6){\cross}[0,2]&&\\
     \ar@{.}[0,2]&&\\
     \ar@{--}|(.6){\cross}[0,2]&&\\
     \ar@{--}|(.15){\cross}|(.9){\cross}[0,2]&&\\
     \ar@{.}[0,2]&&\\
     \ar@{--}|(.2){\cross}|(.7){\cross}[0,2]&&\\
     \ar@{--}|(.5){\cross}[0,2]&&\\
     \ar@{.}[0,2]&&\\
     \ar@{--}|(.2){\cross}|(.5){\cross}|(.7){\cross}[0,2]&&
   }}
\end{xalignat*}

\noindent
with the columns respectively corresponding to $\cK_{0,1}$,
$\cK_{0,2}$, $\cK_{0,3}$, and so on ($2^{n-1}$ cells in the $n$th
column).

\medskip

\noindent
\textit{\thesubsubsection.3. \ \ Example} \ There are $\binom{n+m}{n}$
cells of the maximum dimension $2n$.  To be more specific, we can
assume that the punctures $w_1,\dots,w_m$ are chosen such that
$a_1\equiv\mathrm{Im}\,w_1 < a_2\equiv\mathrm{Im}\,w_2< \dots
<a_m\equiv\mathrm{Im}\,w_m$.  In $\mathbb{R}^{2n}$ with coordinates
$x_i, y_i$, $i=1,\dots,n$, the list of maximum-dimension cells is then
given by
\begin{gather*}
  -\infty<y_1<y_2<\dots<y_n<a_1<\dots<a_m<+\infty,\\
  -\infty<y_1<y_2<\dots<y_{n-1}<a_1<y_n<\dots<a_m<+\infty,\\
  -\infty<y_1<y_2<\dots<y_{n-1}<a_1<a_2<y_n<\dots<a_m<+\infty,\\
  \vdots\\
  -\infty<y_1<y_2<\dots<a_1<y_{n-1}<a_2<\dots<a_{m-1}<y_n<a_m<+\infty,\\
  \vdots\\
  -\infty<a_1<\dots<a_m<y_1<y_2<\dots<y_n<+\infty
\end{gather*}
with $-\infty<x_i<+\infty$, $1\leq i\leq n$, in each row (each cell).
This explicitly shows that the cells are open neighborhoods in
$\mathbb{R}^{2n}$.

\medskip

It follows that the dimension-$k$ cells, $n\leq k\leq 2n$, are in
bijective correspondence with the data
\begin{equation}\label{cell-def}
  ([j_1],[j_2],\dots,[j_{k+m-n}]),
\end{equation}
where
\begin{enumerate}
\item each $[j_a]$ is either an integer $j_a>0$ or a pair of integers
 $(j_a'\geq0,j_a''\geq0)$,

\item there are exactly $m$ pairs of integers among the $[j_a]$,

\item $\sum_{a=1}^{k+m-n}[j_a]=n$, with each summand understood as
  either $j_a$ or $j_a'+j_a''$ depending on the type.
\end{enumerate}
In~\eqref{cell-def}, each $[j_a]=j_a$ represents a movable line
carrying $j_a$ crosses, and each $[j_a]=(j_a',j_a'')$ represents a
fixed line carrying $j_a'$ crosses before and $j_a''$ crosses after
the puncture (relative to the chosen orientation).  The \textit{total}
number of lines is $k+m-n$, of which $k-n$ are movable lines.

For fixed punctures, any cell represented as in~\eqref{cell-def} is
uniquely assigned to one of the $\cK_{m,n}$ spaces: the number of
punctures $m$ is just the number of pairs of integers among the
$[j_a]$, and the number of crosses is $n=\sum_{a=1}^{k+m-n}[j_a]$, as
noted above.  This allows us to omit explicit indications of the
$\confS_{m,n}$ space to which a given cell belongs.

The two cells in~\textbf{\thesubsubsection.1} are thus written as
$((1,0),2,(1,1),2)$ and $(1,(0,0),(0,2),1,3)$.

\textit{Oriented} cells are $([j_1],[j_2],\dots,[j_s];\pm1)$.

\subsubsection{The $F^{k-1}_b$ in the closure of $F^k_a$: line
  merger}\label{sec:incidence}
We now (indirectly) describe all pairs $(F^k_a, F^{k-1}_b)$ of a
dimension-$k$ and a dimension-$(k-1)$ cell such that the
codimension-$1$ part of $\partial F^k_a$ contains $F^{k-1}_b$.  This
amounts to describing all cases of a merger of either two neighboring
movable lines or a movable line with a neighboring fixed line.

For each cell $([j_1],[j_2],\dots,[j_s])$ written as in
\eqref{cell-def}, the procedure is as follows.  Find all pairs
$[j_r],[j_{r+1}]$ such that at least one of $[j_r]$ and $[j_{r+1}]$ is
an integer (not a pair of integers).
\begin{enumerate}
\item If both are integers, then replace the chosen pair with the
  integer $j_r+j_{r+1}$ (merger of two movable lines).  This produces
  the new cell
  \begin{equation*}
    ([j_1],\dots,[j_{r-1}],j_r+j_{r+1},[j_{r+2}],\dots,[j_s]). 
  \end{equation*}
  
\item If $[j_{r+1}]=(j_{r+1}',j_{r+1}'')$ (merger of a movable and a
  fixed line), then $j_r + 1$ new cells are produced.  For each $0\leq
  b\leq j_r$, replace $(j_{r+1}',j_{r+1}'')$ with $(b + j_{r+1}',
  j_{r+1}'' + j_r - b)$, $0\leq b\leq j_r$, which gives the cell
 \begin{equation*}
   ([j_1],\dots,[j_{r-1}],(b + j_{r+1}', j_{r+1}'' + j_r  - b),
   [j_{r+2}],\dots,[j_s]).
 \end{equation*}
 If $[j_{r}]=(j_{r}',j_{r}'')$, then do the same with $r$ and $r+1$
 swapped.
\end{enumerate}
The cell dimension decreases by unity in either case, simply because
fewer movable lines are left for the crosses to slide along.  The new
cell is assigned the same orientation as $F^k_a$ if $r$ is odd, and
the reversed orientation if $r$ is even.

\medskip

\noindent
\textit{\thesubsubsection.1.  Example.} \ \ The codimension-$1$ part
of the boundary of the dimension-$6$ cell
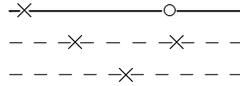
\begin{equation*}
  \xymatrix@R=6pt@C=40pt@R=4pt{
   \ar@*{[|(1.6)]}@{-}|(.1){\cross}|(.67){\punct}[0,2]&&\\
   \ar@{--}|(.3){\cross}|(.7){\cross}[0,2]&&\\
   \ar@{--}|{\cross}[0,2]&&
 }
\end{equation*}
is the union of the cells
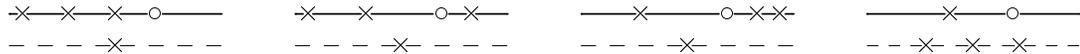
\begin{equation*}
  \xymatrix@R=6pt@C=36pt@R=4pt{
   \ar@*{[|(1.6)]}@{-}|(.1){\cross}|(.3){\cross}|(.5){\cross}
   |(.67){\punct}[0,2]&&\\
   \ar@{--}|{\cross}[0,2]&&
 }\quad
 \xymatrix@R=6pt@C=36pt@R=4pt{
   \ar@*{[|(1.6)]}@{-}|(.1){\cross}|(.35){\cross}
   |(.67){\punct}|(.80){\cross}[0,2]&&\\
   \ar@{--}|{\cross}[0,2]&&
 }\quad
 \xymatrix@R=6pt@C=36pt@R=4pt{
   \ar@*{[|(1.6)]}@{-}|(.3){\cross}
   |(.67){\punct}|(.80){\cross}|(.9){\cross}[0,2]&&\\
   \ar@{--}|{\cross}[0,2]&&
 }\quad
  \xymatrix@R=6pt@C=36pt@R=4pt{
   \ar@*{[|(1.6)]}@{-}|(.4){\cross}|(.67){\punct}[0,2]&&\\
   \ar@{--}|(.3){\cross}|(.5){\cross}|(.7){\cross}[0,2]&&
 }
\end{equation*}

\medskip

\noindent
\textit{\thesubsubsection.2.  Example.} \ \ For each cell in the list
in~\bref{cells:1}\textbf{.3}, its boundary is obtained when $y_i$ and
$y_{i+1}$ coincide or when a $y_i$ coincides with a neighboring $a_j$.
For example, for the first cell in the list, taking $y_1$ equal to
$y_2$ gives the $(2n-1)$-dimensional cell
\begin{gather*}
  \left\{
    \begin{aligned}
      &-\infty<y_1<y_3\dots<y_{n-1}<y_n<a_1< a_2<\dots<a_n<+\infty,\\
      &-\infty<x_1<x_2<+\infty,\\
      &-\infty<x_i<+\infty,\quad i\geq3
    \end{aligned}
  \right.
\end{gather*}
on the boundary.

\subsubsection{Local system: the restriction morphisms}
We choose two braided linear spaces $(X,\Psi)$ and $(Y,\Psi)$.  With
each cell \eqref{cell-def}, we associate the braided linear space
\begin{equation}\label{sections}
  \Gamma\bigl(([j_1],[j_2],\dots,[j_s])\bigr)=
  (X^{\tensor[j_1]})\tensor(X^{\tensor[j_2]})\tensor\dots\tensor(X^{\tensor[j_s]}),
\end{equation}
where $X^{\tensor[j]}$ is either $X^{\tensor j}$ or $X^{\tensor
  j'}\tensor Y\tensor X^{\tensor j''}$ depending on whether $[j]$ is
an integer or a pair of integers.  Although all spaces
in~\eqref{sections} are linearly isomorphic to $X^{\otimes n}\otimes
Y^{\otimes m}$, \textit{the brackets must not be removed} in this type
of expressions; the bracket positions actually characterize the
underlying cell.


To define the local system $\lS_{m,n}$, it remains to specify the
morphisms between vector spaces of form~\eqref{sections} for the
mergers described in~\bref{sec:incidence}\,---\,the restriction maps
$\Gamma(F^k_a)\to\Gamma(F^{k-1}_b)$ corresponding to each case where
$F^{k-1}_b$ is in the closure of $F^k_a$.  These maps can be regarded
as a procedure to ``collectivize'' the crosses populating these lines.

\begin{figure}[htb]
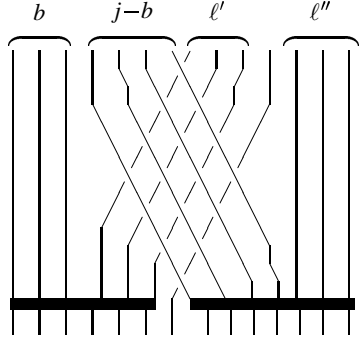

   \centering
   \begin{equation*}
     \begin{tangles}{l}
       \object{\kern20pt\overgroup{\rule{25pt}{0pt}}^{b}}
       \object{\kern90pt\overgroup{\rule{35pt}{0pt}}^{j-b}}
       \object{\kern155pt\overgroup{\rule{25pt}{0pt}}^{\ell'}}
       \object{\kern233pt\overgroup{\rule{30pt}{0pt}}^{\ell''}}
       \\[-3pt]
       \vstr{67}\hstr{67}\id\step[1.5]\id\step[1.5]\id\step[1.5]\id
       \step[1.5]\dh\step[1]\dh\step[1]\hx\step[1]\ddh\step[1]\ddh
       \step[1.5]\id\step[1.5]\id\step[1.5]\id\step[1.5]\id\\
       \vstr{67}\hstr{67}\id\step[1.5]\id\step[1.5]\id\step[1.5]\dh
       \step[1.5]\dh\step[1]\hx\step[1]\hx\step[1]\ddh\step[1.5]\ddh
       \step[1.5]\id\step[1.5]\id\step[1.5]\id\\
       \vstr{67}\hstr{67}\id\step[1.5]\id\step[1.5]\id\step[1.5]\step[.5]
       \d\step[1]\hx\step[1]\hx\step[1]\hx\step[1]\dd\step[2]\id
       \step[1.5]\id\step[1.5]\id\\
       \vstr{67}\hstr{67}\id\step[1.5]\id\step[1.5]\id\step[1.5]
       \step[1.5]\hx\step[1]\hx\step[1]\hx\step[1]\hx\step[3]\id
       \step[1.5]\id\step[1.5]\id\\
       \vstr{67}\hstr{67}\id\step[1.5]\id\step[1.5]\id\step[1.5]
       \step[.5]\dd\step[1]\hx\step[1]\hx\step[1]\hx\step[1]\d\step[2]
       \id\step[1.5]\id\step[1.5]\id\\
       \vstr{67}\hstr{67}\id\step[1.5]\id\step[1.5]\id\step[1.5]
       \step[.5]\id\step[1.5]\hdd\step[1]\hx\step[1]\hx\step[1]\d
       \step[1]\hd\step[1.5]\id\step[1.5]\id\step[1.5]\id\\
       \vstr{67}\hstr{67}\object{\kern53pt\rule[-4pt]{55pt}{4pt}}
       \id\step[1.5]\id\step[1.5]\id\step[1.5]\step[.5]\id\step[1.5]
       \id\step[1.5]\id\step[1]\hx\object{\kern62pt\rule[-4pt]{62pt}{4pt}}
       \vstr{67}\hstr{67}\step[1]\d\step[1]\hd\step[1]\hd\step[1]
       \id\step[1.5]\id\step[1.5]\id\\[-4pt]
       \vstr{67}\hstr{67}\id\step[1.5]\id\step[1.5]\id\step[1.5]
       \id\step[1.5]\id\step[1.5]\id\step[1.5]\id\step[2]\id\step[1.3]
       \id\step[1.3]\id\step[1.3]\id\step[1.3]\id\step[1.3]\id
       \step[1.5]\id
     \end{tangles}
   \end{equation*}
   \caption[Illustration of formula~\eqref{eq:mu}]{\small Illustration
     of formula~\eqref{eq:mu}.  Read right to left (which corresponds
     to reading the diagram from top down), this formula tells the
     following: first, $j-b$ crosses (representing
     $(X^{\tensor(j-b)})$) ``travel past'' $\ell'$ crosses and the
     puncture (representing~$X^{\tensor \ell'}\tensor Y$); next, the
     $j-b$ and $\ell''$ are collectivized by taking a sum over
     $\binom{j-b+\ell''}{j-b}$ shuffle permutations, which yields
     $\Bbin{}{j - b, \ell''}$ (the right dark strip in the picture);
     and finally, the remaining $b$ are similarly collectivized with
     the $\ell'$ crosses, the sum over the $\binom{b+\ell'}{b}$ shuffle
     permutations yielding the braided binomial $\Bbin{}{b, \ell'}$ (the
     left dark strip).}
   \label{fig:mu}
 \end{figure}

\begin{enumerate}
\item In the case of a merger of two movable lines that originally
  carried $\ell$ and $j$ crosses, the crosses can be ``collectivized''
  in $\binom{\ell+j}{\ell}$ different ways, labeled by shuffle
  permutations
  $\sigma\in\mathbb{S}^{\ell+j}_{\ell,j}\subset\mathbb{S}_{\ell+j}$.
  For each such $\sigma$, let $\hat{\sigma}$ be its Matsumoto section
  (see~\bref{sec:q-shuffles}),
  \begin{equation*}
    \hat{\sigma}:\Gamma(F^k_a)\to\Gamma(F^{k-1}_b)
  \end{equation*}
  Here, $\Gamma(F^k_a)$ is a product~\eqref{sections}, necessarily of
  the form $\dots(X^{\tensor\ell}) \tensor(X^{\tensor j})\dots$; the
  map $\hat{\sigma}$ is nontrivial only on these two factors, and we
  do not write the rest of the tensor factors.  Then the restriction
  map $\Gamma(F^k_a)\to\Gamma(F^{k-1}_b)$ is given by the sum of the
  above $\hat{\sigma}$ over all shuffle permutations
  $\sigma\in\mathbb{S}^{\ell+j}_{\ell,j}$.  By definition, this is the
  braided binomial $\Bbin{}{\ell,j}$ (see~\bref{sec:q-shuffles}):
 \begin{equation*}
   \Bbin{}{\ell,j}\equiv
   \sum_{\sigma\in\mathbb{S}^{\ell+j}_{\ell,j}}\!\!\hat{\sigma}
   :
   (X^{\tensor\ell}) \tensor(X^{\tensor j})
   \to(
   X^{\tensor(\ell+j)}).
 \end{equation*}

\item In the case of a merger of a movable line carrying $j$
  crosses with a fixed line carrying $\ell'$ crosses on the left and
  $\ell''$ crosses on the right of the puncture, any number $b$, $0\leq
  b\leq j$, of the $j$ crosses may go to the left of the puncture
  (again, we understand ``left'' and ``right'' in terms of the chosen
  orientation on the lines).  This gives $j+1$ distinct
  $F^{k-1}_b$; for each $b$, the number of possibilities of
  collectivizing the crosses is
  $\binom{b+\ell'}{b}\binom{j-b+\ell''}{j-b}$.  The relevant map
  $\Gamma(F^k_a)\to\Gamma(F^{k-1}_b)$ is different from the identity
  only in the component
 \begin{equation*}
   \mu_{b;j,\ell',\ell''}:(X^{\tensor j}) \tensor(X^{\tensor \ell'}\tensor Y\tensor
   X^{\tensor \ell''})
   \to
   (X^{\tensor(\ell'+b)}\tensor Y\tensor X^{\tensor(\ell''+j-b)}),
 \end{equation*}
 where it is given by (see
 Fig.~\ref{fig:mu})
 \begin{equation}\label{eq:mu}
   \mu_{b;j,\ell',\ell''} =
   \Bbin{}{b, \ell'}
   \ccirc\xShift{b + \ell' + 1}{
     \Bbin{\shift(b + \ell' + 1)}{j - b, \ell''}}
   \ccirc\Shift{b}{\Bbraid_{j - b, \ell' + 1}}.
 \end{equation}
\end{enumerate}

\subsubsection{}
A useful \textit{algebraic} simplification in the second case of line
merger can be achieved by considering the sum of the
$\mu_{b;j,\ell',\ell''}$ maps over all $0\leq b\leq j$.  The braided
binomial identity\footnote{Which is a ``quantum'' counterpart of the
  identity
  \begin{equation*}
    \sum_{b=0}^{j}
    \mbox{\footnotesize$\displaystyle
      \binom{b+\ell'}{b}\binom{j-b+\ell''}{j-b}$}
    =\mbox{\footnotesize$\displaystyle
      \binom{j + \ell' + \ell'' + 1}{j}$},
  \end{equation*}
  where the left-hand side is the total number of different
  possibilities of collectivizing the crosses for all $0\leq b\leq j$,
  and the right-hand side is the number of shuffle permutations that
  ``mix'' $j$ crosses with $\ell' + 1 + \ell''$
  ``points.''}
\begin{equation*}
  \sum_{b=0}^{j}
  \mu_{b;j,\ell',\ell''}
  =\Bbin{}{j, \ell' + 1 + \ell''}
\end{equation*}
then shows that the $j+1$ cases of different collectivization
types, taken together, are conveniently described by the map
$\Gamma(F^{k}_a)\to\bigoplus_{b=0}^{j}\Gamma(F^{k-1}_b)$ whose
nontrivial component~is
\begin{equation*}
 \Bbin{}{j, \ell' + \ell'' + 1}:
 (X^{\tensor j}) \tensor(X^{\tensor \ell'}\tensor Y\tensor X^{\tensor \ell''})
 \to
 \bigoplus_{b=1}^{j}(X^{\tensor(\ell'+b)}\tensor Y\tensor
 X^{\tensor(\ell''+j-b)}).
\end{equation*}

\subsection{Homology complex}\label{sec:homology}
The above construction is neatly summarized if our local system
$\lS_m$ is used, indeed, as a local system \textit{of coefficients}.
The homology complex of $\cK_{m,n}$ with coefficients in $\lS_m$ is
\begin{equation}\label{the-complex}
  0\xleftarrow{}\Cepx_{m,n}\xleftarrow{\partial}\Cepx_{m,n+1}
  \xleftarrow{\partial}\dots\xleftarrow{\partial}\Cepx_{m,2n}
  \xleftarrow{}0,
\end{equation}
where
\begin{equation*}
  \Cepx_k=\bigoplus_{([j_1],[j_2],\dots,[j_{k-n}])}
  (X^{\tensor[j_1]})\tensor(X^{\tensor[j_2]})
  \tensor\dots\tensor
  (X^{\tensor[j_{k-n}]})
\end{equation*}
with the sum taken over all sets $([j_1],[j_2],\dots,[j_{k-n}])$ such
that $\sum_{r=1}^{k-n}[j_r]=n$ and each set contains exactly $m$
\textit{pairs} $(j'_i,j''_i)$.
The differential acts on each space~\eqref{sections} by the braid
group algebra element
\begin{equation}\label{hom-diff}
 \partial=\sum_{i=2}^{s}(-1)^{i}
 \xShift{\#[j_1]+\dots+\#[j_{i-2}]}{
   \Bbin{\shift(\#[j_1]+\dots+\#[j_{i-2}])}{[j_{i-1}],[j_{i}]}}
\end{equation}
where we define the ``length''
\begin{equation}
 \#[j]=
 \begin{cases}
   j,& [j]=j,\\
   j'+j''+1,& [j]=(j',j'')
 \end{cases}
\end{equation}
and set
\begin{equation}
 \Bbin{}{[j],[\ell]}=
 \begin{cases}
   0,& [j]=(j',j'') \text{ and } [\ell]=(\ell',\ell''),\\
   \Bbin{}{\#[j],\#[\ell]},&\text{otherwise,}
 \end{cases}
\end{equation}
which map between the relevant spaces as follows:
\begin{equation*}
 (X^{\tensor[j]})\tensor(X^{\tensor[\ell]})\to
 \begin{cases}
   X^{\tensor(j+\ell)},& [j]=j \text{ and } [\ell]=\ell,\\
   \bigoplus\limits_{a=0}^j X^{\tensor(\ell'+a)}\tensor Y\tensor
   X^{\tensor(j-a+\ell'')},&
   [j]=j \text{ and } [\ell]=(\ell',\ell''),\\
   \bigoplus\limits_{a=0}^{\ell} X^{\tensor(j'+a)}\tensor Y\tensor
   X^{\tensor(\ell-a+j'')},&
   [j]=(j',j'') \text{ and } [\ell]=\ell.
 \end{cases}
\end{equation*}

\subsubsection{Example}
We write explicitly how the differential acts on the space
\begin{equation*}
  (Y\tensor X)\tensor(X^2)\tensor(X\tensor Y)=\Gamma(F^{5}),
\end{equation*}
whose underlying cell is $F^{5}=((0,1),2,(1,0))$.  In codimension $1$,
the cell boundary is the union of cells
$\bigl((0,1),\,(3,0)\bigr)\cup\bigl((0,3),\,(1,0)\bigr)
\cup\bigl((0,1),\,(2,1)\bigr)\cup\bigl((1,2),\,(1,0)\bigr)
\cup{}\bigl((0,1),\,(1,2)\bigr) \cup\bigl((2,1),\,(1,0)\bigr)$.  Under
each term in the differential, we indicate the space \textit{to which}
it maps:
\begin{multline*}
 \partial\bigm|_{(Y\tensor X)\tensor(X^2)\tensor(X\tensor Y)}=
 -\underset{(Y\tensor X)\tensor(X^3\tensor Y)
 }{\id}+
 \underset{(Y\tensor X^3)\tensor(X\tensor Y)}{\id}
 -\underset{(Y\tensor X)\tensor(X^3\tensor Y)}{ \Psi_{4}}
 +
 \underset{(Y\tensor X^3)\tensor(X\tensor Y)}{ \Psi_{2}}
 \\
 -\underset{(Y\tensor X)\tensor(X^3\tensor Y)}{ \Psi_{3} \Psi_
   {4}}
 +\underset{(Y\tensor X^3)\tensor(X\tensor Y)}{ \Psi_{3} \Psi_
   {2}}
 -\underset{(Y\tensor X)\tensor(X^2\tensor Y\tensor X)}{\Psi_{5}\Psi_{4}}
 +\underset{(X\tensor Y\tensor X^2)\tensor(X\tensor Y)}{ \Psi_{1}\Psi_{2}}
 \\
 -\underset{(Y\tensor X)\tensor(X^2\tensor Y\tensor X)}{
   \Psi_{3}\Psi_{5} \Psi_{4}}
 +\underset{(X\tensor Y\tensor X^2)\tensor(X\tensor Y)}{
   \Psi_{1} \Psi_{3} \Psi_{2}}
 -\underset{(Y\tensor X)\tensor(X\tensor Y\tensor X^2)}{
   \Psi_{4} \Psi_{3} \Psi_{5} \Psi_{4}}
 +\underset{(X^2\tensor Y\tensor X)\tensor(X\tensor Y)}{
   \Psi_{2} \Psi_{1} \Psi_{3} \Psi_ {2}}
\end{multline*}
(the two identity operators are different operators because they map
to different spaces).

\subsection{Flat bundle}
The local system $\lS_{m,n}$ can be alternatively described in terms
of horizontal sections of a flat bundle over $\cK_{m,n}$.  The total
space of the bundle is the quotient~of
\begin{equation*}
  \bigl(
  (\oC_{w_1,\dots,w_m})^{\times n}\,
  \setminus\,\Delta\bigr)
  \times \bigl(X^{\otimes n}\otimes Y^{\otimes m}\bigr)
\end{equation*}
by the identifications
\begin{equation*}
  s\bigl(
  (\oC_{w_1,\dots,w_m})^{\times n}\,
  \setminus\,\Delta\bigr)
  \times \bigl(X^{\otimes n}\otimes Y^{\otimes m}\bigr)
  \sim \bigl(
  (\oC_{w_1,\dots,w_m})^{\times n}\,
  \setminus\,\Delta\bigr)
  \times\hat{s}\bigl(X^{\otimes n}\otimes Y^{\otimes m}\bigr),
\end{equation*}
where $\hat{s}$ is the Matsumoto section of $s\in\mathbb{S}_n$, which
is a linear operator defined on $X^{\otimes n}\otimes Y^{\otimes m}$.

For the nonpunctured complex plane, the flat bundle is of course
equivalently specified by a morphism $\pi_1(\cK_{0,n})\to G$, where
$G$ is the structure group, $G\subset\End(X^{\otimes n})$, and
$\pi_1(\cK_{0,n})=\mathbb{B}_n$\,---\,that is, by a braid group
representation (the one we started with, of course).

\subsection{Algebraic structures associated with the
  chains}\label{sec:chain-struct}
The braided linear space
\begin{equation*}
  \cH_+ = \xymatrix@1{
    \ar@{--}|(.5){\cross}[0,2]&&}\;\oplus\;
  \xymatrix@1{
    \ar@{--}|(.3){\cross}|(.7){\cross}[0,2]&&}\;\oplus\;
  \xymatrix@1{\ar@{--}|(.25){\cross}|(.5){\cross}|(.75){\cross}[0,2]&&}
  \;\oplus\;\dots
  =\bigoplus_{j\geq 1} (X^{\tensor j}),
\end{equation*}
plays a special role in what follows.  An associative product can be
constructed on $\cH_+$ as a composition of the tensor product and the
differential, i.e., $x\dotact x'=\partial(x\tensor x')$.  For $x\in
(X^{\tensor j})$ and $x'\in (X^{\tensor j'})$, this gives
\begin{equation}\label{H-mult-d}
  x\dotact x'=\partial(x\tensor x')
  =\Bbin{}{j,j'}(x\tensor x')\in (X^{\tensor (j+j')}).
\end{equation}

By the same pattern, an $x\in (X^{\tensor j})$ \textit{acts} on a
$y\in X^{\tensor (j',j'')}=X^{\tensor j'}\tensor Y\tensor X^{\tensor
  j''}$ from the left and from the right as
\begin{equation}\label{x-x-act}
  x\dotact y=\partial(x\tensor y)
  =\Bbin{}{j,j'+j''+1}(x\tensor y),\quad
  y\dotact x=\partial(y\tensor x)
  =\Bbin{}{j'+j''+1,j}(y\tensor x).
\end{equation}

\subsubsection{Augmentation of $\cH_+$}\label{sec:unit}
We next let $\cH$ be the unital algebra obtained by augmenting
$(\cH_+,{}\dotact{})$ with a formal chain~$1$---the tensor unit in the
monoidal category of braided linear spaces, which can be regarded as
an object of type \eqref{sections} with zero number of $X$ or $Y$,
i.e., as~$(\,)=(X^{\tensor0})$.  Then
\begin{equation}\label{H-def}
  \cH=\cH(X)=\bigoplus_{j\geq0} (X^{\tensor j})
\end{equation}
The action of differential \eqref{hom-diff} is immediately extended to
any chain with any number of $(\,)$, by setting
$\Bbin{}{0,n}=\Bbin{}{n,0}=\id$.  Hence, $(\,)\dotact
x=x\dotact(\,)=x$ for any chain $x$.

\subsubsection{Coproduct}\label{sec:comult}
The algebra $\cH(X)$ can be made into a bialgebra by line cutting (or
splitting, depending on the interpretation).  For any line with $n$
crosses, the coproduct is the result of cutting the line into two and
summing over the possibilities of how many crosses occur to the left
and to the right of the cut,
\begin{align}\label{cutting}
  \Delta:\xymatrix@R=4pt@C=28pt{
    \ar@{--}|(.25){\cross}|(.5){\cross}|(.75){\cross}[0,2]&&
  }\mapsto{}& \xymatrix@R=4pt@C=10pt{
    \ar@{--}|(.25){\cross}|(.5){\cross}|(.75){\cross}[0,3]
    &&&{|}\ar@{--}[0,3]&&& } + \xymatrix@R=4pt@C=10pt{
    \ar@{--}|(.33){\cross}|(.66){\cross}[0,3]
    &&&{|}\ar@{--}|(.5){\cross}[0,3]&&& }
  \\
  \notag &+ \xymatrix@R=4pt@C=10pt{ \ar@{--}|(.5){\cross}[0,3]
    &&&{|}\ar@{--}|(.33){\cross}|(.66){\cross}[0,3]&&& } +
  \xymatrix@R=4pt@C=10pt{ \ar@{--}[0,3]
    &&&{|}\ar@{--}|(.25){\cross}|(.5){\cross}|(.75){\cross}[0,3] &&& }
  \\
  \intertext{with the line cutting symbol $|$ then to be read
    as~$\tensor$,}\notag
  ={}&(X^{\tensor3})\tensor()\oplus(X^{\tensor2})\tensor(X)
  \oplus(X)\tensor(X^{\tensor2})\oplus ()\tensor(X^{\tensor3})
\end{align}
(Figuratively speaking, ``we do not know'' where the cut happens to
break the crosses---the \textit{movable} points---into two groups, and
``therefore'' we sum over all possibilities.)  This of course gives
the well-known deconcatenation coproduct $\Delta:\cH\to\cH\tensor\cH$,
\begin{equation}\label{deconcat}
  \Delta: x_1\tensor x_2\tensor\dots\tensor x_n
  \mapsto
  \sum_{i=0}^n(x_1\tensor\dots\tensor x_i)
  \tensor
  (x_{i+1}\tensor\dots\tensor x_n).
\end{equation}


The above product and coproduct make $\cH(X)$ into a braided
bialgebra---which is in fact a braided Hopf algebra, our main subject
in Sec.~\ref{sec:H}.

\subsection{The legacy of screenings: a braided contour algebra}
Here, mainly for illustrative purposes, we detail the correspondence
between the ``traditional'' picture of screenings as contour integrals
and the operations making $\cH(X)$ in~\eqref{H-def} into a braided
Hopf algebra: in a particular example, we show that if elements of
$\cH(X)$ are represented in terms of contour integrals, then the Hopf
algebra axioms indeed follow from certain rules for manipulating with
these integrals.

We note that with the integral expressions in this subsection
interpreted as elements in representation spaces of a braided Hopf
algebra, the equality signs are typically to be understood as
isomorphisms between different representations of the same object;
``numerical'' identities would then follow by taking matrix elements
of the operators in relevant representations and explicitly inserting
the isomorphisms (which we do not do here).

\subsubsection{$\cH$ and its tensor powers}
We fix a basis $\scr_j$ in the braided space $X$,\footnote{In actual
  CFT{} examples, of course, the logic is just the reverse: $X$ is the
  span of elements of a preferred basis, that of screenings.} and
construct a basis in~$\cH$, $(\scr_{j_1},\dots,\scr_{j_r})\in
X^{\tensor r}$, as the integrals
\begin{equation}\label{main-int}
  (\scr_{j_1},\dots,\scr_{j_r})=
  \idotsint\limits_{-\infty<z_1<\dots<z_r<\infty} f_{j_1}(z_1)\dots f_{j_r}(z_r),
\end{equation}
where, by definition, the $f_j(z)$ (with the arguments ordered along
a line) satisfy the exchange relations
\begin{equation}\label{braiding-rel}
  f_i(u)f_j(z)=\Psi_{i,j}^{k,l}f_k(z)f_l(u),\qquad
  u>z.
\end{equation}

Elements of $\cH\tensor\cH$ are, by definition, two objects of
form~\eqref{main-int} ``placed next to each other.''
A useful convention is to represent them in the form with the two
lines joined into a single one through an intermediate point $a$:
\begin{multline}\label{main-int-tensor}
  \cH\tensor\cH\ni(\scr_{i_1},\dots,\scr_{i_r}\mid \scr_{j_1},\dots,\scr_{j_s})
  =\\
  \idotsint\limits_{-\infty<z_1<\dots<z_r<a<u_1<\dots<u_s<\infty}
  f_{i_1}(z_1)\dots f_{i_r}(z_r)\,
  f_{j_1}(u_1)\dots f_{j_s}(u_s).
\end{multline}
All such integrals form a basis of $\cH\tensor\cH$.  The point $a$ is
fixed here, and expressions \eqref{main-int-tensor} with different $a$
are isomorphic representations of the same $\cH\otimes\cH$ element
(identifying~\eqref{main-int-tensor} with the same object where $a$ is
replaced with~$b$ means identifying two equivalent representations of
the same algebraic object).

This generalizes to multiple tensor products; for example, a basis in
$\cH\tensor\cH\tensor\cH$ is given by
\begin{multline*}
  (\scr_{i_1},\dots,\scr_{i_r}\mid \scr_{j_1},\dots,\scr_{j_s}\mid \scr_{k_1},\dots,\scr_{k_p})=
  \\
  \idotsint\limits_{-\infty<z_1<\dots<z_r<a<u_1<\dots<u_s<
    b<v_1<\dots<v_p<\infty}\kern-60pt
  f_{i_1}(z_1)\dots f_{i_r}(z_r)\,
  f_{j_1}(u_1)\dots f_{j_s}(u_s)\,
  f_{k_1}(v_1)\dots f_{k_p}(v_p).
\end{multline*}

\subsubsection{Multiplication of the integrals}
We now reformulate multiplication~\eqref{H-mult-d} in terms of the
above integrals, i.e., express $m:\cH\tensor\cH\to\cH$ in the
bases~\eqref{main-int-tensor} and~\eqref{main-int}.  The
multiplication is the map in
\begin{equation*}
  m:(\scr_{i_1},\dots,\scr_{i_r}\mid \scr_{j_1},\dots,\scr_{j_s})\mapsto
  \!\!\!\idotsint\limits_{-\infty<z_1<\dots<z_r<\infty}\!\!\!\!\!
  f_{i_1}(z_1)\dots f_{i_r}(z_r)
  \!\!\!\!\!\idotsint\limits_{-\infty<u_1<\dots<u_s<\infty}\!\!\!
  f_{j_1}(u_1)\dots f_{j_s}(u_s),
\end{equation*}
where braiding relations \eqref{braiding-rel} are to be used 
to express the result in terms of basis elements~\eqref{main-int}.
In the simplest case with $r=s=1$, this becomes
\begin{align*}
  m:(\scr_i\mid \scr_j)&\mapsto
  \int_{-\infty}^\infty f_i(z)\int_{-\infty}^\infty f_j(u)
  =\iint\limits_{-\infty<z<u<\infty}\kern-10pt
  f_i(z)f_j(u)
  +
  \iint\limits_{-\infty<z<u<\infty}\kern-10pt f_i(u)f_j(z)
  \\
  \intertext{which we rewrite using the braiding as}
  &=\iint\limits_{-\infty<z<u<\infty}\kern-10pt
  f_i(z)f_j(u)
  +
  \Psi_{i,j}^{k,l}\kern-10pt\iint\limits_{-\infty<z<u<\infty}\kern-10pt
  f_k(z)f_l(u)
  =(\delta_i^k \delta_j^l + \Psi_{i,j}^{k,l})(\scr_k,\scr_l)\\  
  &=(\Bbin{}{1,1})_{i,j}^{k,l}(\scr_k,\scr_l),
\end{align*}
where $\id + \Psi = \Bbin{}{1,1}$ is the simplest quantum shuffle
(see~\bref{sec:Bbin}).
In general,
\begin{equation}\label{m-def}
  m:(\scr_{i_1},\dots,\scr_{i_r}\mid \scr_{j_1},\dots,\scr_{j_s})\mapsto
  (\Bbin{}{r,s})_{i_1,\dots,i_r,
    j_1,\dots,j_s}^{k_1,\dots,k_{r+s}}
  (\scr_{k_1},\dots,\scr_{k_{r+s}}),
\end{equation}
which reproduces the product \eqref{H-mult-d} in the bases given
in~\eqref{main-int-tensor} and~\eqref{main-int}.

\subsubsection{Comultiplication of the integrals}
The coproduct $\Delta:\cH\to\cH\tensor\cH$ defined
in~\bref{sec:comult} can be conveniently represented as
in~\eqref{cutting}, i.e., by inserting a cut into the integration
domain in~\eqref{main-int}.  The simplest case if merely
\begin{equation*}
  \Delta:\int_{-\infty}^\infty f_i(z)
  \mapsto
  \int_{-\infty}^a f_i(z)+\int_a^\infty f_i(z).
\end{equation*}
As in~\bref{sec:unit}, where the void was interpreted as the unit for
the tensor product, working now in terms of bases we interpret an
empty line, semi-infinite interval, or finite interval as the algebra
unit.  This means that $\int_{-\infty}^a f_i(z)$ is identified with
$\xymatrix@1@C=6pt{\ar@{--}|(.3){\cross}|(.5){}[0,3]
  &&&{|}\ar@{--}[0,3]&&&}$ and consequently with $\scr_i\tensor1$.
Hence, $\Delta \scr_i = \scr_i\tensor1+1\tensor \scr_i$.  In general,
cutting the line carrying $(\scr_{i_1},\dots,\scr_{i_r})$ at a point
$a$ yields $r+1$ cases $z_1<\dots<z_r<a\;$, \
$\;z_1<\dots<z_{n-1}<a<z_r\;$, \ \;\dots,\; \ $\;a<z_1<\dots<z_r$, and
hence the deconcatenation coproduct
\begin{multline}\label{Delta-def}
  \Delta: (\scr_{i_1},\dots,\scr_{i_r})
  \mapsto
  (\scr_{i_1},\dots,\scr_{i_r})\tensor 1\\
  +(\scr_{i_1},\dots,\scr_{i_{r-1}})\tensor (\scr_{i_r})
  +(\scr_{i_1},\dots,\scr_{i_{r-2}})\tensor (\scr_{i_{r-1}},\scr_{i_r})
  +\dots
  +1\tensor (\scr_{i_1},\dots,\scr_{i_{r-1}}).
\end{multline}

\subsubsection{The antipode: contour reversal}The antipode map on
integrals is naturally given by contour reversal.  Applied to
(\ref{main-int}), this operation yields
\begin{equation*}
  \A\bigl((\scr_{j_1},\dots,\scr_{j_r})\bigr)=
  (-1)^r\idotsint\limits_{-\infty<z_r<\dots<z_1<\infty}
  f_{j_1}(z_1)\dots f_{j_r}(z_r).
\end{equation*}
Using the braiding relations, we then reorder the $f_{j_k}$ such that
they occur in the same order as the integration contours.  This is
achieved by the product of braid group generators $\Psi_1
(\Psi_2\Psi_1)(\Psi_3\Psi_2\Psi_1)\dots
(\Psi_{r-1}\Psi_{r-2}\dots\Psi_1)$, a Matsumoto lift of the longest
element of the symmetric group.

\subsubsection{Verifying Hopf algebra axioms: an example}
We now show by rearranging the integrals that (a very particular
instance of) the braided bialgebra axiom \eqref{Hopf-axiom} holds for
product~\eqref{m-def} and coproduct~\eqref{Delta-def}.  The purpose is
to illustrate the fact that the braided Hopf algebra $\cH(X)$
introduced in~\bref{sec:chain-struct} is a natural formalization of
contour manipulation rules.

We evaluate both sides of~\eqref{Hopf-axiom} on the element
\begin{equation}\label{SxS}
  (\scr_i\mid \scr_j) =
  \int_{-\infty}^a f_i(z)\int_a^\infty f_j(w)\in\cH\tensor\cH,
\end{equation}
The left-hand side of~\eqref{Hopf-axiom} then becomes
\begin{multline}\label{b1-left}
  \Delta\ccirc m ((\scr_i\mid \scr_j))=\\
  (\delta_i^k\delta_j^l+\Psi_{i,j}^{k,l})
  \Bigl(\int\limits_{-\infty<z<w<a}\kern-16pt f_k(z)f_l(w)
  +\kern-12pt
  \int\limits_{-\infty<z<a<w<\infty}\kern-16pt
  f_k(z)f_l(w)
  +\kern-12pt
  \int\limits_{a<z<w<\infty}\kern-16pt f_k(z)f_l(w)\Bigr),
\end{multline}
which is simply
$(\delta_i^k\delta_j^l+\Psi_{i,j}^{k,l})\bigl((\scr_k\mid\scr_l)\tensor
1+(\scr_k)\tensor(\scr_l) +1\tensor
(\scr_k\mid\scr_l)\bigr)\in\cH\tensor\cH$.

The right-hand side of \eqref{Hopf-axiom} is $(m\tensor m)
\ccirc(\id\tensor\Psi\tensor\id)\ccirc(\Delta\tensor\Delta)$.
Applying first the coproduct to each factor in~\eqref{SxS} 
yields
\begin{multline}\label{expanded}
  (\Delta\tensor\Delta)(\scr_i\mid \scr_j)
  =\int_{-\infty}^{a_1} f_i(z)\int_{a}^{a_3}
  f_j(w)+\int_{-\infty}^{a_1} f_i(z)\int_{a_3}^{\infty} f_j(w)\\
  +\int_{a_1}^{a} f_i(z)\int_{a}^{a_3} f_j(w)+\int_{a_1}^{a}
  f_i(z)\int_{a_3}^{\infty} f_j(w),
\end{multline}
where the integrations in each term in
are split by
the three points $-\infty<a_1<a<a_3<\infty$, and hence each term is an
element of $\cH\tensor\cH\tensor\cH\tensor\cH$.  
Next, in applying $\id\tensor\Psi\tensor\id$ to~\eqref{expanded}, the
braiding is nontrivial in only one term, and hence
\begin{multline}\label{b1-preend}
  (\id\tensor\Psi\tensor\id)\ccirc(\Delta\tensor\Delta)(\scr_i\mid \scr_j)\\
  =\int_{-\infty}^{a_1}\!\!\! f_i(z)\int_{a_1}^{a}\!\!\! f_j(w)
  +\int_{-\infty}^{a_1}\!\!\! f_i(z)\int_{a_3}^{\infty}\!\!\! f_j(w)
  +\Psi_{i,j}^{k,l}\int_{a_1}^{a}\!\!\! f_k(z)\int_{a}^{a_3}\!\!\! f_l(w)
  +\int_{a}^{a_3}\!\!\! f_i(z)\int_{a_3}^{\infty}\!\!\! f_j(w).
\end{multline} 
It remains to apply $m\tensor m$ to~\eqref{b1-preend}, which amounts
to taking $a_1\to-\infty$ and $a_3\to a$ \textit{in all the lower
  integration limits}, and $a_1\to a$ and $a_3\to+\infty$ \textit{in
  all the upper limits}.  Finally,
\begin{multline*}
  (m\tensor m)
  \ccirc(\id\tensor\Psi\tensor\id)\ccirc(\Delta\tensor\Delta)(\scr_i\mid \scr_j)
  \\
  =\int_{-\infty}^{a} f_i(z)\int_{-\infty}^{a} f_j(w)
  +\int_{-\infty}^{a} f_i(z)\int_{a}^{\infty} f_j(w)\\
  {}+\Psi_{i,j}^{k,l}\int_{-\infty}^{a} f_k(z)\int_{a}^{+\infty} f_l(w)
  +\int_{a}^{+\infty} f_i(z)\int_{a}^{\infty} f_j(w),
\end{multline*}
which coincides with~\eqref{b1-left} after the use of (two) formulas
of the type $\int_{-\infty}^{a} f_i(z)\int_{-\infty}^{a}
f_j(w)=(\delta_i^k\delta_j^l+\Psi_{i,j}^{k,l})\iint\limits_{-\infty<z<w<a}
f_i(z) f_j(w)$.

\section{$\cH(X)$ and $\Nich(X)$ modules}\label{sec:H}
In this section, we study (Hopf and Yetter--Drinfeld) modules of the
braided Hopf algebra $\cH(X)$ in~\eqref{H-def}, the shuffle algebra of
a braided vector space~$X$ (whose copies were associated with
screenings$/$crosses in~\bref{sec:LS}).  Constructing $\cH(X)$-modules
requires the other braided vector space, $Y$, that we fixed in the
beginning of Sec.~\ref{sec:start-start}.  The techniques in this
section generalize to the case where $X$ and $Y$ are objects of a
braided monoidal category $(\catC,\Psi)$ (with mild additional
assumptions whenever necessary), but we avoid extra complications by
assuming the category to be that of braided linear spaces.  The
identities derived below for operators acting on the representations
in fact require nothing in addition to braiding, and can therefore be
considered ``universal,'' i.e., holding in the braid group algebra
$\mathbb{B}_N$ with a sufficiently large $N$.

All constructions for~$\cH(X)$ in~\bref{H-defs}--\bref{sec:more-Hopf}
restrict to the Nichols algebra $\Nich(X)$.

\subsection{The braided Hopf algebra $\cH(X)$}\label{H-defs}
The braided bialgebra $\cH=\cH(X)= T^{\bullet}(X)$ with the shuffle
product~\eqref{H-mult-d} and deconcatenation
coproduct~\eqref{deconcat} becomes a braided Hopf algebra with the
antipode $\A_r:X^{\tensor r}\to X^{\tensor r}$ given by a signed the
``half-twist''---\,the braid group element obtained via the Matsumoto
lift of the longest element in the symmetric group:
\begin{equation}\label{antipode-def}
  \A_r=
  (-1)^r\,
  \Psi_1 (\Psi_2\Psi_1)(\Psi_3\Psi_2\Psi_1)\dots
  (\Psi_{r-1}\Psi_{r-2}\dots\Psi_1)
\end{equation}
(with the brackets inserted to highlight the structure). \ For
example,
\begin{equation*}
  \A_5=-
  \ \begin{tangles}{l}
 \hstr{70}\vstr{50}\hx\step[1]\id\step[1]\id\step[1]\id\\
 \hstr{70}\vstr{50}\id\step[1]\hx\step[1]\id\step[1]\id\\
 \hstr{70}\vstr{50}\id\step[1]\id\step[1]\hx\step[1]\id\\
 \hstr{70}\vstr{50}\hx\step[1]\id\step[1]\hx\\
 \hstr{70}\vstr{50}\id\step[1]\hx\step[1]\id\step[1]\id\\
 \hstr{70}\vstr{50}\hx\step[1]\hx\step[1]\id\\
 \hstr{70}\vstr{50}\id\step[1]\hx\step[1]\id\step[1]\id\\
 \hstr{70}\vstr{50}\hx\step[1]\id\step[1]\id\step[1]\id\\
 \end{tangles}\ 
\end{equation*}
Also, $\A_1=-\id$.  The sign can be thought to follow from contour
reversal when the elements of $\cH(X)$ are viewed as multiple
integrals of the screenings, as in Sec.~\ref{sec:start-start}.  The
counit is $\varepsilon:X^{\tensor r}\to \delta_{r,0}k$ (identity
on~$k$).

We note that Eq.~\eqref{antipode-def} defines the action of $\A_r$ on
any tensor product of $r$ objects from~$\catC$.  With this
understanding, we also define the \textit{full twist} as the operator
acting on any $r$-fold tensor product in the category,
\begin{equation}\label{theta}
  \theta_r = \A_r\,\A_r.
\end{equation}

\subsubsection*{Notation}With $X$ fixed in what follows, we use the
short notation $(r)=X^{\tensor r}$.

\subsection{Hopf bimodules}\label{sec:Hopf-bimodule}
We introduce a category of $\cH(X)$-modules, which turn out to also be
comodules, and in fact Hopf bimodules.

With our second braided linear space~$Y$, we consider the space
\begin{align*}
  \Yspace&=\bigoplus_{s,t\geq0}X^{\tensor s}\tensor Y\tensor X^{\tensor t}\\
  &\equiv\bigoplus_{s,t\geq0}(s;Y;t),
\end{align*}
and make it into an $\cH(X)$-module bimodule as in~\eqref{x-x-act}.
In our current notation, the left action
$\cH(X)\tensor\Yspace\to\Yspace$ restricted to $(r)\tensor(s;Y;t)$ is
the map
\begin{equation*}
  \Bbin{}{r, s + 1 + t}:
  (r)\tensor(s;Y;t)\to\bigoplus_{i=0}^{r}(s+i;Y;t+r-i).
\end{equation*}
We often use the ``infix notation'' $(r)\dotact(s;Y;t)$ for it, and
write
\begin{equation}\label{r.sYt}
  (r)\dotact(s;Y;t)
  = \Bbin{}{r, s + 1 + t}\bigl((r)\tensor(s;Y;t)\bigr)
  \subset\bigoplus_{i=0}^{r}(s+i;Y;t+r-i).
\end{equation}
This map is illustrated in Fig.~\ref{fig:shuffle-intro} in the case
where $r = 2$, $s = 1$, and $t = 1$.  In each term in the figure, the
``incoming'' strands (at the top) then correspond to $X$, $X$, $X$,
$Y$, and $X$.

The right $\cH(X)$-action on $\Yspace$, also in accordance
with~\eqref{x-x-act}, is
\begin{equation}\label{sYt.r}
  (s;Y;t)\dotact (r)=\Bbin{}{s + 1 + t, r}\bigl((s;Y;t)\tensor(r)\bigr)
  \subset\bigoplus_{i=0}^{r}(s+i;Y;t+r-i).
\end{equation}

Moreover, $\Yspace$ is also a bicomodule via left and right
deconcatenations:
\begin{equation*}
  \deltaL(s;Y;t) = \bigoplus_{i=0}^{s}\,(i)\tensor({s-i};Y;t),
  \qquad
  \deltaR(r;Y;s)
  = \bigoplus_{i=0}^{s}\,({r};Y;{s-i})\tensor (i).
\end{equation*}

The following fact is well known.
\begin{prop}\label{prop:biHopf}
  With the above $\cH(X)$ action and coaction, $\Yspace$ is a Hopf
  bimodule.
\end{prop}

We note that the \textit{action} property for both left and right
actions is expressed by one of the ``basic'' quantum shuffle
identities~\eqref{BB}.  Some other identities expressing the statement
in~\bref{prop:biHopf} are given in~\bref{sec:biHopf-shuffles}.

\begin{prop}\label{prop-anti-H}
  On each graded component $(s;Y;t)$ of $\Yspace$, the ``relative
  antipode'' defined in~\eqref{rel-anti-Hopf} acts as
  \begin{equation*}
    \sigma|_{(s;Y;t)}\equiv
    \sum_{i=0}^{s}\sum_{j=0}^{t}
    \Bbin{}{s + 1 + t - j, j} \Bbin{}{i, s - i + 1 + t - j}
    \A_{i} \Shift{(s + 1 + t - j)}{\A_{j}}
    = - \A_{s + 1 + t}.
  \end{equation*}
\end{prop}
\noindent
We emphasize that the right-hand side, defined
in~\eqref{antipode-def}, acts here by a flip of all the $(s+1+t)$
strands representing $X^{\otimes s}\tensor Y\tensor X^{\otimes t}$,
$\A_{s + 1 + t}:(s;Y;t)\to(t;Y,s)$.

\subsection{Yetter--Drinfeld modules and the adjoint
  action}\label{sec:YD}
We use the Hopf bimodules discussed above to construct (left--left)
Yetter--Drinfeld $\cH(X)$ modules.  In any braided monoidal category
(with split idempotents), the space of right coinvariants of a Hopf
bimodule is a left--left Yetter--Drinfeld module under the left
adjoint action (see~\cite{[Besp-TMF],[Besp-next]} and the references
therein).  The relevant definitions are recalled in~\eqref{yd-axiom}
and~\eqref{adja}.  We now translate the adjoint action defined
in~\eqref{adja} into the shuffle language.  The main result is
in~\bref{prop-adja} below.

\subsubsection{Rewriting~\eqref{adja} in terms of shuffles}
On all of $\Yspace$, the adjoint action \eqref{adja} readily
reformulates as the map $\adjoint:\cH(X)\tensor\Yspace\to\Yspace$
whose restriction
\begin{alignat}{5}\notag
  \adjoint{}:
  (r)\tensor(s;Y;t)\to{}&\bigoplus_{i=0}^r(s+i;Y;t+r-i)\kern-160pt
  \\
  \intertext{is (again, in the ``infix notation'') given by}
  \label{r*sYt}
  (r)\adjoint (s;Y;t)
  ={}&
  \sum_{i=0}^{r}\,
  &&\Bbin{}{r-i, s + 1 + t + i}\,
  &&\xShift{}{\Bbin{\shift(r-i)}{s + 1 + t, i}}\,
  &&\Shift{(r-i + s + 1 + t)}{\A_{i}}\,
  &&\Shift{(r-i)}{\Bbraid_{i, s + 1 + t}}
  \bigl((r)\tensor(s;Y;t)\bigr).
  \\*[-6pt]
  \notag
  &
  &&\quad\dotact
  &&\quad\dotact
  &&\quad\A
  &&\quad\Psi\quad \Delta
\end{alignat}
For clarity of comparison with~\eqref{adja}, we placed the ingredients
of that formula under the corresponding pieces of the ``shuffle''
formula; the deconcatenation coproduct $\Delta$ can be conventionally
localized at the first (counting from the right) occurrence of the
index $i$.

\subsubsection{Right coinvariants in the Yetter--Drinfeld module
  $\YDspace$}\label{sec:right-coinv}
With the right coaction given also by deconcatenation, the space of
right coinvariants in $\Yspace$ is simply
\begin{equation*}
  \YDspace\equiv\bigoplus_{s\geq0}(s;Y;\,)=\bigoplus_{s\geq0}(s;Y),
\end{equation*}
where we also somewhat streamlined the notation by omitting a
redundant semicolon separating the void for $t=0$.  We now describe a
nicer formula for $(r)\adjoint (s;Y)$ than just~\eqref{r*sYt} with
$t=0$.

Let $\ \begin{tangles}{l}\vstr{70}\lu
  \object{\raisebox{5.2pt}{\tiny$\bullet$}}
\end{tangles}\ $ be the left action $\cH(X)\tensor\YDspace\to\YDspace$:
\begin{align*}
  \FromLeftii{r,s}&=\Bbin{}{r,s}:(r)\tensor(s;Y)\to(r+s;Y)\\
  \intertext{and $\ \begin{tangles}{l}
      \vstr{70}\object{\raisebox{5.1pt}{\tiny$\bullet$}}\ru
    \end{tangles}\ $ be the right action
    $\YDspace\tensor\cH(X)\to\YDspace$:}
  \FromRightii{s,r}&=\Bbin{}{s,r}\Shift{s}{\Bbraid_{1,r}}
  :(s;Y)\tensor(r) \to (s+r;Y).
\end{align*}
These left and right actions commute with each other, and
\textit{restrict to the spaces of right coinvariants
  $\YDspace\subset\Yspace$} (which we have just used by specifying the
targets as $\YDspace$, not~$\Yspace$).  We also note for the future
use that $\FromLeftii{}$ and $\FromRightii{}$ respectively satisfy the
left--left and \raisebox{1.5pt}{left}--\raisebox{-1.5pt}{right} Hopf
module axioms
\begin{equation}\label{eq:Hopf-2-2}
  \begin{tangles}{c}
    \lu\object{\raisebox{8pt}{\tiny$\bullet$}}\\[-4pt]
    \ld
  \end{tangles}
  \ = \
  \begin{tangles}{lcr}
    \hcd&&\hld\\
    \id&\hx&\id\\[-.8pt]
    \hcu&&\hlu\object{\raisebox{13pt}{\tiny$\bullet$}}
  \end{tangles}
  \qquad\text{and}\quad
  \begin{tangles}{l}
    \step\object{\raisebox{8pt}{\tiny$\bullet$}}\ru\\[-4pt]
    \ld
  \end{tangles}
  \ \ = \ \
  \begin{tangles}{l}
    \ld\step\hcd\\
    \vstr{67}\id\step\hx\step[1]\id\\
    \hcu\step\object{\raisebox{8pt}{\tiny$\bullet$}}\ru
  \end{tangles}
\end{equation}

\begin{prop}\label{prop-adja}
  The left adjoint action of $\cH(X)$ on $\YDspace$ restricted to
  $(r)\tensor (s;Y)$ is the map $(r)\tensor (s;Y)\to(s+r;Y)$ given by
  \begin{align}\notag
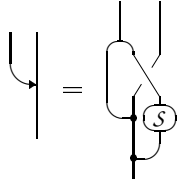
    \begin{tangles}{l}
      \vstr{200}\lu
      \object{\raisebox{18pt}{\kern-4pt\tiny$\blacktriangleright$}}
    \end{tangles}
    \ \ &= \ \
    \begin{tangles}{l}
      \hcd\step\id\\
      \vstr{80}\id\step\hx\\
      \vstr{80}\lu\object{\raisebox{6pt}{\tiny$\bullet$}}
      \step\O{\medA}\\
      \vstr{75}\step\object{\raisebox{5.8pt}{\tiny$\bullet$}}\ru
    \end{tangles}
    \\[-3pt]
    \intertext{in terms of the $\ \begin{tangles}{l}\vstr{70}\lu
        \object{\raisebox{5.2pt}{\tiny$\bullet$}}\end{tangles}\ $ and
      $\ \begin{tangles}{l}
        \vstr{70}\object{\raisebox{5.1pt}{\tiny$\bullet$}}\ru
      \end{tangles}\ $ actions,
      that is,}
    (r)\adjoint (s;Y)&=
    \sum_{i=0}^{r}
    \Bbin{}{r - i, s + i}
    \xShift{}{\Bbin{\shift(r - i)}{s, i}}
    \Shift{(r - i + s)}{(\A_{i}\,\Bbraid_{1, i}\,\Bbraid_{i, 1})}\,
    \Shift{(r - i)}{\Bbraid_{i, s}}
    \,(s+r;Y)
    \label{r-adja}
    \\
    &=\sum_{i=0}^{r} \Bbin{}{r - i, s + i} \bigl(\Bbin{}{s, i}
    \Shift{s}{\Bbraid_{1,i}} \,\Shift{(s+1)}{\A_i}\, \Bbraid_{i,s+1}
    \bigr)^{\shift(r-i)}\,(s+r;Y).\notag
  \end{align}
\end{prop}

\begin{figure}[tbp]
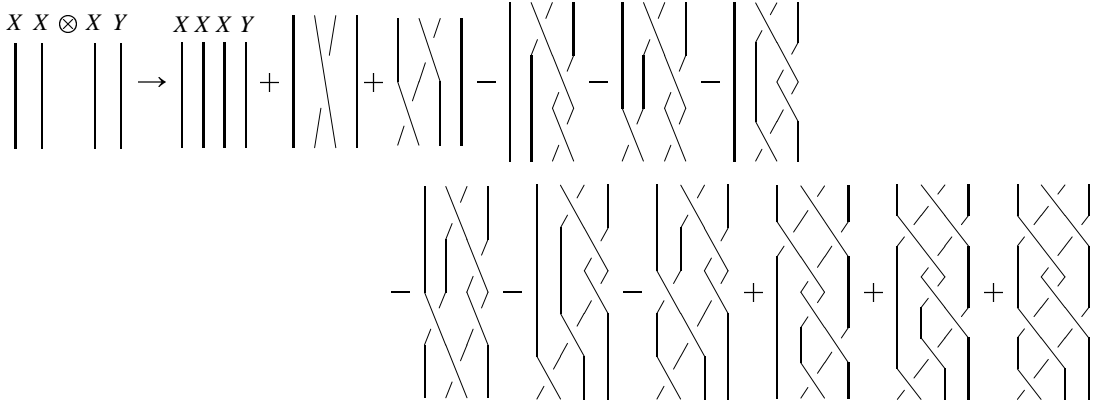

  \centering
  \begin{multline*}
\begin{tangles}{l}
      \fobject{X}\step\fobject{X}\step\fobject{\otimes}\step
      \fobject{X}\step\fobject{Y}\\
      \id\step\id\step[2]\id\step\id\\
      \id\step\id\step[2]\id\step\id\\
    \end{tangles}\
    \to\
    \begin{tangles}{l}
      \hstr{80}\fobject{X}\step\fobject{X}\step\fobject{X}\step
      \fobject{\;Y}\\
      \hstr{80}\id\step\id\step\id\step\id\\
      \hstr{80}\id\step\id\step\id\step\id\\
    \end{tangles}
    \;
      + 
      \ \begin{tangles}{l}
        \hstr{80}\vstr{250}\id\step[1]\hx\step[1]\id\\
      \end{tangles} + 
      \ \begin{tangles}{l}
        \hstr{80}\vstr{120}\id\step[1]\hx\step[1]\id\\
        \hstr{80}\vstr{120}\hx\step[1]\id\step[1]\id\\
      \end{tangles}\  - 
      \ \begin{tangles}{l}
        \hstr{80}\vstr{100}\id\step[1]\hx\step[1]\id\\
        \hstr{80}\vstr{100}\id\step[1]\id\step[1]\hx\\
        \hstr{80}\vstr{100}\id\step[1]\id\step[1]\hx\\
      \end{tangles}\  - 
      \ \begin{tangles}{l}
        \hstr{80}\vstr{100}\id\step[1]\hx\step[1]\id\\
        \hstr{80}\vstr{100}\id\step[1]\id\step[1]\hx\\
        \hstr{80}\vstr{100}\hx\step[1]\hx\\
      \end{tangles}\  - 
      \ \begin{tangles}{l}
        \hstr{80}\vstr{75}\id\step[1]\hx\step[1]\id\\
        \hstr{80}\vstr{75}\id\step[1]\id\step[1]\hx\\
        \hstr{80}\vstr{75}\id\step[1]\id\step[1]\hx\\
        \hstr{80}\vstr{75}\id\step[1]\hx\step[1]\id\\
      \end{tangles}
    \\[4pt]
    - 
    \ \begin{tangles}{l}
      \hstr{80}\vstr{100}\id\step[1]\hx\step[1]\id\\
      \hstr{80}\vstr{100}\id\step[1]\id\step[1]\hx\\
      \hstr{80}\vstr{100}\hx\step[1]\hx\\
      \hstr{80}\vstr{100}\id\step[1]\hx\step[1]\id\\
    \end{tangles}\
    - 
    \ \begin{tangles}{l}
      \hstr{90}\vstr{81}\id\step[1]\hx\step[1]\id\\
      \hstr{90}\vstr{81}\id\step[1]\id\step[1]\hx\\
      \hstr{90}\vstr{81}\id\step[1]\id\step[1]\hx\\
      \hstr{90}\vstr{81}\id\step[1]\hx\step[1]\id\\
      \hstr{90}\vstr{81}\hx\step[1]\id\step[1]\id\\
    \end{tangles}\  - 
    \ \begin{tangles}{l}
      \hstr{90}\vstr{81}\id\step[1]\hx\step[1]\id\\
      \hstr{90}\vstr{81}\id\step[1]\id\step[1]\hx\\
      \hstr{90}\vstr{81}\hx\step[1]\hx\\
      \hstr{90}\vstr{81}\id\step[1]\hx\step[1]\id\\
      \hstr{90}\vstr{81}\hx\step[1]\id\step[1]\id
    \end{tangles}\  + 
    \ \begin{tangles}{l}
      \hstr{90}\vstr{67}\id\step[1]\hx\step[1]\id\\
      \hstr{90}\vstr{67}\hx\step[1]\hx\\
      \hstr{90}\vstr{67}\id\step[1]\hx\step[1]\id\\
      \hstr{90}\vstr{67}\id\step[1]\hx\step[1]\id\\
      \hstr{90}\vstr{67}\id\step[1]\id\step[1]\hx\\
      \hstr{90}\vstr{67}\id\step[1]\hx\step[1]\id
    \end{tangles}\  + 
    \ \begin{tangles}{l}
      \hstr{90}\vstr{58}\id\step[1]\hx\step[1]\id\\
      \hstr{90}\vstr{58}\hx\step[1]\hx\\
      \hstr{90}\vstr{58}\id\step[1]\hx\step[1]\id\\
      \hstr{90}\vstr{58}\id\step[1]\hx\step[1]\id\\
      \hstr{90}\vstr{58}\id\step[1]\id\step[1]\hx\\
      \hstr{90}\vstr{58}\id\step[1]\hx\step[1]\id\\
      \hstr{90}\vstr{58}\hx\step[1]\id\step[1]\id
    \end{tangles}\  + 
    \ \begin{tangles}{l}
      \hstr{90}\vstr{58}\id\step[1]\hx\step[1]\id\\
      \hstr{90}\vstr{58}\hx\step[1]\hx\\
      \hstr{90}\vstr{58}\id\step[1]\hx\step[1]\id\\
      \hstr{90}\vstr{58}\id\step[1]\hx\step[1]\id\\
      \hstr{90}\vstr{58}\hx\step[1]\hx\\
      \hstr{90}\vstr{58}\id\step[1]\hx\step[1]\id\\
      \hstr{90}\vstr{58}\hx\step[1]\id\step[1]\id
    \end{tangles}
  \end{multline*}
  \caption[Adjoint action $(r)\protect\adjoint(s;Y)$ for $r=2$
    and $s=1$.]
    {\small Adjoint action $(r)\protect\adjoint(s;Y)$ for $r=2$ and
      $s=1$.  The ``$XY$'' assignment shown in the first term in the
      right-hand side is the same in all terms.  The first three terms
      on the right are the expansion of the first term
      in~\eqref{2-adja} and represent the shuffling of the two ``new''
      strands with the $s$ \ $X$-strands in $(s;Y)$; the last, $Y$
      strand then remains ``passive.''  The term with the minus sign
      in~\eqref{2-adja} expands into the ``minus'' terms in the
      figure, where the \textit{second} of the new strands winds
      around the $Y$ strand and returns (by the leftmost
      $\Psi_{s+2}$), to be shuffled (by
      $\xShift{}{\Bbin{\shift1}{s,1}}$) with the $s$ \ $X$-strands;
      the resulting $s+1$ strands are then shuffled
      (by~$\Bbin{}{1,s+1}$) with the \textit{first} of the new
      strands.  In the last term in~\eqref{2-adja}, both new strands
      wind around the last ($Y$-)strand: the second goes to the right
      of $Y$ via $\Psi_{s+2}\dots\Psi_{1}$ and ``waits'' there for the
      first to wind around $Y$ and return (the inner brackets), and
      then returns itself ($\Psi_{s+1}\Psi_{s+2}$), after which the
      two are shuffled by $\Bbin{}{s,2}$ with the $s$ original
      $X$-strands.  (That the two ``traveler'' strands swap their
      order in the final position reflects the presence of the
      antipode in~\eqref{r-adja}.)}
  \label{fig:adj-braids2}
\end{figure}

The structure of~\eqref{r-adja} is already seen from the lowest cases:
\begin{align}
  \label{1-adja}
  (1)\adjoint (s;Y)
  &=\bigl(\Bbin{}{1, s} - 
  \Bbin{}{s, 1} \Psi_{s + 1}\Psi_{s + 1}\Psi_{s}\dots\Psi_1
  \bigr)(s+1;Y),
  \\
  \label{2-adja}
  (2)\adjoint(s;Y)&=
  \Bigl(\Bbin{}{2, s} - 
  \Bbin{}{1, s + 1} 
  \xShift{}{\Bbin{\shift1}{s, 1}} \Psi_{s + 2} 
  \Psi_{s + 2}\dots\Psi_2
  \\*
  &\quad{}+
  \Bbin{}{s, 2} 
  \Psi_{s + 1} \Psi_{s + 2}
  \bigl(\Psi_{s + 1} \Psi_{s+1}\dots\Psi_1\bigr)
  \Psi_{s + 2}\dots\Psi_2
  \Bigr)(s+2;Y).\notag
\end{align}
The right-hand side of~\eqref{1-adja} expands into $2s+2$ terms
illustrated in Fig.~\ref{fig:adj-braids} for $s=2$, when the two
braided integers $\Bbin{}{1, s}$ and $\Bbin{}{s, 1}$ become
$\Bbin{}{1, 2}=\id + \Psi_{1} + \Psi_{2}\Psi_{1}$ and $\Bbin{}{2,
  1}=\id + \Psi_{2} + \Psi_{1}\Psi_{2}$.  The leftmost strand
representing the space $(1)=X$ in Fig.~\ref{fig:adj-braids} braids
with each of the $X$ spaces in $(s)=X^{\tensor s}$ one by one (terms
with the plus sign, all following from the expansion of
$\Bbin{}{1,s}$), or braids with the last, $(s+1)$th strand
representing~$Y$, but does not stay to the right of it; instead, it is
immediately braided with it again and, ``on its way back,'' is braided
with each of the $s$ strands (the ``minus'' terms).
The right-hand side of~\eqref{2-adja} expands similarly but more
interestingly, as is illustrated in Fig.~\ref{fig:adj-braids2}
 for $s=1$.

\subsubsection{Proof of~\eqref{r-adja}} The proof is in fact
elementary.  From~\eqref{r*sYt}, we know that the left adjoint action
on $(s;Y)$ is
\begin{equation*}
  (r)\adjoint (s;Y)
  ={}
  \sum_{i=0}^{r}\,
  \Bbin{}{r-i, s + 1 + i}\,
  \xShift{}{\Bbin{\shift(r-i)}{s + 1, i}}\,
  \Shift{(r-i + s + 1)}{\A_{i}}\,
  \Shift{(r-i)}{\Bbraid_{i, s + 1}}
  \bigl((r)\tensor(s;Y)\bigr).
\end{equation*}
The claim is therefore equivalent to the
identity
\begin{align*}
  \sum_{i=0}^{r}\, \Bbin{}{r-i, s + 1 + i}\,
  \xShift{}{\Bbin{\shift(r-i)}{s + 1, i}}\,&\Shift{(r-i + s +
    1)}{\A_{i}}\, \Shift{(r-i)}{\Bbraid_{i, s + 1}}
  \\[-6pt]
  &=
  \sum_{i=0}^{r} \Bbin{}{r - i, s + i} \bigl(\Bbin{}{s, i}
  \Shift{s}{\Bbraid_{1,i}} \,\Shift{(s+1)}{\A_i}\, \Bbraid_{i,s+1}
  \bigr)^{\shift(r-i)}.
  \\[-3pt]
  \intertext{But $\Bbin{}{s, i}\Shift{s}{\Bbraid_{1,i}}
    =\Bbin{}{s+1,i}-\Bbin{}{s+1, i-1}$ by virtue of~\eqref{Sh(s+1)},
    and we can therefore continue as}
  &=  \Bbin{}{r,s}
  +
  \sum_{i=1}^{r} \Bbin{}{r - i, s + i}
  \bigl(\Bbin{\shift{(r-i)}}{s+1,i}-\Bbin{\shift{(r-i)}}{s+1, i-1}\bigr)
  \Shift{(r-i+s+1)}{\A_i}\, \Shift{(r-i)}{\Bbraid_{i,s+1}}. 
\end{align*} %
Similar terms in the left- and the right-hand side are now combined
using~\eqref{Sh(s+1)} once again, $\Bbin{}{r-i,s+1+i}-\Bbin{}{r-i,s+i}
=\Bbin{}{r-i-1,s+1+i}\Shift{(r-i-1)}{\Bbraid_{1,s+1+i}}$, and hence
the claim equivalently reformulates as
\begin{multline*}
  \Bbin{}{r,s+1} +
  \sum_{i=1}^{r - 1}
  \Bbin{}{r - i - 1, s + 1 + i} 
  \Shift{(r - i - 1)}{\Bbraid_{1, s + 1 + i}} 
  \xShift{r - i}{\Bbin{\shift(r - i)}{{s + 1, i}}}
  \Shift{(r - i + s + 1)}{\A_{i}}
  \Shift{(r - i)}{\Bbraid_{i, s + 1}}
  \\[-6pt]
  {}=
  \Bbin{}{r,s}
  - \sum_{i=1}^{r} \Bbin{}{r - i, s + i}
  \Shift{(r - i)}{\bigl(
    \Bbin{}{s + 1, i - 1}\Shift{(s + 1)}{\A_{i}} 
    \Bbraid_{i, s + 1}
    \bigr)}.
\end{multline*}
The two isolated $\Bbin{}{}$-terms combine with the sum in the
left-hand side, extending it to $i=0$; we then shift the summation
index to the range from $1$ to $r$, as in the right-hand side.  The
resulting equality is
a term-by-term identity.\qed

\subsubsection{}
In what follows, we use the notation $\Adj_{r,s}:
(r)\tensor(s;Y)\to(r+s;Y)$ for the ``adjoint action operator''
in~\eqref{r-adja}:
\begin{equation*}
  \Adj_{r,s}(r+s;Y) = (r)\adjoint(s;Y).
\end{equation*}
The fact that the adjoint action is an \textit{action} is expressed as
\begin{equation*}
  \Adj_{r, s + t}\Shift{r}{\Adj_{s, t}}
  = \Adj_{r + s, t} \Bbin{}{r, s}
\end{equation*}
for all $r,s,t\geq0$.

\subsubsection{Remark}
We make contact with~\cite{[HS-LMS]}, where a formula for the left
adjoint action was derived based on the operator $T_r:(r; Y) \to (r;
Y)$ given by
\begin{equation*}
  T_r=(\id - \Psi_r\Psi_r\Psi_{r-1}\dots\Psi_1)\dots
  (\id - \Psi_r\Psi_r\Psi_{r-1})(\id - \Psi_r\Psi_r),\quad r\geq
  1,
\end{equation*}
and $T_0=\id$.  This operator is related to the adjoint-action
operator introduced above as
\begin{equation*}
  \Bfac{r} T_{r} = \Adj_{r, 0} \Bfac{r},
\end{equation*}
with the ``braided factorial'' defined in~\bref{sec:Bbin}.

\bigskip

Examining the operators that act only on the last, $Y$ strand in
$(r)\adjoint(s;Y)$, and then using induction on the number of strands
counted from the right, we can restate~\bref{prop-adja} in a slightly
more general form, which shows how the $\adjoint$ action on $(s+t;Y)$
``factors'' (of course, via the coproduct) through $\adjoint$ acting
on the \textit{last} several strands, $(t;Y)$.  The following
proposition expresses a formula of the type
$(r)\adjoint(s+t;Y)=\sum\dots ((r')\adjoint(t;Y))\dots$.

\begin{prop}\label{prop:Adj2}
  The left adjoint action satisfies the identity
  \begin{equation}\label{eq:Adj2}
    \begin{tangles}{l}
      \fobject{(r)}\step[3]\fobject{\ (s+t;Y)}\\
      \vstr{550}\lu[3]
      \object{\raisebox{53pt}{\kern-4pt\tiny$\blacktriangleright$}}
    \end{tangles}\quad=\ \
    \begin{tangles}{l}
      \step\fobject{(r)}\step[3]\fobject{(s)}\step[2]\fobject{\ (t;Y)}\\
      \cd\step[2]\id\step[2]\id\\
      \vstr{70}\id\step[2]\x\step[2]\id\\
      \vstr{70}\lu[2] \object{\raisebox{5pt}{\tiny$\bullet$}}
      \step\cd\step[1]\id\\
      \vstr{50}\step[2]\id\step\id\step[2]\hx\\
      \vstr{70}\step[2]\d\lu[2]
      \object{\raisebox{5pt}{\kern-4.5pt\tiny$\blacktriangleright$}}
      \step\id\\
      \vstr{50}\step[3]\id\step[2]\hx\\
      \vstr{70}\step[3]\d\step[1]\O{\medA}\step\id\\
      \vstr{70}\step[4] \object{\raisebox{5pt}{\tiny$\bullet$}}\ru
      \step[1]\id
    \end{tangles}
  \end{equation}
  with the same left and right actions $\FromLeftiiD$ and
  $\FromRightiiD$  
  as above.
\end{prop}

\begin{prop}[cf.~\cite{[Besp-next]}]\label{s2-eval}
  On each graded component $(s;Y)$ of $\YDspace$, the ``squared
  relative antipode''~\eqref{rel-anti-YD} acts as the full twist
  $\theta_{s+1}$:
  \begin{equation*}
    \sigma_2\bigr|_{(s;Y)}
    \equiv \sum_{i=0}^{s}\Adj_{i, s - i} \A_{i} = \theta_{s+1}.
  \end{equation*}
\end{prop}

We emphasize that the full twist defined in~\eqref{theta} is here
applied to all of the $s+1$ strands in $(s;Y)$; the relation nicely
illustrates why $\sigma_2$ is a \textit{squared} ``antipode'' in the
terminology of~\cite{[Besp-TMF],[Besp-next]}.

\subsection{Multivertex Hopf bimodules}\label{sec:more-Hopf}
Hopf bimodules can be further constructed as $\cY\tensor\cZ$, where
$\cY$ is a Yetter--Drinfeld module and $\cZ$ is a Hopf
bimodule~\cite{[Besp-Dr-(Bi)]}.  The tensor product is a Hopf bimodule
under the right action and right coaction induced from those of $\cZ$,
and under the diagonal left action and codiagonal left coaction.  For
a Hopf bimodule $\Zspace$ and a Yetter--Drinfeld module $\YDspace$,
the resulting Hopf bimodule has the graded components $(s;Y;t;Z;u)$.

Evidently, the construction can be iterated by further taking tensor
products with Yetter--Drinfeld module(s); it is associative in the
sense that the Yetter--Drinfeld modules can be multiplied first (to
produce a Yetter--Drinfeld module under the diagonal action and
codiagonal coaction).  The resulting Hopf bimodules have the graded
components
\begin{equation}\label{homo-tprod}
  (s_1;Y_1;s_2;Y_2;\dots;s_n;Y_n;s_{n+1})
  =X^{\tensor s_1}\tensor Y_1\tensor X^{\tensor s_2}\tensor
  Y_2\tensor\dots\tensor X^{s_n}\tensor Y_n\tensor X^{s_{n+1}}
\end{equation}
(in our setting, we only encounter cases where $Z=Y_i=Y$, but it is
useful to allow ``different $Y$'' spaces in order to clarify the
structure of formulas).

\subsection{Introducing $\cH^*$}\label{sec:H*}
We introduce a ``conjugate'' $E$ to the adjoint action
$F\equiv(1)\adjoint{}$.  Although this does not allow obtaining even a
decent associative algebra in general, it is interesting to see that
an ``$E F - F E$'' formula can nevertheless be written with a
meaningful right-hand side.

For $\cH(X)$ in~\eqref{H-def}, we introduce the dual Hopf algebra $\cA$
built on the tensor algebra of $A=X^*$:
\begin{align*}
  \cA&=\bigoplus_{r\geq0}A^{\tensor n}, \\
  &=(\,)\oplus(-1)\oplus(-2)\oplus\dots.
\end{align*}
The pairing $\rho:A\tensor X\to k$ is denoted by \raisebox{-6pt}{\
  $\begin{tangles}{l}\ev
  \end{tangles}$\ }.  It is extended to the tensor products as
\begin{equation*}
  \rho(A^{\tensor n},X^{\tensor n}) =
  \ \begin{tangles}{l}
    \vstr{100}\id
    \step\object{\dots}\step
    \id\step\id\step[2]\id\step[1]\id\step[1]\object{\dots}\step[1]\id
    \\[-4pt]
    \vstr{67}\step[1]
    {\makeatletter\@ev{0,\hm@de}{30,\hm@detens}{80}b\makeatother}
    \step[2]{\makeatletter\@ev{0,\hm@de}{10,\hm@detens}{40}b\makeatother}
    \ev
  \end{tangles}
\end{equation*}

\bigskip

\noindent

Then any left $\cH(X)$ comodule becomes a left $\cA$ module under the
action
\begin{equation*}
  \begin{tangles}{l}
    \fobject{\cA}\step[3]\fobject{\cH}\\
    \id\step[1]\ld\\
    \hh\ev\step\id
  \end{tangles}
  \ \ : \ \alpha\tensor h\mapsto \alpha\leftii h =
  \rho(\alpha,h\mone)h\zero,
  \qquad
  \alpha\in\cA,\quad h\in\cH.
\end{equation*}
With the deconcatenation coaction on $\Yspace$,
we then have the map $(-1)\tensor (s;Y)\to(s-1;Y)$ given by
\begin{equation}\label{-1-adja}
  \begin{tangles}{l}
    \fobject{A}
    \step[3]\fobject{\ X}\step\fobject{\ X}
    \step[3]\fobject{\ X}\step\fobject{\ \ Y}
    \step[3]\fobject{\ A}\step[2]\fobject{\ X}
    \step\fobject{\ \ X}\step[3]\fobject{\ \ X}\step\fobject{\ \ \ Y}\\
    \vstr{160}\id\step\raisebox{12pt}{$\tensor$}\step
    \id\step\id\step\raisebox{12pt}{\dots\ \ }\id\step\id
    \raisebox{12pt}{\ \ $\to$\ }
    \step\ev\step\id\step\raisebox{12pt}{\dots\ \ }\id\step\id
  \end{tangles}
\end{equation}
(with $s$ \ $X$-strands in the left-hand side and $s-1$ ``free'' ones
in the right-hand side).

We introduce a new (and hopefully suggestive) notation for the maps
in~\eqref{-1-adja} and~\eqref{1-adja}, writing them as
\begin{alignat*}{2}
  \mathsf{e}_s &:\ & (-1)\tensor (s;Y)&\to(s-1;Y),\\
  \mathsf{f}_s&:&(1)\tensor (s;Y)&\to (s+1;Y),
\end{alignat*}
and ``compare'' the results of applying $\mathsf{e}$ and $\mathsf{f}$
in different orders.  In considering the product $A\tensor X\tensor
(s;Y)=(-1)\tensor(1)\tensor (s;Y)$, we \textit{enumerate} the factors
also as $-1$, $1$, $2$, \dots, $s+2$.  Then, in particular, $\Psi_1$
(see~\bref{sec:notation})
acts here as $\id\tensor\Psi\tensor\id^{\tensor s}$, and the pairing
between the $A$ strand and the leftmost $X$ strand is conveniently
denoted as~$\rho_{-1,1}$.

\begin{lemma}\label{lemma:K2}
The following identity holds for maps
  $(-1)\tensor(1)\tensor(s;Y)\to\tensor(s;Y)$\textup{:}
  \begin{equation*}
      \mathsf{e}_{s+1}\ccirc\mathsf{f}_s -
      \mathsf{f}_{s-1}\ccirc
      \mathsf{e}_s\ccirc\Psi_1
      = \rho_{-1,1}\tensor\id^{\tensor(s+1)}
      -
      \KK(s+1),
  \end{equation*}
  where $\KK(s+1):(-1)\tensor(1)\tensor(s;Y)\to(s;Y)$ is the map
  \begin{equation*}
    \KK(s+1)
    =\rho_{-1,1}\ccirc(\Psi_1\dots\Psi_{s+1}\Psi_{s+1}\dots\Psi_1)
    \ : \
    \begin{tangles}{l}
    \fobject{A}\step[2]\fobject{X}\step[1]\fobject{X}\step\fobject{X}
    \step[3]\fobject{X}\step\fobject{Y}\\
    \vstr{62}\id\step[2]\hx\step[1]\id\step[3]\id\step[1]\id\\
    \vstr{62}\id\step[2]\id\step[1]\hx\step[3]\id\step[1]\id\\
    \vstr{100}\id\step[2]\id\step[1]\id\step[2]\object{\ \dots\ }\step[2]
    \id\step\id\\
    \vstr{62}\id\step[2]\id\step[1]\id\step[3]\hx\step[1]\id\\
    \vstr{62}\id\step[2]\id\step[1]\id\step[3]\id\step\hx\\
    \vstr{62}\id\step[2]\id\step[1]\id\step[3]\id\step\hx\\
    \vstr{62}\id\step[2]\id\step[1]\id\step[3]\hx\step\id\\
    \vstr{100}\id\step[2]\id\step[1]\id\step[2]\object{\ \dots\ }\step[2]
    \id\step\id\\
    \vstr{62}\id\step[2]\id\step\hx\step[3]\id\step\id\\
    \vstr{62}\id\step[2]\hx\step\id\step[3]\id\step\id\\
    \vstr{62}\ev\step\id\step\id\step[3]\id\step\id\\
    \step[3]\fobject{2}\step\fobject{3}\step[3]\fobject{s+1\quad}
    \step\fobject{\quad s+2}
  \end{tangles}
  \end{equation*}
\end{lemma}
(the leftmost $X$ strand is passed around the bunch of $s+1$ strands,
and upon returning is contracted with the $A$ strand).

The proof is straightforward.  
On $(-1)\tensor(1)\tensor(s;Y)$, we
first apply $\mathsf{e}_s$ to $(s;Y)$, for which we need to contract
$(-1)$ with the first $X$ strand in $(s;Y)$.  For this, we first move
this strand to the left by $\Psi_1$:
\begin{equation*}
  \rho_{-1,1}\ccirc\Psi_1:\
  \begin{tangles}{l}
    \vstr{30}\id\step\id\step\id\step\id\step[2]\id\step\id\\
    \vstr{80}\id\step\hx\step\id\step[1]\object{\dots}
    \step[1]\id\step\id\\
    \hev\step[1]\id\step\id\step[2]\id\step\id
  \end{tangles}
\end{equation*}
Next applying $\mathsf{f}_{s-1}$ gives the map
\begin{equation*}
  (\rho_{-1,1}\tensor\Shift{1}{\mathsf{f}_{s-1}})\ccirc\Psi_1
  =\rho_{-1,1}\ccirc(\id\tensor\Shift{1}{\mathsf{f}_{s-1}})\ccirc\Psi_1
  :
  (-1)\tensor(1)\tensor(s;Y)\to(s;Y).
\end{equation*}

The first-$\mathsf{f}$-then-$\mathsf{e}$ variant simply gives the map
\begin{equation*}
  \rho_{-1,1}\ccirc(\id\tensor\mathsf{f}_s)
  :
  (-1)\tensor(1)\tensor(s;Y)\to(s;Y).
\end{equation*}
The maps to the right of $\rho$ in the last two formulas satisfy the
identity
$\mathsf{f}_s
-
\Shift{1}{\mathsf{f}_{s-1}} \Psi_1  
= \id
-
\Psi_1\dots\Psi_{s+1}\Psi_{s+1}\dots\Psi_1$,
as is easily verified from the definition.  The formula in the lemma
now follows immediately.

In the case of \textit{diagonal braiding},
the formula for $\KK$ in the lemma reduces nicely, to the product of
squared braiding with each of the modules in~$(s;Y)=X^{\tensor
  s}\tensor Y$.  We encounter a particular example in~\eqref{qcomm}.

\section{Fusion product}\label{sec:FUSION}
Our aim in this section is to define a product of Yetter--Drinfeld
modules of right coinvariants introduced in~\bref{sec:YD}.  As a
linear space, this is the tensor product, but it carries the
\textit{adjoint} representation of~$\cH(X)$ (the ``cumulative''
adjoint, as we discuss in what follows).  This construction plays a
role in CFT{} applications, but may also be interesting in the
``abstract'' setting.  Geometrically, in the spirit of
Sec.~\ref{sec:start-start}, fusion corresponds to a degeneration of
the picture as in Fig.~\ref{fig:intro}, when several punctures are
allowed to sit on the same line.  On the algebraic side, fusion gives
an ample supply of Yetter--Drinfeld $\cH(X)$ and $\Nich(X)$ modules.

\subsection{Fusion product of Hopf bimodules}\label{sec:fusion-H}
For Hopf bimodules described in the preceding section, their fusion
product is the ``$\mapX$'' map defined in~\eqref{chi-map}.  Its
categorial meaning is that it projects onto the cotensor product.  We
now write~\eqref{chi-map} in terms of the braid group algebra
generators (as a ``quantum shuffle'' formula).  The governing idea is
that for subspaces $(s_1;Y;s_2)$ and $(t_1;Z;t_2)$ of two Hopf
bimodules (see~\bref{sec:Hopf-bimodule}), their product
$(s_1;Y;s_2)\odot(t_1;Z;t_2)$ is given by ``propagation'' of the
$X$-strands that ``face the other module'' to that module: the $s_2$
strands split as $(s_2-i)+i$, and the $t_1$, similarly, as
$j+(t_1-j)$; then $(i)$ and $(j)$ are braided with each other, after
which $(i)$ acts to the right, via \eqref{r.sYt}
($\Bbin{}{i,t_1-j+1+t_2}$), and $(j)$ acts to the left, via
\eqref{sYt.r} ($\Bbin{}{s_1+1+s_2-i,j}$):
\begin{multline}\label{2-fusion}
  (s_1;Y;s_2)\odot(t_1;Z;t_2)=
  \\
  =\sum_{i=0}^{s_{2}}\sum_{j=0}^{t_1}
  \Bbin{}{s_1 + 1 + s_2 - i, j}\,
  \xShift{}{\Bbin{\shift(j + s_1 + 1 + s_2 - i)}{i,
     t_1 + 1 + t_2 - j}}\,
  \Shift{(s_1 + 1 + s_2 - i)}{\Bbraid_{i,j}}
  (s_1;Y_1;s_2+t_1;Z_1;t_2).
\end{multline}

More generally, for two multivertex Hopf bimodules
from~\bref{sec:more-Hopf}, the fusion product of their respective
subspaces $(s_1;Y_1;\dots;s_n;Y_n;s_{n+1})$ and
$(t_1;Z_1;\dots;t_m;Z_m;t_{m+1})$ is given by
\begin{multline}\label{gen-fusion}
  (s_1;Y_1;\dots;s_n;Y_n;s_{n+1})\odot
  (t_1;Z_1;\dots;t_m;Z_m;t_{m+1})=
  \\
  =\sum_{i=0}^{s_{n+1}}\sum_{j=0}^{t_1}
  \Bbin{}{n + \sum_{a=1}^{n+1}s_a - i, j}
  \xShift{}{\Bbin{\shift(n + j + \sum_{a=1}^{n+1}s_a - i)}{i,
      m + \sum_{a=1}^{m+1}t_a - j}}
  \Shift{(n + \sum_{a=1}^{n+1}s_a - i)}{\Bbraid_{i,j}}\\
  \times(s_1;Y_1;\dots;s_n;Y_n;s_{n+1}+t_1;Z_1;\dots;t_m;Z_m;t_{m+1})
\end{multline}

\subsection{Fusion product of Yetter--Drinfeld modules}
For Yetter--Drinfeld modules given by right coinvariants in Hopf
bimodules, Eq.~\eqref{2-fusion} is simplified to
\begin{align}\label{YDfusion}
  (s;Y)\odot(t;Z)&=
  \sum_{j=0}^{t}
  \Bbin{}{s, j} \Shift{s}{\Bbraid_{1, j}}(s;Y;t;Z)
  \\
  \intertext{(with $(s;Y;t;Z)=X^{\tensor s}\tensor Y\tensor X^{\tensor
      t}\tensor z$ as before), where we recognize the right action
    introduced in~\bref{sec:right-coinv},}
  &= \sum_{j=0}^{t}\FromRightii{s,j} (s;Y;t;Z)
  = \ \
  \begin{tangles}{l}
    \fobject{(s;Y)}\step[2]\fobject{\ (t;Z)}\\[-2pt]
    \id\step[1]\ld\\
    \object{\raisebox{8pt}{\tiny$\bullet$}}\ru\step[1]\id
  \end{tangles}\notag
\end{align}

\noindent
The last diagram is of course a ``dotted'' version of the $\mapI$ map
studied in~\bref{sec:mapI}, and we use the relevant results
from~\bref{sec:mapI} shortly.

The next proposition describes how $\cH(X)$ (and $\Nich(X)$) coacts
and acts on the right coinvariants obtained via the $\odot$ product.

\begin{prop}\label{prop:Adja2}
  On the graded components $(s;Y;t;Z)=X^{\otimes s}\tensor Y\tensor
  X^{\otimes t}\tensor Z$ of a fusion product, the $\cH(X)$ coaction
  is ``by deconcatenation up to the first vertex,''
  \begin{align*}
    (s;Y;t;Z)&\to\sum_{i=0}^{s}(i)\tensor(s-i;Y;t;Z)\\
    \intertext{and the $\cH(X)$ action is the ``cumulative'' adjoint
      action
    }
    (r)\tensor(s;Y;t;Z)&\to
    \Adj_{r,s+1+t}(r+s;Y;t;Z)
    \subset\bigoplus_{j=0}^{r}(s+j;Y;t+r-j;Z).
  \end{align*}
\end{prop}
These statements are the result of a calculation showing that the
above coaction and action are intertwined by $\odot$ with the standard
coaction and action on tensor products of Yetter--Drinfeld modules:
  \begin{equation}\label{2-diags}
  \begin{tangles}{l}
    \ld\step[1]\ld[1]\\
    \id\step[1]\hx\step[1]\d\\
    \hcu\step[1]\id\step\ld\\
    \step[.5]\id\step[1.5]
    \object{\raisebox{8pt}{\tiny$\bullet$}}\ru\step\id
  \end{tangles}
  \ \ = \ \
  \begin{tangles}{l}
    \step\id\step\ld\\
    \step\object{\raisebox{8pt}{\tiny$\bullet$}}\ru\step\id\\
    \vstr{133}\ld\step[2]\id
  \end{tangles}
  \qquad\text{and}\qquad
    \begin{tangles}{l}
      \step[.25]
      \fobject{(r)}\step[1.5]\fobject{(s;Y)}\step[1.75]\fobject{(t;Z)}\\
      \hcd\step[1]\id\step\id\\
      \vstr{80}\id\step[1]\hx\step\id\\
      \vstr{80}\lu\object{\raisebox{6pt}{\kern-4pt\tiny$\blacktriangleright$}}
      \step[1]\lu\object{\raisebox{6pt}{\kern-4pt\tiny$\blacktriangleright$}}\\
      \vstr{80}\step[1]\id\step\ld\\
      \vstr{80}\step[1]\object{\raisebox{6pt}{\tiny$\bullet$}}\ru\step\id
    \end{tangles}
    \ \ = \ \ 
    \begin{tangles}{l}
      \vstr{150}\id\step[1]\id\step\ld\\
      \vstr{150}\id\step[1]\object{\raisebox{13pt}{\tiny$\bullet$}}\ru
      \step\id\\
      \lu\object{\raisebox{8pt}{\kern-4pt\tiny$\blacktriangleright$}}
      \step[2]\id\object{\kern-20pt\rule[8pt]{20pt}{4pt}}
    \end{tangles}
  \end{equation}
  where
  $\begin{tangles}{l}
    \lu\object{\raisebox{8pt}{\kern-4pt\tiny$\blacktriangleright$}}
    \step[1]\id\object{\kern-10pt\rule[8pt]{10pt}{4pt}}
  \end{tangles}\ $ denotes the ``cumulative'' adjoint action
  $(r)\adjoint(s;Y;t;Z)=
  \Adj_{r,s+1+t}(r+s;Y;t;Z)$.

The first intertwining identity in~\eqref{2-diags}, for coactions,
follows by noticing that it is a ``dotted'' version of
\eqref{comodule-morphism} (for the right action $\FromRightii{}=
\ \begin{tangles}{l}
      \vstr{70}\object{\raisebox{5.1pt}{\tiny$\bullet$}}\ru
    \end{tangles}\ $ from \bref{sec:right-coinv}), and is proved the
    same \textit{because of the second property in~\eqref{eq:Hopf-2-2}
      for $\FromRightii{}$}; just this property underlies the
    identity.  To prove the second identity in~\eqref{2-diags}, we
    combine~\bref{prop-adja} with the calculations in~\bref{sec:mapI}
    to obtain the chain of identities
\begin{equation*}
  \begin{tangles}{l}
    \vstr{133}\hcd\step[1]\id\step\id\\
    \vstr{133}\id\step[1]\hx\step\id\\
    \vstr{133}\lu
    \object{\raisebox{11.2pt}{\kern-4pt\tiny$\blacktriangleright$}}
    \vstr{133}\step[1]\lu
    \object{\raisebox{11.2pt}{\kern-4pt\tiny$\blacktriangleright$}}\\
    \vstr{133}\step[1]\id\step\ld\\
    \vstr{133}\step[1]\object{\raisebox{12pt}{\tiny$\bullet$}}\ru\step\id
  \end{tangles}
  \ \ \stackrel{\bref{prop-adja}}{=} \ \
  \begin{tangles}{l}
      \step[.5]\cd\step[1]\id\step[1.5]\id\\
      \hcd\step[1.5]\hx\step[1.5]\id\\
      \id\step\id\step\ddh\step[.5]\hcd\step[1]\id\\
      \id\step\hx\step[1]\id\step\hx\\
      \lu\object{\raisebox{8pt}{\tiny$\bullet$}}
      \step\O{\medA}\step[.5]\step[.5]\lu
      \object{\raisebox{8pt}{\tiny$\bullet$}}\step\O{\medA}\\
      \step[1]\id\object{\raisebox{8pt}{\tiny$\bullet$}}\ru\step[1]
      \step\object{\raisebox{8pt}{\tiny$\bullet$}}\ru\step[-1.5]\hld\\
      \step[1]\hd\step[1]\ne1\step[.5]\id\\
      \step[1.5]\object{\raisebox{8pt}{\tiny$\bullet$}}\ru\step[1.5]\id
    \end{tangles}
  \ \ \ \stackrel{\eqref{adja-intertwine}}{=} \ \
  \begin{tangles}{l}
    \vstr{80}\step[.5]\id\step[1.5]\id\step[1]\ld\\
    \hcd\step[1]\object{\raisebox{8pt}{\tiny$\bullet$}}\ru\step[1]\id\\
    \vstr{70}\id\step[1]\hx\step[2]\id\\
    \lu\object{\raisebox{8pt}{\tiny$\bullet$}}\step[.5]\hcd\step[1]\ddh\\
    \vstr{70}\step[1]\id\step[.5]\id\step[1]\hx\\
    \step[1]\id\step[.5]\lu\object{\raisebox{8pt}{\tiny$\bullet$}}
    \step[.5]\hcd\\
    \vstr{70}\step[1]\id\step[1]\hdd\step[.5]\O{\medA}\step[1]\id\\
    \step[1]\id\step[1]\object{\raisebox{8pt}{\tiny$\bullet$}}\ru\dd\\
    \vstr{70}\step[1]\id\step[1]\hx\\
    \vstr{70}\step[1]\id\step[1]\O{\medA}\step[1]\hd\\
    \vstr{60}\step[1]
    \object{\raisebox{4pt}{\tiny$\bullet$}}\ru\step[1.5]\id\\
  \end{tangles}
  \ \ \ = \ \
  \begin{tangles}{l}
    \vstr{80}\step[.5]\id\step[1.5]\id\step[1]\ld\\
    \hcd\step[1]\object{\raisebox{8pt}{\tiny$\bullet$}}\ru\step[1]\id\\
    \vstr{70}\id\step[1]\hx\step[2]\id\\
    \lu\object{\raisebox{8pt}{\tiny$\bullet$}}\step[.5]\hcd\step[1]\ddh\\
    \vstr{70}\step[1]\id\step[.5]\id\step[1]\hx\\
    \step[1]\id\step[.5]\lu
    \object{\raisebox{8pt}{\kern-4pt\tiny$\blacktriangleright$}}
    \step[1]\O{\medA}\\
    \step[1]\hd\step[1]\hx\\
    \step[1.5]\object{\raisebox{8pt}{\tiny$\bullet$}}\ru\step\id
  \end{tangles}
  \ \ \stackrel{\bref{prop:Adj2}}{=} \ \
  \begin{tangles}{l}
    \vstr{150}\id\step[1]\id\step\ld\\
    \vstr{150}\id\step[1]\object{\raisebox{13pt}{\tiny$\bullet$}}\ru
    \step\id\\
    \vstr{150}
    \lu\object{\raisebox{12.8pt}{\kern-4pt\tiny$\blacktriangleright$}}
    \step[2]\id\object{\kern-20pt\rule[13pt]{20pt}{4pt}}
  \end{tangles}
\end{equation*}
Here, invoking \eqref{adja-intertwine} refers to its ``dotted''
version (with $\FromRightii{}$ from \bref{sec:right-coinv} for all
right actions and with $\FromLeftii{}$ for all left actions except the
adjoint), whose proof is identical to the proof given in
Fig.~\ref{fig:weary}, again because properties~\eqref{eq:Hopf-2-2}
hold (both are needed this time).  The next, unlabeled equality in the
above chain once again uses \bref{prop-adja} to recover the adjoint
action.

\subsubsection{Two-vertex Yetter--Drinfeld
  modules}\label{sec:YD-Adja2}
\textit{The Yetter--Drinfeld axiom holds for the coaction and action
  in~\bref{prop:Adja2}},
\begin{equation}\label{eq:YD2}
  \begin{tangles}{l}
    \hcd\step\id\step\id\\
    \vstr{50}\id\step\hx\step\id\\
    \vstr{50}\id\step\id\step\hx\\[-1pt]
    \vstr{50}\lu
    \object{\raisebox{3pt}{\kern-4pt\tiny$\blacktriangleright$}}
    \step[1]\object{\kern-10pt\rule[3pt]{10pt}{4pt}}\id\step\id\\
    \vstr{67}\ld\step\hx\\
    \vstr{67}\id\step\hx\step\id\\
    \hcu\step\id\step\id
  \end{tangles}
  \ \ \ \ = \ \
  \begin{tangles}{l}
    \vstr{150}\cd\step[1]\ld\step\id\\
    \id\step[2]\hx\step\id\step\id\\
    \vstr{150}\cu\step[1]
    \lu\object{\raisebox{13pt}{\kern-4pt\tiny$\blacktriangleright$}}
    \step[1]\id\object{\kern-10pt\rule[13pt]{10pt}{4pt}}
  \end{tangles}
\end{equation}
simply because both diagram identities in~\eqref{2-diags} are
\textit{intertwining} formulas.  We thus have \textit{two-vertex}
Yetter--Drinfeld modules with graded components $(s;Z;t;Y)=X^{\otimes
  s}\tensor Z\tensor X^{\otimes t}\tensor Y$.



\subsubsection{Multivertex Yetter--Drinfeld modules}\label{multi-YD}
More Yetter--Drinfeld modules can be constructed via iterated fusion
product, using the associativity of the $\mapI$ map
(see~\bref{sec:mapI}).  For example, in evaluating the product
$\ZDspace\odot\bigl(\YDspace\odot \UDspace\bigr)$, both $\ZDspace$ and
$\YDspace\odot\UDspace$ are Yetter--Drinfeld modules under the left
adjoint action; the construction is iterated smoothly to yield
Yetter--Drinfeld modules with graded components
$(s_1;Y_1,s_2;Y_2;\dots;s_n;Y_n)$ carrying the ``cumulative'' adjoint
action
\begin{equation*}
  (r)\adjoint(s_1;Y_1,s_2;Y_2;\dots;s_n;Y_n)=
  \Adja_{r,s_1+\dots+s_n + n -
    1}(r+s_1;Y_1,s_2;Y_2;\dots;s_n;Y_n)
\end{equation*}
and the coaction by ``deconcatenation up to $Y_1$.''

\subsubsection{A left-action reformulation of the fusion}
Another useful reformulation of fusion product~\eqref{YDfusion} is
\begin{align}\notag
  (s;Y)\odot(t;Z)
  &= \FromLeftii{s+1,t}\bigl((s;Y)\tensor(t;Z)\bigr)
  \equiv
  (s;Y)\FromLeftiiD (t;Z),
\end{align}
where we use the definition of $\FromLeftii{}$
from~\bref{sec:right-coinv}, and where \textit{all the $s+1$ strands
  in $(s;Y)$ ``act'' on $(t;Z)$}, i.e., $\FromLeftii{}$ 
takes $X^{\tensor s}\tensor Y$ as its first argument.

\subsection{Braiding of the fusion product of Yetter--Drinfeld
  modules}
\begin{prop}\label{prop:FYDB}
  For two $\YDspace$-type Yetter--Drinfeld modules, the braiding
  \begin{equation*}
    \Fbraid:\YDspace\odot\ZDspace\to\ZDspace\odot\YDspace
  \end{equation*}
  is given on braided components by the maps
  \begin{gather*}
    \Fbraid_{s, t}:
    (s;Y;t;Z)\to\bigoplus_{i=0}^{s}(i;Z;s-i+t;Y),\\
    \Fbraid_{s, t}
    =-\FromRightii{s, t + 1} \Shift{(s + 1)}{\A_{t + 1}}.
  \end{gather*}  
\end{prop}

We emphasize that the antipode here formally acts,
via~\eqref{antipode-def}, on all of the $t+1$ strands of $(t;Z)$. \
All the $t+1$ strands, moreover, then act from the right on $(s;Y)$
via an obvious extension of the right action
$\FromRightii{}=\!\!\!\FromRightiiD$ defined
in~\bref{sec:right-coinv}.

\subsubsection{Proof of~\bref{prop:FYDB}}
We find the braiding directly, by calculating how the standard
Yetter--Drin\-feld braiding (see~\eqref{eq:YDbraiding}) is intertwined
with the $\mapI$ map (see~\bref{sec:mapI}).  For this, we represent
$\mapI\ccirc
(\text{standard Yetter--Drin\-feld braiding})$ in a form
$\Fbraid\ccirc\mapI$ via the following transformations:
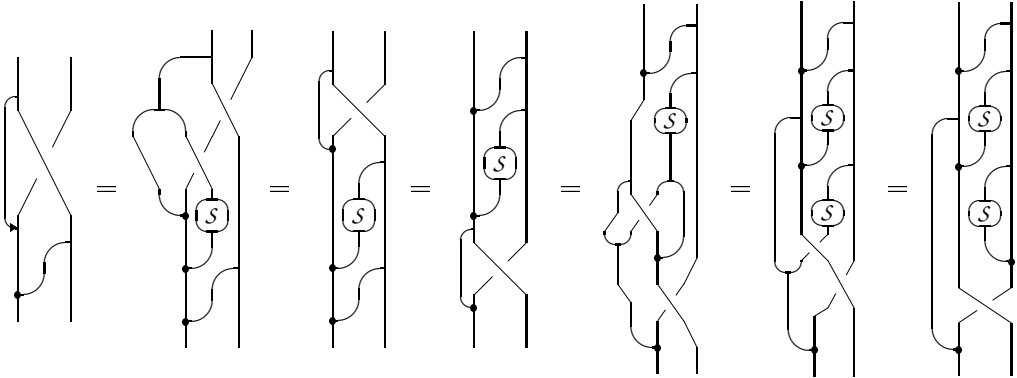
\begin{equation*}
  \begin{tangles}{l}
    \hld\step[2]\id\\
    \vstr{200}\id\step[.5]\x\\[-1pt]
    \hlu\object{\raisebox{13.3pt}{\kern-4pt\tiny$\blacktriangleright$}}
    \step[1]\ld\\
    \step[.5]\object{\raisebox{8pt}{\tiny$\bullet$}}\ru\step[1]\id
  \end{tangles}
  \ \ = \
  \begin{tangles}{l}
    \step[1]\ld[2]\step[1]\ddh\\
    \cd\step[1]\hx\\
    \d\step[1]\hx\step[1]\id\\
    \step[1]\lu[1]\object{\raisebox{8pt}{\tiny$\bullet$}}\step\O{\medA}
    \step[1]\id\\
    \step[2]\object{\raisebox{8pt}{\tiny$\bullet$}}\ru
    \ld\\
    \step[2]\object{\raisebox{8pt}{\tiny$\bullet$}}\ru\step[1]\id
  \end{tangles}
  \ = \
  \begin{tangles}{l}
    \step[.5]\hld\step[2]\id\\
    \step[.5]\id\step[.5]\x\\[-1pt]
    \step[.5]\hlu\object{\raisebox{13pt}{\tiny$\bullet$}}
    \step[1]\ld\\
    \step[1]\id\step[1]\O{\medA}\step[1]\id\\
    \step[1]\object{\raisebox{8pt}{\tiny$\bullet$}}\ru\ld\\
    \step[1]\object{\raisebox{8pt}{\tiny$\bullet$}}\ru\step\id
  \end{tangles}
  \ \ = \
  \begin{tangles}{l}
    \step[1]\id\step[1]\ld\\
    \step[1]\object{\raisebox{8pt}{\tiny$\bullet$}}\ru\ld\\
    \step[1]\id\step[1]\O{\medA}\step\id\\
    \step[.5]\hld\object{\raisebox{8pt}{\tiny$\bullet$}}\ru\step[1]\id\\
    \step[.5]\id\step[.5]\x\\[-1pt]
    \step[.5]\hlu\object{\raisebox{13pt}{\tiny$\bullet$}}\step[2]\id
  \end{tangles}
  \ \ \ = \ \
  \begin{tangles}{l}
    \vstr{80}\step[1.5]\id\step[1]\ld\\
    \step[1.5]\object{\raisebox{8pt}{\tiny$\bullet$}}\ru\ld\\
    \vstr{80}\step[1]\hdd\step[1]\O{\medA}\step\id\\
    \step[.5]\hld\step[1]\hcd\step[.5]\id\\
    \vstr{70}\ddh\step[.5]\hx\step\id\step[.5]\id\\
    \hcu\step\object{\raisebox{8pt}{\tiny$\bullet$}}\ru\ddh\\
    \vstr{70}\step[.5]\hd\step[1]\hx\\
    \step\lu\object{\raisebox{8pt}{\tiny$\bullet$}}\step[1]\hd
  \end{tangles}
  \ \ \ = \ \
  \begin{tangles}{l}
    \vstr{80}\step[1]\id\step[1]\ld\\
    \step[1]\object{\raisebox{8pt}{\tiny$\bullet$}}\ru\ld\\
    \vstr{80}\ld\step\O{\medA}\step\id\\
    \id\step\object{\raisebox{8pt}{\tiny$\bullet$}}\ru\ld\\
    \vstr{80}\id\step\id\step\O{\medA}\step\id\\
    \vstr{50}\id\step\hx\step\id\\
    \vstr{90}\hcu\step\hx\\
    \vstr{30}\step[.5]\id\step\hdd\step\id\\
    \step[.5]\lu\object{\raisebox{8pt}{\tiny$\bullet$}}\step[1.5]\id
  \end{tangles}
  \ \ \ = \ \
  \begin{tangles}{l}
    \vstr{80}\step[1]\id\step[1]\ld\\
    \step[1]\object{\raisebox{8pt}{\tiny$\bullet$}}\ru\ld\\
    \vstr{80}\ld\step\O{\medA}\step\id\\
    \id\step\object{\raisebox{8pt}{\tiny$\bullet$}}\ru\ld\\
    \vstr{80}\id\step\id\step\O{\medA}\step\id\\
    \id\step\id\step\lu\object{\raisebox{8pt}{\tiny$\bullet$}}\\
    \vstr{67}\id\step\x\\
    \lu\object{\raisebox{8pt}{\tiny$\bullet$}}\step[2]\id
  \end{tangles}
\end{equation*}

\begin{figure}[tbp]
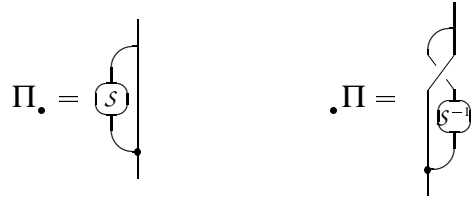

  \centering
  \begin{equation*}
    \Pi_{_{\bullet}} = \ \ \
    \begin{tangles}{l}
      \ld\\
      \O{\medA}\step\id\\
      \lu\object{\raisebox{8pt}{\tiny$\bullet$}}\\
    \end{tangles}
    \qquad\qquad\qquad
    {}_{_{\bullet}}\Pi = \ \ \
    \begin{tangles}{l}
    \ld\\
    \vstr{67}\hxx\\
    \id\step\O{\medA^{\scriptscriptstyle-1}}\\
    \object{\raisebox{8pt}{\tiny$\bullet$}}\ru
  \end{tangles}
  \end{equation*}
  \caption{\small Two projectors on left coinvariants.}
  \label{fig:thePi}
\end{figure}

\noindent
(we used the action and coaction associativity and coassociativity to
isolate and then kill a ``bubble'' as in~\eqref{antipode-axiom} in the
third diagram, and then insert such a bubble into the fourth).
Everything below the $\mapI$-map at the top of the last diagram is
therefore the desired braiding:
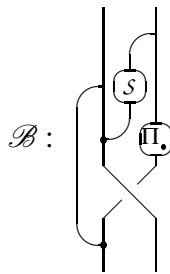
\begin{equation*}
  \Fbraid:\ \
  \begin{tangles}{l}
    \step[1]\id\step\ld\\
    \ld\step\O{\medA}\step\id\\
    \id\step\object{\raisebox{8pt}{\tiny$\bullet$}}\ru
    \step\O{\Pi_{_{\bullet}}}\\
    \id\step\x\\
    \lu\object{\raisebox{8pt}{\tiny$\bullet$}}\step[2]\id
  \end{tangles}
\end{equation*}

\noindent
where $\Pi_{_{\bullet}}$ is defined in Fig.~\ref{fig:thePi}.
But there is more to it, because $\Pi_{_{\bullet}}$ is a projector
onto the left coinvariants in each $\YDspace$:
\begin{equation*}
  \Pi_{_{\bullet}}(s;Y)=
  \begin{cases}
    0,&s\geq 1,\\
    (0;Y),&s=0
  \end{cases}
\end{equation*}
(by the nature of the coaction and of the $\FromLeftii{}$ action, this
property is just the one in~\eqref{antipode-axiom}).  It hence follows
that on graded components, the braiding is given by
\begin{equation}\label{Fbraiding}
  \Fbraid_{s,t}=
  \sum_{i=0}^{s}
  \Shift{i}{\Bbraid_{s - i + t, 1}} \xShift{}{\Bbin{\shift i}{s - i, t}}
  \Shift{s}{\Bbraid_{1, t + 1}} \Shift{(s + 1)}{\A_{t}}
  :
  (s;Y;t;Z)\to\bigoplus_{i=0}^{s}(i;Z;s-i+t;Y).
\end{equation}
Finally, because $\Bbraid_{1, t}\Shift{1}{\A_{t}} =-\A_{t + 1}$, it is
not difficult to see that the braiding can be equivalently written as
stated in the proposition.\qed

\subsubsection{}
The inverse braiding is found similarly, with the result
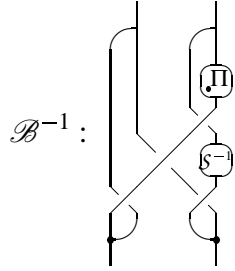
\begin{equation*}
  \Fbraid^{-1}:\ \
  \begin{tangles}{l}
    \ld\step[2]\ld\\
    \id\step\id\step[2]\id\step\O{{}_{_{\bullet}}\!\!\Pi}\\
    \vstr{50}\id\step\id\step[2]\hxx\step\\
    \id\step\xx\step\O{\medA^{\scriptscriptstyle-1}}\\
    \vstr{50}\hxx\step[2]\hxx\\
    \object{\raisebox{8pt}{\tiny$\bullet$}}\ru
    \step[2]\lu\object{\raisebox{8pt}{\tiny$\bullet$}}
  \end{tangles}
\end{equation*}

\noindent
where ${}_{_{\bullet}}\Pi$ is also a projector on left coinvariants,
defined in Fig.~\ref{fig:thePi}

\subsubsection{Squared braiding}
Standard diagram manipulations allow calculating the squared braiding
$\Fbraid^2$.  In terms of components, the $\Fbraid^2$ operation is
by definition the map
\begin{equation*}
  \Fbraid^2\equiv
  \sum_{i=0}^{s} \Fbraid_{i, s - i + t}\Fbraid_{i, s, t}:
  (s;Y;t;Z)\to\bigoplus_{i=0}^{s}(i;Y;s-i+t;Z),
\end{equation*}
where $\Fbraid_{i, s, t}$ is the summand in~\eqref{Fbraiding}.  The
calculation standardly uses the axioms relating the ``dotted'' actions
to left coaction and also the projector property of
$\Pi_{_{\bullet}}$, plus the identity
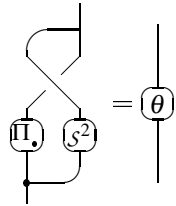
\begin{equation*}
  \begin{tangles}{l}
    \ld[2]\\
    \x\\
    \O{\Pi_{_{\bullet}}}\step[2]\O{\medA^2}\\
    \vstr{80}\object{\raisebox{6pt}{\tiny$\bullet$}}\ru[2]
  \end{tangles}
  \ \ \ = \ \
  \begin{tangles}{l}
    \id\\
    \O{\mbox{\footnotesize$\theta$}}\\
    \id
  \end{tangles}
\end{equation*}

\noindent
where $\theta$ is defined in~\eqref{theta} (and the identity itself is
a reformulation of the one in~\bref{s2-eval}).  The result of diagram
evaluation is
\begin{equation}\label{B2}
  \Fbraid^2:\
  \begin{tangles}{l}
    \step[1.5]\ld\step[2]\O{\mbox{\footnotesize$\theta$}}\\
    \step[1]\hdd\step\O{\Pi_{_{\bullet}}}\step[2]\id\\
    \step[1]\id\step[1.5]\x\\
    \step[.5]\hcd\step[1]\x\\
    \vstr{50}\step[.5]\id\step[1]\hx\step[2]\id\\
    \vstr{80}\step[.5]\lu\object{\raisebox{6pt}{\tiny$\bullet$}}
    \step[1]\lu[2]\object{\raisebox{6pt}{\kern-4pt\tiny$\blacktriangleright$}}
  \end{tangles}
\end{equation}

\noindent
Componentwise, if the two input strands in the diagram represent
$(s;Y)$ and $(t;Z)$, then the $\theta$ map here is exactly
$\Shift{(s+1)}{\theta_{t+1}}$, the full twist of \textit{all} the
$t+1$ factors in $(t;Z)=X^{\otimes t}\tensor Z$.  Furthermore, because
$\Pi_{_{\bullet}}$ projects onto left coinvariants, we can reformulate
the last diagram
as the following proposition.

\begin{prop}\label{FB2}
  $\Fbraid^2$ restricted to $(s;Y;t;Z)$ is the map
  \begin{equation*}
    \sum_{i=0}^{s}\Shift{(i + 1)}{\Adja_{s - i, t}}\Shift{i}{\Bbraid_{s - i, 1}}
    \Shift{s}{\theta_{1+t+1}
    } :
    (s;Y;t;Z)\to\bigoplus_{i=0}^{s}(i;Y;s-i+t;Z)
  \end{equation*}
  where $\Shift{s}{\theta_{1+t+1}} =\id^{\tensor
    s}\tensor\theta_{t+2}: (s)\tensor(Y;t;Z)\to(s)\tensor(Y;t;Z)$ is
  the full twist acting \textit{on the last $t+2$ tensor factors}.
\end{prop}

The formula is worth being restated in words.  To apply the squared
braiding to $(s;Y;t;Z)$, we write this as $(s)\tensor(Y;t;Z)$ and act
with the full twist on the ``ribbon'' $(Y;t;Z)$ whose ``edges'' are
the two ``vertex'' linear spaces.  Next, we deconcatenate $(s)$ into
$(i)\tensor(s-i)$ and then take $(s-i)$ past $Y$ via braiding and
evaluate the left adjoint action $(s-i)\adjoint(t;Z)$.

\begin{Rem}
  In the general position assumed in Secs.~\eqref{sec:start-start}
  and~\eqref{sec:H}, no two punctures $w_1$, \dots, $w_m$ lie on the
  same line from the family.  Degenerate cases where two punctures
  happen to be on the same line (``two fixed lines collide'') can be
  studied in terms of another homology complex, based on a
  stratification of the configuration space of the $w_1,\dots,w_m$,
  with
  the $m$ punctures now ``movable.''
  Similarly to what we did above, the differential of the homology
  complex with coefficients in the corresponding local system then
  yields an associative multiplication.  The full-fledged geometric
  construction is not given here, but the fusion product can be
  regarded as its ``algebraic counterpart.''  On the algebraic side,
  in terms of the corresponding ``large'' Nichols algebra $\Nich(Y)$,
  we note that the case of a subspace $X\subset Y$ such that
  $\Psi(X\tensor X)=X\tensor X$ was considered in~\cite{[G-free]},
  where, under some ``nondegeneracy'' assumptions, it was shown that a
  subalgebra $K$ of $\Nich(Y)$ exists such that
  $\Nich(Y)=K\tensor\Nich(X)$ as a left $K$-module and right
  $\Nich(X)$-module.  Specifically, $K$ is the kernel of braided
  derivations $\partial_{\alpha}$ for $\alpha\in X^*\subset Y^*$ (with
  the last inclusion defined in terms of dual bases), and
  $\Nich(Y)=K\tensor\Nich(X)$ as a left $K$-module and right
  $\Nich(X)$-module. This description of $\Nich(Y)$ as
  $K\tensor\Nich(X)$ is very appealing from the CFT{} standpoint: the
  algebra of screenings has to be tensored just with pure
  (``unscreened'') vertex operators.
\end{Rem}

\section{One-dimensional $X$ and the $(p,1)$ logarithmic CFT models}
\label{sec:1p}
CFT{} applications are associated with \textit{diagonal braiding},
which means that a basis $(x_i)$ exists in the braided linear
space~$X$ such that
\begin{equation}\label{diag-braiding}
  \Psi:x_i\tensor x_j\mapsto q_{ij}\,x_j\tensor x_i
\end{equation}
for some numbers $q_{ij}$.  The braiding involving the $Y$ space is
also assumed diagonal (in addition to the $q_{ij}$, it is then defined
by $q_{i\beta}$, $q_{\alpha j}$, and $q_{\alpha \beta}$, where
$\alpha$ labels basis elements of~$Y$).  In
Appendix~\bref{sec:YD-and-diag}, we briefly describe the context in
which diagonal braiding occurs in the theory of Nichols algebras.  It
has been in the focus of recent activity, with profound results
(see~\cite{[Heck-class],[Ag-all]} and the references therein).

We specialize the constructions in Secs.~\ref{sec:H}--\ref{sec:FUSION}
to the simple(st) case of a rank-$1$ Nichols algebra that in the CFT{}
language corresponds to the $(p,1)$ logarithmic conformal
models~\cite{[Kausch],[Gaberdiel-K]} (also
see~\cite{[FHST],[FGST],[AM-3],[NT]} and the references therein).  We
thus take $X$ to be the one-dimensional linear space spanned by a
screening operator; the $Y$ space is then spanned by highest-weight
vectors of representations of the lattice vertex-operator algebra
associated with the 1-dimensional lattice $\sqrt{2p}\,\oZ$. \
We fix an integer $p\geq 2$ and set  
\begin{equation}\label{root-1}
  \q=e^{\frac{i\pi}{p}}.
\end{equation}

\subsection{CFT background}
\subsubsection{Lattice vertex-operator algebra}\label{sec:lvoa}
Let $\varphi(z)$ be a free chiral scalar field with the basic operator
product expansion
\begin{equation*}
  \varphi(z)\varphi(w) = \log (z-w).
\end{equation*}
The space of observables $\cF$ is the direct sum
$\cF=\bigoplus_{j\in\oZ}\cF_j$ of Fock modules $\cF_j$, $j\in\oZ$,
generated by $\partial\varphi(z)$ from the vertices
\begin{equation}\label{Y-basis-general}
  V^j(z)\equiv
  V^j_{0,0}(z)= e^{\frac{j}{\sqrt{2p}}\varphi(z)},\qquad
  j\in\oZ
\end{equation}
(the ``double zero'' notation is redundant at this point, but is
introduced for uniformity in what follows).  A vector from $\cF_j$ has
the form $P(\partial\varphi(z))V^j(z)$ for a a differential
polynomial~$P$.

The maximum local algebra acting in $\cF$ is the lattice
vertex-operator algebra $\lvoa_{\!\!\sqrt{2p}}$ associated with the
lattice $\sqrt{2p}\,\oZ$.  The space $\cF$ decomposes into a direct
sum of $2p$ simple $\lvoa_{\!\!\sqrt{2p}}$-modules
\begin{equation*}
  \cF=\bigoplus_{r\in\oZ_{2p}}\modVF_r
\end{equation*}
where $\modVF_r=\bigoplus_{j\equiv r\,\mbox{\tiny\rm mod}\,2p}\cF_j$.
The vacuum $\lvoa_{\!\!\sqrt{2p}}$-module is generated from $V^0=1$.


Exchange relations for the vertex operators from~$\cF$ are
\begin{equation}\label{exchange-F}
  V^j(z)
  V^{k}(w)
  =
  \q^{\frac{jk}{2}}\,
  V^{k}(w)
  V^j(z),
  \quad z>w,
\end{equation}
where the condition $z>w$ is to be used for $z$ and $w$ lying on the
same line from the family introduced in~\bref{cells:1}, and is
understood in the sense of the orientation chosen on the lines.

\subsubsection{The $(p,1)$ model and the triplet algebra}
The $(p,1)$ logarithmic CFT{} model is based on the triplet vertex
operator algebra originating in~\cite{[Kausch],[Gaberdiel-K]}.  This
algebra can be constructed by taking the kernel of the ``short''
screening operator
in the vacuum representation of~$\lvoa_{\!\!\sqrt{2p}}$~\cite{[FHST]}.
The algebra is generated by three currents $W^-(z)$, $W^0(z)$, and
$W^+(z)$ given by
\begin{equation}\label{the-W}
  W^-(z)= e^{-\sqrt{2p}\varphi(z)},\qquad
  W^0(z) = [S_+,W^-(z)],\qquad
  W^+(z) = [S_+,W^0(z)],
\end{equation}
where $S_+=\int e^{\sqrt{2 p}\,\varphi}$ (the ``long'' screening
operator).  The associated energy--momentum tensor
\begin{equation*}
  T(z)=\ffrac{1}{2}\partial\varphi(z)\partial\varphi(z)
  +\ffrac{p-1}{\sqrt{2p}}\,\partial^2\varphi(z),
\end{equation*}
encodes the Virasoro algebra with central charge
$c=13-6p-\frac{6}{p}$.

The main source of Hopf-algebra structures associated with the model
is the ``short'' screening operator
\begin{equation}\label{thev1}
  \thev{1}=\int\limits_{-\infty<z<\infty}{
    f(z)
  }, \qquad   f(z)=V^{-2}_{0,0}(z)=e^{\alphaminus\varphi(z)}.
\end{equation}
We let $X$ be the $1$-dimensional space with
basis~\eqref{thev1}.

\subsection{$\cH(X)$ and $\Nich(X)$}
The braided Hopf algebra $\cH=\cH(X)$ is the span of the iterated
screenings
\begin{equation}
  \thev{r}=\kern-10pt\mint{-\infty<z_1<\dots<z_r<\infty}{\kern-10pt
    s(z_1)\dots s(z_r)
  },\qquad r\geq1.
\end{equation}
It follows from~\eqref{exchange-F} that
\begin{equation}\label{FFbraid}
  \Psi(\thev{r}\otimes \thev{s})
  =\q^{2 r s} \thev{s}\otimes \thev{r}
\end{equation}
and
\begin{equation*}
  \thev{r}\dotact \thev{s} = \Abin{r+s}{r}\,\thev{r+s},
\end{equation*}
where the $q$-binomial coefficients defined as
\begin{equation}\label{Abin}
  \Abin{r}{s}=\ffrac{\Afac{r}}{\Afac{s}\Afac{r-s}},
  \quad\Afac{r}=\Aint{1}\dots\Aint{r},
  \quad\Aint{r}=\ffrac{q^{2r}-1}{q^2-1}
\end{equation}
are specialized to $q=\q$.  In particular, $\Aint{p}=0$.  The
deconcatenation coproduct is
\begin{align*}
  \Delta(\thev{r})&=\sum_{s=0}^r \thev{s}\otimes \thev{r-s},
  \\
  \intertext{and the antipode follows by
    evaluating~\eqref{antipode-def} for the above braiding:}
  \A(\thev{r}) &=(-1)^r\q^{r(r-1)}\thev{r}.
\end{align*}
The counit is merely $\epsilon(\thev{r})=\delta_{r,0}1$.

The Nichols algebra of $(X,\Psi)$\,---\,the subalgebra $\Nich(X)$ in
$\cH(X)$ generated by $\thev{1}$\,---\,is therefore the linear span of
the $\thev{r}$ with $0\leq r\leq p-1$.

We also write $F=\thev{1}$; then $F^p=0$ and the algebra $\Nich$ has
the basis $F^i$, $0\leq i\leq p-1$.

\subsubsection{The $Y$ space}
We next specify the $Y$ space, the source of Yetter--Drinfeld modules
over the Nichols algebra.

Setting $V^{j}=V^{j}(0)$, we define the braided linear space
\begin{equation*}
  Y= \text{span}(V^j\mid j\in\oZ_{4p})
\end{equation*}
where the braiding is
\begin{equation}\label{braid-simplest}
  \Psi(V^j\tensor V^{k}) = \q^{\frac{jk}{2}}  V^{k}\tensor V^j.
\end{equation}
Also, the braiding of (elements of) $X$ and $Y$ is
\begin{equation}\label{SV}
  \Psi(V^j\tensor\thev{r})=\q^{-j r}\thev{r}\tensor V^j,
  \quad
  \Psi(\thev{r}\tensor V^j)=\q^{-j r}\,V^j\tensor\thev{r}.
\end{equation}


\subsection{Hopf bimodules and Yetter--Drinfeld $\Nich(X)$-modules}
\subsubsection{Hopf bimodules}
For diagonal braiding, the Hopf bimodule $\Yspace$
(see~\bref{sec:Hopf-bimodule}) decomposes as
$\Yspace=\bigoplus_{i}(\!(y_i)\!)$, where the sum is taken over a
basis in $Y$ that diagonalizes the braiding.  Hence, a Hopf
$\Nich(X)$-bimodule $(\!(V^j)\!)$ is generated from each of the $4p$
chosen basis elements of $Y$.  For each $j=1,\dots,4p$, the basis in
$(\!(V^j)\!)$ is
\begin{multline}\label{dressed-vert}
  V^j_{r,t} \equiv (F(r);\, V^j;\, F(t))=
  \\
  \mint{-\infty<x_1<\dots<x_r<0<x_{r+1}<\dots<x_{r+t}<\infty}{
    \kern-50pt s(x_1)\dots s(x_r) V^j(0) s(x_{r+1})\dots s(x_{r+t})},
  \quad 0\leq r,t\leq p-1.
\end{multline}
In this basis, the left and right $\Nich(X)$-actions 
(see~\eqref{r.sYt} and~\eqref{sYt.r}) become
\begin{align*}
  F\dotact V^j_{r,t}(z)&=\aint{r+1}V^j_{r+1,t}(z)
  +\q^{2 r-j}\aint{t+1}V^j_{r,t+1}(z),\\
  V^j_{r,t}(z)\dotact F&=\q^{2 t-j}\aint{r+1}V^j_{r+1,t}(z)
  +\aint{t+1}V^j_{r,t+1}(z).
\end{align*}

The left and right coactions on~$(\!(V^j)\!)$ are given, as usual, by
deconcatenation of the integrals in~\eqref{dressed-vert}:
\begin{align}
  \label{Ldelta}
  \delta_{\mathrm{L}} V^j_{s,t} =\sum_{j=0}^s F(j)\otimes
  V^j_{s-j,t},\qquad
  \delta_{\mathrm{R}} V^j_{s,t} =\sum_{j=0}^t V^j_{s,t-j}\otimes F(j).
\end{align}

\subsubsection{``Multivertex'' Hopf bimodules}\label{multi-Hopf}
The multivertex Hopf bimodules introduced in \bref{sec:more-Hopf} are
now labeled by $\set{j}=\{j_1,j_2,\dots,j_\ell\}$, where each $j_i$
ranges over $4p$ values.  For each such $\set{j}$, the Hopf bimodule
$(\!(\set{j})\!)$ is the span of the
$V^{\set{j}}_{t_1,t_2,\dots,t_{\ell+1}}$ that are obtained by fixing,
once and for all, some $z_1<z_2<\dots<z_\ell$ in the expressions
\begin{align*}
  V^{\set{j}}_{t_1,t_2,\dots,t_{\ell+1}}(z_1,z_2,\dots,z_{\ell}) &=
  (F(t_1); V^{j_1}(z_1); F(t_2); V^{j_2}(z_2); \dots;F({t_{\ell}});
  V^{j_{\ell}}(z_{\ell});F(t_{\ell+1})),
  \\
  \intertext{where $0\leq t_j\leq p-1$.  These expressions are
    $(t_1+t_2+\dots+t_{\ell+1})$-fold integrals:} &=\kern-30pt
  \idotsint\limits_{\substack{ -\infty<x^1_1<\dots<x^1_{t_1}<
      z_1<x^2_1<
      \dots\\
      \dots
      <x^\ell_{t_\ell}<z_{\ell}<x^{\ell+1}_1<\dots<x^{\ell+1}_{t_{\ell+1}}<\infty
    }}  \kern-15pt
  \prod_{n=1}^\ell\Bigl(\prod_{m=1}^{t_n}\bigl(s(x^n_m)\bigr)
  V^{j_n}(z_n)\Bigr)
  \prod_{m=1}^{t_{\ell+1}}\bigl(s(x^{\ell+1}_m)\bigr)
\end{align*}
with the integrations only over the $x_{t_i}^j$, \textit{not over
  the~$z_i$}.

The structure of the $\Nich(X)$-action on
$V^{\set{j}}_{t_1,t_2,\dots,t_{\ell+1}}\in(\!(\set{j})\!)$ is already
clearly seen in the case of two-vertex modules ($\ell=2$).  The left
and right $\Nich(X)$-actions are
\begin{align}
  \notag F\dotact
  V^{\{j_1,j_2\}}_{n_1,n_2,n_3}={}&\,\aint{n_1+1}\,V^{\{j_1,j_2\}}_{n_1+1,n_2,n_3}
  +\q^{2n_1-j_1}\aint{n_2+1}V^{\{j_1,j_2\}}_{n_1,n_2+1,n_3}
  \\
  \notag
  {}&+\q^{2n_1+2n_2-j_1-j_2}\aint{n_3+1}V^{\{j_1,j_2\}}_{n_1,n_2,n_3+1},
  \\
  \notag V^{\{j_1,j_2\}}_{n_1,n_2,n_3}\dotact F={}&
  \q^{2n_2+2n_3-j_1-j_2}\aint{n_1+1}V^{\{j_1,j_2\}}_{n_1+1,n_2,n_3}
  +\q^{2n_3-j_2}\aint{n_2+1}V^{\{j_1,j_2\}}_{n_1,n_2+1,n_3}
  \\
  \notag {}&+\aint{n_3+1}V^{\{j_1,j_2\}}_{n_1,n_2,n_3+1},
\end{align}

\subsubsection{Yetter--Drinfeld $\Nich(X)$-modules}\label{sec:YDp1}
As in~\bref{sec:YD}, we take the space of right coinvariants in
$(\!(V^j)\!)$, that is,
\begin{equation*}
  \modV^j=\text{span}(V^j_{s,0}\mid 0\leq s\leq p-1).
\end{equation*}
It is a (left--left) Yetter--Drinfeld $\Nich(X)$-module under the left
coaction~\eqref{Ldelta} and the left adjoint action described
in~\bref{prop-adja}, which now evaluates as
\begin{align}\label{F(r)-act}
  F(r)\adjoint\VertexII{s}{j} &= \abin{r + s}{r} \prod_{a=s}^{s + r -
    1}\bigl(1 - \q^{2 a - 2 j}\bigr)\, \VertexII{r + s}{j}.
\end{align}

\subsection{
  The module nomenclature}\label{sec:nomencl}
The dependence of the right-hand side of~\eqref{exchange-F} (and
hence~\eqref{braid-simplest}) on $j$ and $k$ is through
$\q^{\frac{jk}{2}}$, and is therefore $4p$-periodic.  By contrast, it
is clear from~\eqref{F(r)-act} that the module structure on $\modV^j$
depends on $j$ modulo $p$ (and in fact so does the comodule
structure).

Hence, if the braiding structure is forgotten, then there are only $p$
nonisomorphic $\Nich(X)$ module comodules among the $\modV^j$ (e.g.,
those with $j=0,\dots,p-1$).  If not only the module comodule
structure but also the braiding
$\Psi:\Nich(X)\tensor\modV^j\to\modV^j\tensor\Nich(X)$ is considered,
then there are $2p$ nonisomorphic objects among the $\modV^j$.  And
finally, if these are considered as objects in the braided category of
Yetter--Drinfeld $\Nich(X)$-modules, then $4p$ of them are
nonisomorphic.

It is therefore convenient to parameterize the ``momenta'' $j$ as
\begin{equation}\label{nu-parametrization}
  j \equiv r-1-\nu
  p\,\;\mathrm{mod}\;\,4p,\quad
  r=1,\dots p, \quad \nu\in\oZ_4=\{0,1,2,3\}
\end{equation}
to single out the range $r=1,\dots,p$ that is only ``seen'' by the
representation theory with no braiding considered, and set
\begin{equation*}
   \modV_{r}(\nu)
   =\modV^{r-1-\nu p},
  \qquad r=1,\dots,p,\quad \nu\in\oZ_4.
\end{equation*}

The following statements are verified directly, using \eqref{F(r)-act}
and \eqref{Ldelta}.
\begin{enumerate}

\item \textit{The Yetter--Drinfeld modules $\modV_p(\nu)$,
    $\nu\in\oZ_4$, are simple}.

\item \textit{For $1\leq r\leq p-1$, the Yetter--Drinfeld module
    $\modV_r(\nu)$ contains a simple $r$-dimensional Yetter--Drinfeld
    submodule}
  \begin{equation*}
    \modX_r(\nu)=\text{span}(V^{r-1-\nu p}_{s,0}\mid 0\leq s\leq r-1).
  \end{equation*}

\item \textit{The quotient $\modV_{r}(\nu)/\modX_r(\nu)$ is isomorphic
    to $\modX_{p-r}(\nu+1)$ as an object in the braided category of
    Yetter--Drinfeld $\Nich(X)$-modules}.
\end{enumerate}

\subsection{Fusion product and related structures}
Formula \eqref{YDfusion} for the fusion product of Yetter--Drinfeld
modules takes the form
\begin{equation*}
  V^{j}_{r,0}\odot
  V^{m}_{s,0}
  = \sum_{i=0}^{s}
  \q^{-i j}
  \abin{r+i}{r}
  \VertexIV{r+i}{j}{s-i}{m}.
\end{equation*}
The $\VertexIV{r_1+i}{j}{s_1-i}{m}$ occurring in the right-hand side
are elements of the two-vertex Yetter--Drinfeld modules
(see~\bref{multi-YD}).  The left adjoint action on these is most
easily calculated from the factorization formula~\eqref{eq:YD2}: it
follows that
\begin{multline}\label{Ad-F(r)}
  F(r)\adjoint
  \VertexIV{t_1}{j_1}{t_2}{j_2} = \sum_{s=0}^{r}\q^{s (2 t_1 - j_1)}
  \abin{t_1 + r - s}{r - s} \abin{t_2 + s}{s}
  \\[-4pt]
  {}\times\prod_{a=s + 1}^{r}\bigl(1 - \q^{2 (t_1 + 2 t_2 + a - 1 -
    j_1 - j_2)}\bigr) \prod_{b=0}^{s - 1}\bigl(1 - \q^{2 (t_2 + b -
    j_2)}\bigr)\, \VertexIV{t_1 + r - s}{j_1}{t_2 + s}{j_2}.
\end{multline}
Calculations with these ingredients support the following conjecture.

\subsubsection{Conjecture}
The fusion products of the $\modX_r(\nu)$, $1\leq r\leq p$,
$\nu\in\oZ_4$, decompose as
\begin{equation}\label{fusionXX}
  \modX_r(\nu)\odot\modX_s(\mu)=
  \bigoplus_{j=|r-s|+1}^{{\mathrm{min}}(r+s-1,2p-r-s-1)}\modX_j(\nu+\mu)
  \,\,\oplus\,\,\bigoplus_{j = 2p-r-s+1}^{p}\modP_j(\nu+\mu),
\end{equation}
where both summations are with step~$2$ (and start at the lower limit,
which may be essential for the second sum), and the $\modP_j(\nu)$
modules are as follows: $\modP_p(\nu)=\modX_p(\nu)$, and each
$\modP_s(\nu)$ with $s=1,2,\dots p-1$
is an indecomposable Yetter--Drinfeld $\Nich(X)$-module with the
structure of subquotients given by
\begin{equation}\label{schem-proj}
  \xymatrix@=12pt{
    &&\modX_{s}(\nu)
    \ar[dl]
    \ar[dr]
    &\\
    &\modX_{p - s}(\nu-1)\ar[dr]
    &
    &\modX_{p - s}(\nu+1)\ar[dl]
    \\
    &&\modX_{s}(\nu)&
  }
\end{equation}

We illustrate this with examples.

\subsubsection{Examples}
We restrict ourself to $p=3$ and consider
$\modX_{3}(0)\odot\modX_{3}(0)$ and $\modX_{2}(0)\odot\modX_{3}(0)$.

\medskip

\noindent
\textit{\thesubsubsection.1. \ $\modX_{3}(0)\odot\modX_{3}(0)$.}  For
$p=3$, we choose $V^{\{2,2\}}_{r,s,0}$ with $0\leq r,s\leq2$ as a
basis in $\modX_{3}(0)\odot\modX_{3}(0)$.  Using (\ref{Ad-F(r)}) and
(\ref{Ldelta}), we then verify that the subspace $\modP_1(0)$ spanned
by the vectors
\begin{equation}\label{Proj1}
  V^{\{2,2\}}_{0,0,0},\quad V^{\{2,2\}}_{1,0,0}-V^{\{2,2\}}_{0,1,0},
  \quad \begin{array}{c}
    \displaystyle V^{\{2,2\}}_{2,0,0}-V^{\{2,2\}}_{1,1,0},\\[10pt]
    \displaystyle V^{\{2,2\}}_{0,2,0},
  \end{array}
  \quad
  V^{\{2,2\}}_{1,2,0},\quad V^{\{2,2\}}_{2,2,0}
\end{equation}
is invariant under action (\ref{Ad-F(r)}) and the left coaction in
(\ref{Ldelta}), and is therefore a Yetter--Drinfeld submodule in
$\modX_{3}(0)\odot\modX_{3}(0)$.  Moreover, the vectors are arranged
in accordance with the structure claimed in~\eqref{schem-proj}.  For
example, coaction applied to the top vector yields
$1\otimes(V^{\{2,2\}}_{2,0,0}-V^{\{2,2\}}_{1,1,0})
+F(1)\otimes(V^{\{2,2\}}_{1,0,0}-V^{\{2,2\}}_{0,1,0}) +F(2)\otimes
V^{\{2,2\}}_{0,0,0}$, which is an element of
$\Nich(X)\tensor\modP_1(0)$.\footnote{To check that the parameters
  $\nu$ of subquotients in \eqref{Proj1} are as shown in
  \eqref{schem-proj}, we represent the total momenta of the relevant
  vertices in form \eqref{nu-parametrization}.  Because the elements
  of the module are related by either the action or the coaction, it
  is straightforward to see that the $\nu$-pattern is always as shown
  in~\eqref{schem-proj}.  It follows similarly that the $\nu$
  parameter is always additive in~\eqref{fusionXX}.}

Another Yetter--Drinfeld submodule in $\modX_{3}(0)\odot\modX_{3}(0)$
is the span of
\begin{equation}\label{steinberg}
  V^{\{2,2\}}_{0,1,0},\quad V^{\{2,2\}}_{1,1,0}+V^{\{2,2\}}_{0,2,0},
  \quad V^{\{2,2\}}_{2,1,0}+V^{\{2,2\}}_{1,2,0}
\end{equation}
Using (\ref{Ad-F(r)}) and (\ref{Ldelta}), it is straightforward to
verify that this last subspace is in fact isomorphic to the
Yetter--Drinfeld module $\modX_{3}(0)$.  Hence,
\begin{equation*}
  \modX_{3}(0)\odot\modX_{3}(0)=\modP_1(0)\oplus\modX_3(0),
\end{equation*}
which is~\eqref{fusionXX} in this particular case.

\medskip

\noindent
\textit{\thesubsubsection.2. \ $\modX_{2}(0)\odot\modX_{3}(0)$}.  This
space has the basis $V^{\{2,2\}}_{r,s,0}$, with $0\leq r\leq1$ and
$0\leq s\leq2$.  We suggestively arrange these vectors as
\begin{equation}\label{proj2}
  V^{\{1,2\}}_{0,0,0},\quad
  \begin{array}{lr}
    \displaystyle
    V^{\{1,2\}}_{1,0,0},&\aint{2}V^{\{1,2\}}_{2,0,0}+\q
    V^{\{1,2\}}_{1,1,0},\\[6pt]
    \displaystyle
    V^{\{1,2\}}_{0,1,0},&\q V^{\{1,2\}}_{1,1,0}+\aint{2} V^{\{1,2\}}_{0,2,0},
  \end{array}
  \quad
  \q V^{\{1,2\}}_{2,1,0}+\aint{2} V^{\{1,2\}}_{1,2,0}
\end{equation}
and then use (\ref{Ad-F(r)}) and (\ref{Ldelta}) to verify that this is
an indecomposable Yetter--Drinfeld module, $\modP_2(0)$ (again, with
the structure following the pattern in~\eqref{schem-proj}), and hence
\begin{equation*}
  \modX_{2}(0)\odot\modX_{3}(0)=\modP_2(0),
\end{equation*}
again in accordance with~\eqref{fusionXX}.

\subsection{$\odot$ is not commutative}
We next show that if the $\Nich(X)$ modules constructed above are
considered not as objects in the braided category of Yetter--Drinfeld
modules but as only $\Nich(X)$ module comodules with the braiding
$\Psi:\Nich(X)\tensor\modV^j\to\modV^j\tensor\Nich(X)$ (and with the
corresponding morphisms; see~\bref{sec:nomencl}), then the $\odot$
product is not commutative, i.e., the product of two modules taken in
two different orders may result in nonisomorphic modules.  As
in~\cite{[KoSa]}, showing this requires modules of the
``$\modO$''~type, labeled by a parameter $z\in\oC\!\mathbb{P}^1$.
Such a $\Nich(X)$ module comodule, denoted by $\modO_2(0)(1,z)$, is
readily identified as a sub(co)module in~\eqref{Proj1}: for each
$z=z_1:z_2$, the relevant subspace is the span of
\begin{gather*}
  t_1=z_1V^{\{2,2\}}_{0,0,0}+z_2V^{\{2,2\}}_{1,2,0},\qquad
  t_2=z_1(V^{\{2,2\}}_{1,0,0}-V^{\{2,2\}}_{0,1,0})+z_2V^{\{2,2\}}_{2,2,0},\\
  b=V^{\{2,2\}}_{0,2,0}.\qquad
\end{gather*}
That this is a submodule subcomodule follows from the action and
coaction formulas
\begin{gather*}
  F\adjoint t_1 = (1-\q^{-2})t_2,\qquad
  F\adjoint t_2 = -z_1(1-\q^{2})b,\qquad F\adjoint b=0,\\
  \delta_{\mathrm{L}} t_1 = 1\otimes t_1+z_2F(1)\otimes b,\quad
  \delta_{\mathrm{L}} t_2 = 1\otimes t_2+F(1)\otimes
  t_1+z_2F(2)\otimes b,\quad \delta_{\mathrm{L}} b = 1\otimes b.
\end{gather*}

We now calculate $\modX_{1}(1)\odot\modO_2(0)(1,z)$ and
$\modO_2(0)(1,z)\odot\modX_{1}(1)$.  From the action of $F$ on
three-vertex objects
\begin{multline}\label{F-on-3}
  F\adjoint \VertexVI{t_1}{a}{t_2}{b}{t_3}{c} = (1 - \q^{2 t_1 + 4 t_2
    + 4 t_3 - 2 a - 2 b - 2 c}) \qInt{t_1 + 1} \VertexVI{t_1 +
    1}{a}{t_2}{b}{t_3}{c}
  \\
  \qquad\qquad\qquad
  {}+ \q^{2 t_1 - a} (1 - \q^{2 t_2 + 4 t_3 - 2 b - 2 c}) \qInt{t_2 + 1}
  \VertexVI{t_1}{a}{t_2 + 1}{b}{t_3}{c}
  \\
  + \q^{2 t_1 + 2 t_2 - a - b} (1 - \q^{2 t_3 - 2 c}) \qInt{t_3 + 1}
  \VertexVI{t_1}{a}{t_2}{b}{t_3 + 1}{c},
\end{multline}
it is straightforward to see that
\begin{equation}\label{XO-prod}
 \modX_{1}(1)\odot\modO_2(0)(1,z)=\modO_2(1)(1,-z),\qquad
 \modO_2(0)(1,z)\odot\modX_{1}(1)=\modO_2(1)(1,z).
\end{equation}

The $\odot$ product is therefore not commutative.  Another way to see
this is to consider the braiding $\Fbraid$ from~\bref{prop:FYDB}.  On
two-vertex Yetter--Drinfeld modules, the braiding operation $\Fbraid$
from~\bref{prop:FYDB} becomes
\begin{align*}
 \Fbraid
 \VertexIV{s}{j_1}{t}{j_2} &= (-1)^t\q^{\frac{1}{2}j_1 j_2-t(j_1+j_2)+t(t-1)}
 \sum_{r=0}^s\q^{-j_2r}\abin{t+r}{t}
 \VertexIV{s-r}{j_2}{t+r}{j_1}.
\end{align*}
A simple calculation with~\bref{prop:FYDB} and~\eqref{F-on-3} shows
that
\begin{equation}\label{B-iso}
 \Fbraid:\modX_{1}(1)\odot\modO_2(0)(1,z)
 \to
 \modO_2(0)(1,-z)\odot\modX_{1}(1) 
\end{equation}
is an isomorphism of module comodules.

The braided Hopf algebra $\Nich(X)$ admits an outer automorphism
defined as
\begin{align*}
  \theaut:F(r)&\mapsto(-1)^{r}F(r)\\
  \intertext{which acts on basis elements of $\Nich(X)$-modules as}
  \theaut:V^{\set{j}}_{s_1,s_2,\dots,s_{\ell},0}&\mapsto
  (-1)^{\sum_{m=1}^\ell s_m}V^{\set{j}}_{s_1,s_2,\dots,s_{\ell},0},
\end{align*}
where $\set{j}=\{j_1,j_2,\dots,j_\ell\}$ and each $j_i$ ranges
over~$\oZ_{4p}$.  It is immediate to see that this preserves the
Yetter--Drinfeld condition~\eqref{yd-axiom}, and also that
\begin{gather*}
 \theaut:\modX_r(\nu)\to \modX_r(\nu),\quad\theaut:\modP_r(\nu)\to
 \modP_r(\nu),
 \\
 \intertext{but}
 \theaut:\modO_2(0)(1,z)\to\modO_2(0)(1,-z).
\end{gather*}
We can then write \eqref{B-iso} as
\begin{equation}
 \FYDBraid:\modX_{1}(1)\odot\modO_2(0)(1,z)
 \to
 \theaut\bigl(\modO_2(0)(1,z)\bigr)\odot\theaut\bigl(\modX_{1}(1)
 \bigr).
\end{equation}

The representation category of the triplet $W$-algebra in $(p,1)$
models is equivalent to the representation category of the $\UresSL2$
quantum group \cite{[FGST2],[NT]} as Abelian categories.  But some
indecomposable $\UresSL2$ representations (the ``$\modO$'' modules
in~\cite{[FGST2]}, from which the notation was borrowed above)
demonstrate ``noncommutativity'' (tensor product in two different
orders gives nonisomorphic modules)~\cite{[KoSa]}.  Not surprisingly,
this is ``the same'' noncommutativity as in \eqref{XO-prod}.  The
category of Yetter--Drinfeld $\Nich(X)$-modules with the braiding
among them forgotten (and with the morphisms being module comodule
morphisms commuting with the braiding
$\Psi:\Nich(X)\tensor\modV^j\to\modV^j\tensor\Nich(X)$) is equivalent
to the category of representations of the triplet $W$-algebra as an
Abelian category.


In another form, the relation to $(p,1)$ models can be expressed in
terms of ``categories with monodromy'' (\textit{entwined} categories
in~\cite{[Brug]})---which on the algebraic side is the category of
Yetter--Drinfeld $\Nich(X)$-modules $\modV$ endowed not with braiding
but with the squared braiding (monodromy)
$\Psi^2:\modV\tensor\modW\to\modV\tensor\modW$ (and with the braiding
$\Psi:\Nich(X)\tensor\modV\to\modV\tensor\Nich(X)$).

\subsection{The dual triangular part of the algebra
}
We recall the ``braided double problem'' mentioned in the Introduction
(see p.~\pageref{doubleproblem})---the nonexistence of an immediate
analogue of the Drinfeld double in the braided case.  We have also
noted that constructing an ordinary Hopf algebra via double
bosonization requires a list of conditions
(cf.~\cite{[Sommerh-deformed]}).  For our current $\Nich(X)$ algebra
of a very simple structure, we now show that introducing the dual
generator as in~\bref{sec:H*} at least leads to an associative
algebra.

Considering $\Nich^*$, we introduce $E\in\Nich^*$ as the linear
functional
\begin{equation}\label{def-E}
  \eval{E}{\thev{r}}=\delta_{r,1}.
\end{equation}
A basis in $\Nich^*$ is then $E^n$ with $0\leq n\leq p-1$.

The action of $E$ on Hopf bimodule elements \eqref{dressed-vert}, also
denoted by $\adjoint$ for uniformity and defined by $E\adjoint V =
\eval{E}{V\mone}V\zero$ with coaction (\ref{Ldelta}), is
\begin{equation}\label{act-E}
  E\adjoint V^{j}_{s_1,s_2} =V^{j}_{s_1-1,s_2}.
\end{equation}
Then the ``commutator'' evaluated in~\bref{lemma:K2} becomes
\begin{equation}\label{qcomm}
  \bigl((E\adjoint)(F\adjoint)-\q^2(F\adjoint)(E\adjoint)\bigr)
  V^{j}_{s_1,s_2}
  = \bigl(1-\q^{4s_1+4s_2-2 j}\bigr)V^{j}_{s_1,s_2}.
\end{equation}
Therefore, $\KK$ defined in~\bref{lemma:K2} can be assigned the
following action on elements of the corresponding Hopf bimodules (or,
for $s_2=0$, Yetter--Drinfeld modules):
\begin{equation}\label{U-def}
  \KK\adjoint V^{j}_{s_1,s_2}=
  \q^{4s_1+4s_2-2 j}V^{j}_{s_1,s_2}.
\end{equation}
We then have the $q$-commutator relation (omitting the $\adjoint$)
\begin{equation}\label{r1}
  E F - \q^2 F E  =  1 - \KK.
\end{equation}
Remarkably, it also follows from (\ref{U-def}) that
\begin{equation}\label{r2}
  \KK F = \q^{4}F\KK,\qquad \KK E =\q^{-4}E\KK.
\end{equation}

We thus have an \textit{associative algebra} $D(\Nich)$ generated by
$F$, $\KK$, and $E$.  (In addition, $F^p=0$ as we saw above, $E^p=0$
either from~\eqref{act-E} or from the fact that
$\Nich(X^*)=\Nich(X)^*$, and $\KK^p=1$ from~\eqref{U-def}).

Interestingly, the operator $\KK$ read off from~\eqref{qcomm} turns
out to be ``too coarse-grained'' compared with what we have in the
nonbraided Hopf algebra that is the Kazhdan--Lusztig-dual to the
triplet $W$ algebra---the $\UresSL2$ quantum group~\cite{[FGST]}.  In
fact, extracting a ``root'' $K$ of $\KK$, $K^2=\KK$, allows making
$D(B)$ into the $\UresSL2$ quantum group; for this, we just take $KE$
and $F$ as generators.

The relevant setting under which the occurrence of a nonbraided
quantum group, such as $\UresSL2$, from a braided Hopf algebra, such
as our $\Nich(X)$, would not seem accidental is offered by ``double
bosonization'' in the version described in~\cite{[Sommerh-deformed]}
(also see\cite{[Sommerh-int]}), which generalizes the Radford
biproduct~\cite{[Radford-bos]} (``bosonization''~\cite{[Mj-braided]})
to triple ``smash'' products $A\Smash H\Smash B$, where $A$ and $B$
are Hopf algebras in $\HHyd$ dual to each other and the Hopf algebra
$H$ is commutative and cocommutative.  The relevant (and rather
intricate) details are discussed in~\cite{[Sommerh-deformed]}; a
crucial necessary requirement is that the braiding~\eqref{YD-braiding}
in $\HHyd$ coincide with the braiding $u\tensor v \mapsto
v\zero\tensor v\mone\leftii u$, which is the braiding in
$\YDname^{\scriptscriptstyle H}_{\!\scriptscriptstyle H}$ for the
right--right Yetter--Drinfeld structure $v\rightii h \mathrel{{:}{=}}
\hA^{-1}(h)\leftii v$ and $v\zero\tensor v_{_{(1)}} \mathrel{{:}{=}}
v\zero\tensor \hA(v\mone)$. \ That is, the condition is $u\mone\leftii
v \tensor u\zero = v\zero\tensor v\mone\leftii u$,
which for diagonal braiding means that $q_{ij}=q_{ji}$.

\section{Conclusions}
We have outlined a strategy to construct Hopf algebra counterparts of
vertex-operator algebras in logarithmic CFT{}. \ Part of the original
motivation was to eliminate, or ``resolve,'' the ad hoc elements of
the approach in~\cite{[FGST],[FGST2],[FGST3],[FGSTq]}. \ The ``purest
and minimalist'' scheme appears to be the one based on
\textit{braided} Hopf algebras; relation to ordinary Hopf algebras
should then follow via some form of bosonization.  The idea is
therefore that the ``braided'' approach is to play a more fundamental
role in the duality of logarithmic CFT models with
Hopf~algebras$/$quantum~groups.

In addition to its ``strategic'' significance, the braided ideology
offers several attractive possibilities at a more practical level.  For
example, the square of the braiding $\Fbraid$ calculated in
Sec.~\ref{FB2} gives a powerful tool for constructing the center of
the category: in a rigid category, one of the two lines in~\eqref{B2}
can be ``closed'' using evaluation and coevaluation, with the result
representing a central element action on the module corresponding to
the other line.  We expect to return to this in the future.

As we have noted, the choice of the basic contour system in
Fig.~\ref{fig:intro} can be interpreted as a space--time decomposition
of the complex plane, with the lines interpreted as spatial slices and
a transverse direction as time.  Topologically equivalent systems of
contours give isomorphic Nichols algebras of screenings. It might be
interesting to consider different nonequivalent choices of contour
systems, among which the one that immediately suggests itself is the
system of concentric circles, corresponding to radial quantization
in~CFT{}.

\subsubsection*{Acknowledgments}
We are grateful to B.~Feigin, A.~Gainutdinov, M.~Grigoriev, A.~Rosly,
I.~Runkel, and C.~Schweigert for the very useful discussions.  In
retrospect, the pressure from J.\;Fuchs, indirect but insistent, is
also gratefully acknowledged.  The above notwithstanding, we also
thank the referee for quite a few very useful comments and questions.
AMS was supported in part by the RFBR grant 10-01-00408 and the
RFBR--CNRS grant 09-01-93105.  IYuT is grateful to H.~Saleur for kind
hospitality in IPhT where a part of the work was made. IYuT was
supported in part by the RFBR--CNRS grant 09-02-93106 and the RFBR
grant 08-02-01118.

\appendix
\section{Braids and shuffles}\label{app:shuffles}
We here recall the definitions pertaining to the braid group and
``quantum'' shuffles, and fix our conventions.  In the last subsection
in this appendix, \bref{techn1}, we write a number of quantum shuffle
identities that express the properties of Hopf and Yetter--Drinfeld
modules established in Secs.~\ref{sec:H} and~\ref{sec:FUSION}.

\subsection{``Classical'' shuffles}\label{sec:cl-shuffles}
Let $\pi$ be a partition of $n\in\oN$.  We write
$\pi=(\pi_1,\dots,\pi_k)$, where
$\pi_1\cup\dots\cup\pi_k=[1,\dots,n]$.  The subgroup
$\mathbb{S}^n_\pi\subset\mathbb{S}_n$ of \textit{shuffle permutations}
subordinated to $\pi$ consists of those functions
$\mathbb{S}_n\ni\sigma:[1,\dots,n]\to[1,\dots,n]$ that are
monotonically increasing on each~$\pi_i$.

\subsection{Braids and ``quantum'' shuffles}\label{sec:q-shuffles}
Let $(\catC,\tensor,\Psi)$ be a braided monoidal category.
This means, in particular, that $\Psi$ satisfies the braid equation
(the hexagon equation in disguise, with the associativity morphisms
omitted)
\begin{equation*}
  (\Psi\tensor\id)(\id\tensor\Psi)(\Psi\tensor\id)=
  (\id\tensor\Psi)(\Psi\tensor\id)(\id\tensor\Psi).
\end{equation*}

We recall that the Matsumoto section $\mathbb{S}_n\to\mathbb{B}_n$ is
constructed as follows.  Each $\sigma\in\mathbb{S}_n$ decomposes into
a product of elementary permutations: $\sigma=\prod\tau_{i,i+1}$.
Whenever this product is \textit{of the minimum length},
$\sigma\mapsto\prod\Psi_{i,i+1}$.

For $n$ objects $X_i$ in $(\catC,\tensor,\Psi)$
and a partition $\pi$ of $n$, the (``quantum'') shuffle
$\Bbin{n}{\pi}\in\End(
X_1\tensor\dots\tensor X_n)$ is
\begin{equation}\label{sh-Cat}
  \Bbin{n}{\pi} =
  \Bbin{n}{\pi}[\catC
  ]
  = \sum_{\sigma\in \mathbb{S}^n_\pi}
  \sigma_{\catC},
\end{equation}
where $\sigma\mapsto
\sigma_{\catC}:\mathbb{S}_n\to\mathbb{B}_n\to\End(
X_1\tensor\dots\tensor X_n)$ is the composition of the Matsumoto
section and the braid group representation on
$X_1\tensor\dots\tensor X_n$, and $\ell(\sigma)$ is the reduced length
of a $\sigma\in\mathbb{S}_n$.  We omit $\catC$ from the notation for
the $\Bbin{}{}$.
In fact, we mostly think of the shuffles in~\eqref{sh-Cat} and similar
formulas below in ``universal'' terms, as elements of the braid group
algebra (for braids on a ``sufficiently large'' number of strands,
i.e., for an inductive-limit braid group $\mathbb{B}_\infty
=\lim\limits_{\longrightarrow}\mathbb{B}_n$
).

\subsection{Notation}\label{sec:notation}
In a braided monoidal category $(\catC,\tensor,\Psi)$, we use the
standard notation
\begin{multline*}
  \Psi_i\equiv\id^{\tensor(i-1)}\tensor\Psi\tensor\id^{\tensor(n-i-1)}:
  X_1\tensor\dots
  \tensor X_n \\ 
  \to
  X_1\tensor\dots X_{i-1}\tensor X_{i+1}\tensor X_i
  \tensor X_{i+2}\tensor\dots \tensor X_n,
\end{multline*}
The braid equation is then expressed as
\begin{equation*}
  \Psi_s\Psi_{s+1}\Psi_s = \Psi_{s+1}\Psi_{s}\Psi_{s+1}.
\end{equation*}
The braiding of $X_1\tensor\dots\tensor X_m$ with
$Y_1\tensor\dots\tensor Y_n$
is effected by
\begin{gather*}
  \Bbraid_{m,n}
  =
  (\Psi_{n}\dots\Psi_{1})
  (\Psi_{n+1}\dots\Psi_{2})
  \dots
  (\Psi_{n+m-1}\dots\Psi_{m}).
\end{gather*}

Another useful notation is the shift $\Shift{m}{}$: for any morphism
of the form
$f=\id^{\tensor i}\tensor F\tensor\id^{\tensor j}$, we set
\begin{equation*}
  \Shift{m}{f}=\id^{\tensor(i+m)}\tensor F\tensor\id^{\tensor(j-m)}.
\end{equation*}
The two notational conventions are related via
$\Shift{j}{\Psi_i}=\Psi_{i+j}$.

\subsection{Braided factorial and binomials}\label{sec:Bbin}
Particular cases of quantum shuffles are the ``braided binomials'' and
``braided factorials'' (total symmetrizers)
\begin{equation*}
  \Bbin{}{r,s}=\Bbin{r+s}{(r,s)}\in\End(X^{\otimes(r+s)})
  \quad\text{and}\quad
  \Bfac{n}=\Bbin{n}{(1,\dots,1)}\in\End(X^{\otimes n}),
\end{equation*}
associated with the respective partitions $(r+s)=r+s$ and
$n=1+\dots+1$.  The braided binomials satisfy (and are in fact
determined by) the recursive relations
\begin{align}
  \label{Sh(s+1)}
  \Bbin{}{r, s + 1} &= \Bbin{}{r, s} +
  \Bbin{}{r - 1, s + 1} \Shift{(r - 1)}{\Bbraid_{1, s + 1}}
  \\
  \label{Sh(r+1)}
  \Bbin{}{r+1, s} &=
  \xShift{}{\Bbin{\shift1}{r, s}} +
  \xShift{}{\Bbin{\shift1}{r+1, s - 1}} 
  \Bbraid_{r+1,1}
\end{align}
(with $\Bbin{}{0,r}=\Bbin{}{r,0}=\id$).  Among the great many
relations satisfied by these maps, we note
\begin{gather}\label{BB}
  \Bbin{}{m + n,k} \Bbin{}{m,n} =
  \Bbin{}{m,n + k} \xShift{}{\Bbin{\shift m}{n,k}},\\
  \Bfac{m + n} =
  \Bbin{}{m, n} \Bfac{m} \Shift{m}{\Bfac{n}},
\end{gather}
and
\begin{align*}
  \sum_{i=0}^{r}
  \Bbin{}{i, s} \xShift{}{\Bbin{\shift(i + s + 1)}{r - i, t - 1}}
  \Shift{i}{\Bbraid_{r - i, s + 1}}
  &=\Bbin{}{r, s + t},
  \\
  \sum_{j=0}^{t}
  \Bbin{}{s, j} \Shift{s}{\Bbraid_{1, j}}
  &=\Bbin{}{s + 1, t}.
\end{align*}


\subsection{Examples}\label{Bint-examples}
The elements $\sigma$ of $\mathbb{S}^5_{(3,2)}$ are
\begin{multline*}
  (1, 2, 3, 4, 5),\
  (1, 2, 4, 3, 5),\
  (1, 2, 4, 5, 3),\
  (1, 4, 2, 3, 5),\
  (1, 4, 2, 5, 3),\
  (4, 1, 2, 3, 5),\\
  (1, 4, 5, 2, 3),\
  (4, 1, 2, 5, 3),\  (4, 1, 5, 2, 3),\
  (4, 5, 1, 2, 3),
\end{multline*}
and the corresponding $\sigma_{\catC}$ are given~by
\begin{multline*}
  \id,\
  \Psi_{3},\
  \Psi_{4}\Psi_{3},\
  \Psi_{2}\Psi_{3},\
  \Psi_{2}\Psi_{4}\Psi_{3},\
  \Psi_{1}\Psi_{2}\Psi_{3},\\
  \Psi_{3}\Psi_{2}\Psi_{4}\Psi_{3},\
  \Psi_{1}\Psi_{2}\Psi_{4}\Psi_{3},\
  \Psi_{1}\Psi_{3}\Psi_{2}\Psi_{4}\Psi_{3},\
  \Psi_{2}\Psi_{1}\Psi_{3}\Psi_{2}\Psi_{4}\Psi_{3}.
\end{multline*}
The sum of these 
gives $\Bbin{}{3,2}$, whose graphical representation is
\begin{multline*}
  \Bbin{}{3,2} =
  \begin{tangles}{l}
    \hstr{75}\vstr{200}\id\step[1]\id\step[1]\id\step[1]\id\step[1]\id
  \end{tangles}
  + 
  \ \begin{tangles}{l}
    \hstr{75}\vstr{200}\id\step[1]\id\step[1]\hx\step[1]\id\\
  \end{tangles}\
  + 
 \ \begin{tangles}{l}
   \hstr{75}\vstr{100}\id\step[1]\id\step[1]\hx\step[1]\id\\
   \hstr{75}\vstr{100}\id\step[1]\id\step[1]\id\step[1]\hx\\
 \end{tangles}\ 
  + 
  \ \begin{tangles}{l}
    \hstr{75}\vstr{100}\id\step[1]\id\step[1]\hx\step[1]\id\\
    \hstr{75}\vstr{100}\id\step[1]\hx\step[1]\id\step[1]\id\\
  \end{tangles}\  + 
  \ \begin{tangles}{l}
    \hstr{75}\vstr{100}\id\step[1]\id\step[1]\hx\step[1]\id\\
    \hstr{75}\vstr{100}\id\step[1]\hx\step[1]\hx\\
  \end{tangles}\  + 
  \ \begin{tangles}{l}
    \hstr{75}\vstr{67}\id\step[1]\id\step[1]\hx\step[1]\id\\
    \hstr{75}\vstr{67}\id\step[1]\hx\step[1]\id\step[1]\id\\
    \hstr{75}\vstr{67}\hx\step[1]\id\step[1]\id\step[1]\id\\
  \end{tangles}
  \\  + 
  \ \begin{tangles}{l}
    \hstr{75}\vstr{100}\id\step[1]\id\step[1]\hx\step[1]\id\\
    \hstr{75}\vstr{100}\id\step[1]\hx\step[1]\hx\\
    \hstr{75}\vstr{100}\id\step[1]\id\step[1]\hx\step[1]\id\\
  \end{tangles}\ \ 
  + 
  \ \begin{tangles}{l}
    \hstr{75}\vstr{100}\id\step[1]\id\step[1]\hx\step[1]\id\\
    \hstr{75}\vstr{100}\id\step[1]\hx\step[1]\hx\\
    \hstr{75}\vstr{100}\hx\step[1]\id\step[1]\id\step[1]\id\\
  \end{tangles}\ \ 
  + 
  \ \begin{tangles}{l}
    \hstr{75}\vstr{100}\id\step[1]\id\step[1]\hx\step[1]\id\\
    \hstr{75}\vstr{100}\id\step[1]\hx\step[1]\hx\\
    \hstr{75}\vstr{100}\hx\step[1]\hx\step[1]\id\\
  \end{tangles}\ \  + 
  \ \begin{tangles}{l}
    \hstr{75}\vstr{75}\id\step[1]\id\step[1]\hx\step[1]\id\\
    \hstr{75}\vstr{75}\id\step[1]\hx\step[1]\hx\\
    \hstr{75}\vstr{75}\hx\step[1]\hx\step[1]\id\\
    \hstr{75}\vstr{75}\id\step[1]\hx\step[1]\id\step[1]\id
  \end{tangles}
\end{multline*}

The ``braided integers'' $\Bint{r+1}=\Bbin{}{1,r}$ are explicitly
given by
\begin{alignat*}{2}
  \Bbin{}{1,1}&=\id+
  \Psi_{1}
  &&=
  \ \ \begin{tangles}{l}
    \hstr{75}\id\step\id
  \end{tangles} \ \
  +\ \ \begin{tangles}{l}
    \hstr{75}\hx
  \end{tangles}\ ,
  \\[4pt]
  \Bbin{}{1,2}&=\id
  +
  \Psi_{1}
  +
  \Psi_{2}\Psi_{1}
  &&= \ \ \begin{tangles}{l}
    \hstr{75}\vstr{134}\id\step\id\step\id
  \end{tangles} \ \
  +\ \ \begin{tangles}{l}
    \hstr{75}\vstr{134}\hx\step\id
  \end{tangles} \ \
  +\ \ \begin{tangles}{l}
    \hstr{75}\vstr{67}\hx\step\id\\
    \hstr{75}\vstr{67}\id\step\hx
  \end{tangles}\ ,
  \\[4pt]
  \Bbin{}{1,3}&=\id +
  \Psi_{1}
  +
  \Psi_{2}\Psi_{1}
  +
  \Psi_{3}\Psi_{2}\Psi_{1}
  &&= \ \ \begin{tangles}{l}
    \hstr{75}\vstr{200}\id\step\id\step\id\step\id
  \end{tangles} \ \
  +\ \ \begin{tangles}{l}
    \hstr{75}\vstr{200}\hx\step\id\step\id
  \end{tangles} \ \
  +\ \ \begin{tangles}{l}
    \hstr{75}\hx\step\id\step\id\\
    \hstr{75}\id\step\hx\step\id
  \end{tangles}\ \
  +\ \ \begin{tangles}{l}
    \hstr{75}\vstr{67}\hx\step\id\step\id\\
    \hstr{75}\vstr{67}\id\step\hx\step\id\\
    \hstr{75}\vstr{67}\id\step\id\step\hx
  \end{tangles}\ ,
\end{alignat*}

\noindent
and so on, with an obvious pattern.  The ``mirror integers'' are
\begin{alignat*}{2}
  \Bbin{}{2, 1} &=\id  
  +
  \Psi_{2}
  +
  \Psi_{1}\Psi_{2}
  &&= \ \ \begin{tangles}{l}
    \hstr{75}\vstr{134}\id\step\id\step\id
  \end{tangles} \ \
  + \ \ \begin{tangles}{l}
    \hstr{75}\vstr{134}\id\step\hx
  \end{tangles} \ \
  + \ \ \begin{tangles}{l}
    \hstr{75}\vstr{67}\id\step\hx\\
    \hstr{75}\vstr{67}\hx\step\id
  \end{tangles} \ ,
  \\[4pt]
  \Bbin{}{3, 1} &=\id
  +
  \Psi_{3}
  +
  \Psi_{2}\Psi_{3}
  +
  \Psi_{1}\Psi_{2}\Psi_{3}
  &&= \ \ \begin{tangles}{l}
    \hstr{75}\vstr{200}\id\step\id\step\id\step\id
  \end{tangles} \ \
  + \ \ \begin{tangles}{l}
    \hstr{75}\vstr{200}\id\step\id\step\hx
  \end{tangles} \ \
  + \ \ \begin{tangles}{l}
    \hstr{75}\id\step\id\step\hx\\
    \hstr{75}\id\step\hx\step\id
  \end{tangles} \ \
  + \ \ \begin{tangles}{l}
    \hstr{75}\vstr{67}\id\step\id\step\hx\\
    \hstr{75}\vstr{67}\id\step\hx\step\id\\
    \hstr{75}\vstr{67}\hx\step\id\step\id
  \end{tangles} \ ,
\end{alignat*}
and so on.

\subsection{Shuffle identities from braided Hopf algebra
  structures}\label{techn1}
We show how the axioms of braided Hopf algebras and their different
modules in Sec.~\ref{sec:H} and~\ref{sec:FUSION} are reformulated as
identities in the braid group algebra, in terms of the operators
introduced in~\bref{sec:q-shuffles}--\bref{sec:Bbin}.

\subsubsection{Braided Hopf algebra axioms in terms of shuffles}
For the product and coproduct in $\cH(X)$ introduced in~\bref{H-defs},
the fundamental Hopf-algebra axiom in Eq.~\eqref{Hopf-axiom}, read
from right to left, is equivalent to the identities
\begin{equation*}
  \sum_{i=0}^{r}\sum_{j=0}^{s}
  \Bbin{}{i, j} \xShift{}{\Bbin{\shift(i + j)}{r - i, s - j}}
  \Shift{i}{\Bbraid_{r-i,j}}
  = (r + s + 1) \Bbin{}{r, s} 
\end{equation*}
The two sums in the left-hand side are the two coactions in the
right-hand side of~\eqref{Hopf-axiom}, here in the form of the
deconcatenations of $(r)$ into $(i)\tensor(r-i)$ and of $(s)$ into
$(j)\tensor(s-i)$; the resulting groups of strands are then braided
and shuffle-multiplied as prescribed by the diagram.  The numerical
factor $r+s+1$ in the right-hand side is what remains in the shuffle
language of the coproduct
$\Delta:(r+s)\mapsto\sum_{i=0}^{r+s}(i)\tensor(r+s-i)$ in the
left-hand side of~\eqref{Hopf-axiom}

The defining properties of the antipode in~\eqref{antipode-axiom} are
equivalent to the identities
\begin{equation*}
  \sum_{s=0}^{r} \Bbin{}{s, r - s} \A_{s} = 0
  \quad\text{and}\quad
  \sum_{s=0}^{r} \Bbin{}{s, r - s} \Shift{s}{\A_{r - s}}=0
  \quad\text{for all}\quad r\geq 1.
\end{equation*}

\subsubsection{Hopf bimodule axioms in terms of
  shuffles}\label{sec:biHopf-shuffles}
The axioms of Hopf $\cH(X)$-bimod\-ules can also be stated as
identities in the braid group algebra.  For example, to reformulate
the axiom relating the left action and left coaction (the first
diagram in~\bref{Hopf-1}), we first decompose the left action
in~\eqref{r.sYt} into a sum
\begin{gather}\label{FromLeft-detail}
 (r)\dotact(s;Y;t)=
 \bigoplus_{i=0}^{r}
 \FromLeft{i,r,s,t}\,\bigl((r)\tensor(s;Y;{t})\bigr),\\[-6pt]
 \intertext{where each $\FromLeft{i,r,s,t}$ maps in a particular
   graded component:}
 \FromLeft{i,r,s,t}=\Bbin{}{i,s}\xShift{}{\Bbin{\shift(i+s+1)}{r-i,t}}
 \Shift{i}{\Bbraid_{r-i,s+1}}:
 \bigl((r)\tensor(s;Y;t)\bigr)\to(i + s; Y; r - i + t).
 \notag
\end{gather}
Then the left--left Hopf module axiom reformulates as
\begin{align*}
  \sum_{i=0}^{r}(i + s + 1)\FromLeft{i, r, s, t} &=
  \sum_{i=0}^{r}\sum_{j=0}^{s} \Bbin{}{i, j} \xShift{}{\Bbin{\shift(i
      + j)}{r - i, s - j + 1 + t}}
  \Shift{i}{\Bbraid_{r - i, j}}.\\
  \intertext{The left-module--right-comodule Hopf axiom is similarly
    expressed as the identity} \sum_{i=0}^{r} (t + r - i + 1)
  \FromLeft{i, r, s, t} &= \sum_{i=0}^{r}\sum_{j=0}^{t} \Bbin{}{i, s +
    1 + j} \xShift{}{\Bbin{\shift(i + s + 1 + j)}{r - i, t - j}}
  \Shift{i}{\Bbraid_{r - i, s + 1 + j}}.
\end{align*}

The right action in~\eqref{sYt.r} can similarly be represented in the
form that explicitly shows the graded components in the target,
\begin{align}\label{FromRight-detail}
  (s;Y;t)\dotact (r) =
  \bigoplus_{i=0}^{r}\FromRight{i,s,t,r}\bigl((s;Y;t)\tensor(r)\bigr),
\end{align}
where
\begin{equation*}  
  \FromRight{i,s,t,r}=
  \Bbin{}{s, i} \xShift{}{\Bbin{\shift(s + i + 1)}{t, r - i}}
  \Shift{s}{\Bbraid_{t + 1, i}}:(s;Y;t))\tensor(r)\to
  (s + i; Y; t + r - i)
\end{equation*}
satisfies the identities similar to those for $\FromLeft{i,r,s,t}$.

\subsubsection{The Yetter--Drinfeld axiom in terms of shuffles}
Referring to the construction of Yetter--Drinfeld modules
in~\bref{sec:YD}, we next write the ``quantum shuffle'' identity that
is a braid-algebra translation of the Yetter--Drinfeld axiom satisfied
by the adjoint action and the deconcatenation coaction.  Diagram
\eqref{yd-axiom} expressing the Yetter--Drinfeld axiom translates into
the shuffle identity
\begin{equation*}
  \sum_{i=0}^{r}\sum_{j=0}^{i + s}
  \Bbin{}{j, r - i} \Shift{j}{\Bbraid_{i + s - j + 1, r - i}} 
  \Adj_{i, s} \Shift{i}{\Bbraid_{r - i, s + 1}} =
  \sum_{i=0}^{r}\sum_{j=0}^{s}
  \Bbin{}{i, j} \Shift{(i + j)}{\Adj_{r - i, s - j}} 
  \Shift{i}{\Bbraid_{r - i, j}}.
\end{equation*}
In the left-hand side here, the two $\Bbraid$ are the two braidings in
the left-hand side of~\eqref{yd-axiom}; the sum over $i$ corresponds
to the coproduct and the sum over $j$ to the coaction in the left-hand
side of~\eqref{yd-axiom}; $\Adj_{i, s}$ is, clearly, the
$\begin{tangle} \vstr{70}\lu
  \object{\raisebox{5pt}{\kern-4pt\tiny$\blacktriangleright$}}
\end{tangle}$ action, and $\Bbin{}{j, r - i}$ is the product in the
algebra.  In the right-hand side, similarly, the sum over $i$
corresponds to the coproduct, the sum over $j$ to the coaction, and
the three operators from right to left are the braiding, the adjoint
action, and the product in the right-hand side of~\eqref{yd-axiom}.

\section{Braided Hopf algebras}\label{app:Hopf}
For basics on Hopf algebras, we refer the reader
to~\cite{[Mont],[Kassel]}.  The main definitions pertaining to braided
Hopf algebras~\cite{[Mj-braided]} (also
see~\cite{[Besp-TMF],[Besp-next],[Besp-Dr-(Bi)],[Majid-book]}) are
recalled in what follows.

\subsection{Standard notation and definitions}\mbox{}

\begin{flushright}\small
  Eine besondere Bezeichnungsweise mag unwichtig sein,\\[-3pt]
  aber wichtig ist es immer, dass diese eine\\[-3pt]
  \textit{m\"ogliche} Bezeichnungsweise ist.\\
  \textsc{L.~Wittgenstein}, \textsl{Tractatus}, 3.3421
\end{flushright}

We use the standard graphical notation for braided
categories~\cite{[Majid-book]}.  The braiding $\Psi:\cX\tensor\cY\to
\cY\tensor\cX$ of two objects $\cX$ and $\cY$ in a braided monoidal
category $\catC$ is denoted as
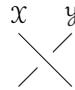
\begin{equation*}
  \begin{tanglec}
    \fobject{\cX}\step[2]\fobject{\cY}\\
    \x
  \end{tanglec}
\end{equation*}

\smallskip

\noindent
The diagrams are read from top down.  In most cases, $\cX$ and $\cY$
are omitted from the notation.

\subsubsection{Braided Hopf algebra axioms}\label{H-axioms}
For a bialgebra $\cH$ in $\catC$, the product $\cH\tensor \cH\to \cH$
and coproduct $\cH\to \cH\tensor \cH$ are respectively denoted as
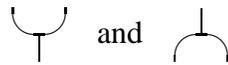
\begin{equation*}
  \begin{tanglec}
    \cu
  \end{tanglec}
  \quad\text{and}\quad
  \begin{tanglec}
    \cd
  \end{tanglec}
\end{equation*}

\smallskip

\noindent
The ``trademark'' of braided Hopf algebras is the braided bialgebra
axiom expressing the compatibility condition between the above two
maps.  It is given by the following identity for maps
$\cH\tensor\cH\to\cH\tensor\cH$:
\begin{equation}\label{Hopf-axiom}
  \begin{tanglec}
    \vstr{120}\cu\\
    \vstr{120}\cd
  \end{tanglec}
  \ = \
  \begin{tangles}{lcr}
    \cd&&\cd\\
    \id&\hx&\id\\
    \cu&&\cu
  \end{tangles}
\end{equation}

The axiom for the antipode $\A:\cH\to\cH$ is
\begin{equation}\label{antipode-axiom}
  \begin{tangles}{l}
    \cd\\
    \id\step[2]\O{\medA}\\
    \cu
  \end{tangles}\quad =\quad
  \begin{tangles}{l}
    \cd\\
    \O{\medA}\step[2]\id\\
    \cu
  \end{tangles}\ \ = \ \
  \begin{tangles}{c}
    \counit\\
    \unit
  \end{tangles}
\end{equation}
where \ $\begin{tanglec} \counit
\end{tanglec}$ \ is the counit and \ $\begin{tanglec} \unit
\end{tanglec}$ \ is the unit of~$\cH$.  It then follows that
\begin{equation*}
  \begin{tangles}{l}
    \step\O{\medA}\\
    \cd
  \end{tangles}\ = \
  \begin{tangles}{l}
    \cd\\
    \x\\
    \O{\medA}\step[2]\O{\medA}
  \end{tangles}\qquad\text{and}\qquad
  \begin{tangles}{l}
    \cu\\
    \step\O{\medA}
  \end{tangles}\ = \
  \begin{tangles}{l}
    \O{\medA}\step[2]\O{\medA}\\
    \x\\
    \cu
  \end{tangles}
\end{equation*}

\subsubsection{Modules and comodules}\label{app:mod-comod-axioms}
Left module, left comodule, right module, and right comodule
structures are expressed as
\begin{equation*}
  \begin{tangles}{ccccccc}
    \fobject{\cH}\step\fobject{\cX}&&\step\fobject{\cX}&&
    \fobject{\cX}\step\fobject{\cH}&&\fobject{\cX}\step\\
    \lu&\qquad\qquad&\ld&\qquad\qquad&\ru&\qquad\qquad&\rd\\
    \step\fobject{\cX}&&\fobject{\cH}\step\fobject{\cX}&&
    \fobject{\cX}\step&&\fobject{\cX}\step\fobject{\cH}
  \end{tangles}
\end{equation*}

\smallskip

\noindent
The left module axiom (``associativity'' of the action) is
\begin{equation*}
  \begin{tangles}{l}
    \id\step\lu\\
    \lu[2]\step
  \end{tangles}
  \ = \
  \begin{tangles}{ll}
    \cu\step\id\\
    \step\lu[2]\id
  \end{tangles}
\end{equation*}

\smallskip

\noindent
The left comodule axiom (``coassociativity'' of the coaction) is an
upside-down of this, and the right module and right comodule axioms
are formulated similarly (mirror symmetrically).

A left--left (left--right) Hopf $\cH$-module in $\catC$ is a left
module and left (right) comodule with the action and coaction
satisfying the respective axiom
\begin{equation}\label{Hopf-1}
  \begin{tangles}{c}
    \lu\\
    \ld
  \end{tangles}
  \ = \
  \begin{tangles}{lcr}
    \cd&&\ld\\
    \id&\hx&\id\\
    \cu&&\lu
  \end{tangles}
  \qquad\text{and}\quad
  \begin{tangles}{c}
    \lu\\
    \step[2]\rd
  \end{tangles}
  \ = \
  \begin{tangles}{lcr}
    \cd&&\rd\\
    \id&\hx&\id\\
    \lu[2]&&\hcu
  \end{tangles}
\end{equation}

\smallskip

\noindent
Mirror-reflected versions define the Hopf-module axioms for right
action and right coaction, and for right action and left coaction:
\begin{equation}\label{Hopf-2}
  \begin{tangles}{l}
    \ru\\
    \rd
  \end{tangles}
  \ \ = \ \
  \begin{tangles}{l}
    \rd\step\cd\\
    \id\step\hx\step[2]\id\\
    \ru\step\cu
  \end{tangles}
  \qquad\text{and}\quad
  \begin{tangles}{l}
    \step\ru\\
    \ld
  \end{tangles}
  \ \ = \ \
  \begin{tangles}{l}
    \ld\step\cd\\
    \id\step\hx\step[2]\id\\
    \hcu\step\ru[2]
  \end{tangles}
\end{equation}

\subsubsection{Hopf bimodules}\label{app:Hopf-bi}
A Hopf bimodule (``bi-Hopf module'') in a braided category $\catC$ is
a left--left Hopf module that is also a right--right, left--right, and
right--left Hopf module, as well as a bimodule and a bicomodule.
Hence, both~\eqref{Hopf-1} and~\eqref{Hopf-2} hold, and the left
action commutes with the right action, and similarly for the
coactions.

Hopf bimodules in braided categories were studied
in~\cite{[Besp-next],[Besp-Dr-(Bi)]}, where further references can be
found (Hopf bimodules appeared in~\cite{[Wor]} as ``bicovariant
bimodules''; also see~\cite{[Sch-H-YD]}).

\subsubsection{``Relative antipode.''}
Bespalov's ``relative antipode''~\cite{[Besp-TMF]} is defined for any
Hopf bimodule $\cX$ as the map $\sigma:\cX\to\cX$ given by
\begin{equation}\label{rel-anti-Hopf}
  \sigma:\ \ \
  \begin{tangles}{l}
    \ld\rd\\
    \O{\medA}\step\id\step\O{\medA}\\
    \lu\ru
  \end{tangles}
\end{equation}

\subsubsection{The ``adjoint'' action}
Any $\cH$-bimodule carries \textit{the (left) ``adjoint'' action} of
$\cH$ given by
\begin{equation}\label{adja}
  \begin{tangles}{l}
    \vstr{200}\lu\object{\raisebox{18pt}{\kern-4pt\tiny$\blacktriangleright$}}
  \end{tangles}
  \ \ = \ \
  \begin{tangles}{l}
    \hcd\step\id\\
    \vstr{80}\id\step\hx\\
    \vstr{80}\lu\step\O{\medA}\\
    \vstr{75}\step\ru
  \end{tangles}\quad
  = \
  {}\dotact{}
  ({}\dotact{}\tensor\A)
  (\id\tensor\Psi)
  (\Delta\tensor\id)
\end{equation}
(where $\dotact$ stands for both left and right actions).
The quotation marks in ``adjoint'' are to remind us that this is not
necessarily \textit{the} adjoint action of $\cH$ on itself; we omit
them in most cases, however.

\subsubsection{Yetter--Drinfeld modules}\label{sec:YD-axiom}
In a braided category, a left--left Yetter--Drinfeld (a.k.a.\
``crossed'') module is a left module under an action \ $\begin{tangle}
  \vstr{80}\lu
  \object{\raisebox{5.9pt}{\kern-4pt\tiny$\blacktriangleright$}}
\end{tangle}\ :\cH\tensor\cY\to\cY$ and left comodule under a coaction
\ $\begin{tangle}\vstr{80}\ld\end{tangle}\ :\cY\to\cH\tensor\cY$ that
are related as
\begin{equation}\label{yd-axiom}
  \begin{tangles}{l}
    \step\fobject{\cH}\step[2]\fobject{\cY}\\
    \vstr{90}\cd\step\id\\
    \vstr{50}\id\step[2]\hx\\
    \vstr{90}\lu[2] \object{\raisebox{7pt}{\kern-4pt\tiny$\blacktriangleright$}}
    \step\hd\\
    \vstr{90}\ld[2]\step\ddh\\
    \vstr{50}\id\step[2]\hx\\
    \vstr{90}\cu\step\id
  \end{tangles}\ \ = \ \
  \begin{tangles}{l}
    \step\fobject{\cH}\step[3]\fobject{\cY}\\
    \cd\step\ld\\
    \id\step[2]\hx\step\id\\
    \cu\step\lu\object{\raisebox{8.5pt}{\kern-4pt\tiny$\blacktriangleright$}}
  \end{tangles}
\end{equation}

Yetter--Drinfeld modules in braided categories were studied in
\cite{[Besp-TMF],[Besp-next]}, where further references can be found.
They form a braided monoidal category: the action on $\cY\tensor\cZ$
is diagonal, which means that
\begin{equation}\label{adj-diag}
  \begin{tangles}{l}
    \vstr{50}\id\step[2]\id\step[2]\id\\
    \lu[2]\object{\raisebox{8.2pt}{\kern15pt\tiny$\blacktriangleright
        \kern.5pt\boxed{\rule[3pt]{16pt}{0pt}}$}}\step[2]\id\\
    \vstr{50}\step[2]\id\step[2]\id
  \end{tangles} \ \
  = \ \
  \begin{tangles}{l}
    \hcd\step[1]\id\step[1]\id\\
    \vstr{80}\id\step[1]\hx\step[1]\id\\
    \lu\object{\raisebox{8pt}{\kern-4pt\tiny$\blacktriangleright$}}
    \step[1]
    \lu\object{\raisebox{8pt}{\kern-4pt\tiny$\blacktriangleright$}}
  \end{tangles}
\end{equation}
and the coaction is codiagonal (the upside down of the above diagram).
The braiding and the inverse braiding are well known to be
respectively given by (see, e.g., \cite{[Besp-TMF]})
\begin{equation}
  \label{eq:YDbraiding}
  \begin{tangles}{l}
    \ld\step[2]\id\\
    \vstr{133}\id\step\x\\
    \lu\object{\raisebox{8pt}{\kern-4pt\tiny$\blacktriangleright$}}
    \step[2]\id
  \end{tangles}
  \qquad\text{and}\qquad
  \begin{tangles}{l}
    \step[1]\xx\\
    \ld\step[2]\id\\
    \vstr{67}\hxx\step[2]\id\\
    \id\step\O{\medA^{\scriptscriptstyle-1}}\step[1]\dd\\
    \id\step\lu\object{\raisebox{8pt}{\kern-4pt\tiny$\blacktriangleright$}}
  \end{tangles}
\end{equation}
(we always assume the antipode to be bijective, i.e., $\A^{-1}$ to
exist).

For a Yetter--Drinfeld module $\cY$, the ``squared relative
antipode''~\cite{[Besp-TMF]} is the map $\sigma_2:\cY\to\cY$ given by
\begin{equation}\label{rel-anti-YD}
  \sigma_2:\ \ \
  \begin{tangles}{l}
    \ld\\
    \O{\medA}\step\id\\
    \lu\object{\raisebox{8.5pt}{\kern-4pt\tiny$\blacktriangleright$}}
  \end{tangles}
\end{equation}
(see~\cite{[Besp-next]} for its properties and use).

\subsubsection{Category equivalence} For a Hopf algebra $\cH$ in a
braided monoidal category with split idempotents, the category
$^{\cH}_{\cH}\mathcal{M}^{\cH}_{\cH}$ of Hopf $\cH$-bimodules is
braided monoidal equivalent to the category $^{\cH}_{\cH}\YDname$ of
Yetter--Drinfeld $\cH$-modules.  The functor
$^{\cH}_{\cH}\mathcal{M}^{\cH}_{\cH}\to{}^{\cH}_{\cH}\YDname$ is given
by taking right coinvariants, and the inverse functor, by induction
(see~\cite{[Besp-TMF],[Besp-next]} and the references therein,
\cite{[Sch-H-YD]} in particular).

\subsection{The $\protect\mapX$ map of Hopf bimodules}
As before, we fix a Hopf algebra $\cH$ in a braided monoidal
category~$\catC$ and consider Hopf bimodules over~$\cH$.

The tensor product of two Hopf bimodules can be given the structure of
a Hopf bimodule under the ``left'' left action, codiagonal left
coaction, ``right'' right action, and codiagonal right coaction,
respectively given by
\begin{equation}\label{Hopf-standard}
  \ast_\Left:\quad\begin{tangles}{rrr}
    \id&\id&\id\\
    &\lu&\step\id\\
    &\id&\step\id
  \end{tangles}
  \qquad\qquad
  \delta'_\Left:\quad\begin{tangles}{lrr}
    \ld&&\ld\\
    \id&\hx&\id\\
    \hcu&\id&\id
  \end{tangles}
  \qquad\qquad
  \ast_\Right:\quad
  \begin{tangles}{lll}
    \id&\step\id&\step\id\\
    &\id&\ru\\
    &\id&\id
  \end{tangles}
  \qquad\qquad
  \delta'_\Right:\quad
  \begin{tangles}{llr}
    \rd&&\rd\\
    \id&\hx&\id\\
    \id&\id&\hcu
  \end{tangles}
\end{equation}
(see, e.g.,~\cite{[Besp-next],[Besp-Dr-(Bi)]}).

\subsubsection{}
For any two Hopf bimodules $\cX$ and $\cY$, we define
$\mapX(\cX\tensor \cY)$ as
\begin{equation}\label{chi-map}
  \begin{tangles}{llr}
    \fobject{\cX}&&\fobject{\cY}\\
    \rd&&\ld\\
    \id&\hx&\id\\
    \ru&&\lu
  \end{tangles}
\end{equation}
The map $\mapX$ intertwines the actions and coactions
in~\eqref{Hopf-standard} respectively with the diagonal left action,
``left'' left coaction, diagonal right action, and ``right'' right
coaction:
\begin{equation}\label{Hopf-dot}
  \cdot_\Left:\quad\begin{tangles}{lllr}
    \hcd&&\id&\id\\
    \id&\hx&&\id\\
    \lu&&\hlu[2]
  \end{tangles}
  \qquad\qquad
  \delta_\Left:\quad
  \begin{tangles}{rrr}
    &\id&\step\id\\
    &\ld&\step\id\\
    \id&\id&\id
  \end{tangles}
  \qquad\qquad
  \cdot_\Right:\quad
  \begin{tangles}{lrrr}
    \id\step&\id&&\hcd\\
    \id&&\hx&\id\\
    \hru[2]&&&\ru
  \end{tangles}
  \qquad\qquad
  \delta_\Right:\quad
  \begin{tangles}{lll}
    \id\step&\id\\
    \id\step&\rd\\
    \id&\id&\id
  \end{tangles}
\end{equation}
i.e.,
\begin{multline*}
  \mapX\ccirc\ast_\Left=\cdot_\Left\ccirc(\id\tensor\mapX),
  \quad
  (\id\tensor\mapX)\ccirc\delta'_\Left=\delta_\Left\ccirc\mapX,
  \\*
  \mapX\ccirc\ast_\Right=\cdot_\Right\ccirc(\mapX\tensor\,\id),
  \quad
  (\mapX\tensor\,\id)\ccirc\delta'_\Right=\delta_\Right\ccirc\mapX.
  \kern-10pt 
\end{multline*}
This is shown by simple manipulations with the diagrams.  For example,
to prove the second identity, we start with its right-hand side and
use the relevant Hopf-module axiom and then the coaction property and
the commutativity of left and right coactions:
\begin{equation}\label{co-id}
  \delta_\Left\ccirc\mapX = \ \
  \begin{tangles}{l}
    \step\hrd\step\hld\\
    \step\id\step[.5]\hx\step[.5]\id\\
    \vstr{200}\ld\hru\step\hlu
  \end{tangles}
  \ \ \ = \ \
  \begin{tangles}{l}
    \step\rd[2]\step\ld\\
    \ld\step[2]\hx\step\id\\
    \id\step\id\step[1]\cd\lu\\
    \id\step\hx\step\dd\step\id\\
    \hcu\step\ru\step[2]\id
  \end{tangles}
  \ \ \ = \ \
  \begin{tangles}{l}
    \ld\step[1]\ld[2]\step[-.5]\hld\\
    \id\step\hx\step[1.5]\id\step[.5]\id\\
    \id\step\id\step\hrd\step\id\step[.5]\id\\
    \id\step\id\step\id\step[.5]\hx\step[.5]\id\\
    \hcu\step[1]\hru\step\hlu
  \end{tangles}
  \ \ = (\id\tensor\mapX)\ccirc\delta'_\Left.
\end{equation}

\smallskip

\subsubsection{}It follows that $\mapX$ is associative in the
following sense:
\begin{equation*}
  \mapX(\mapX(\cX\tensor \cY)\tensor \cZ)=\mapX(\cX\tensor\mapX(\cY\tensor \cZ)),
\end{equation*}
where \textit{each of the inner $\mapX(U,V)$ is understood as a Hopf
  bimodule under
  $(\cdot_\Left,\delta_\Left,\cdot_\Right,\delta_\Right)$}.  This is
also readily verified diagrammatically, with both sides of the last
formula being given by
\begin{equation*}
  \begin{tangles}{llr}
    \vstr{67}\rd[2]&\vstr{67}\ld\rd&\vstr{67}\step\ld[2]\\
    \vstr{67}\rd\step\hx&\vstr{67}\step\id&\vstr{67}\hx\step\ld\\
    \vstr{67}\id\step\hx\step&\vstr{67}\hx\step&\vstr{67}
    \id\step\hx\step\id\\
    \vstr{67}\ru\step\id\step\id&\vstr{67}\step\hx&\vstr{67}
    \step\id\step\lu\\
    \vstr{67}\id\step[2]\id&\vstr{67}\hx&\vstr{67}\hx\step[2]\id\\
    \vstr{67}\id\step[2]\hx&\vstr{67}\step\id&\vstr{67}\id\step\d\step\id\\
    \vstr{67}\ru[2]&\vstr{67}\lu\ru&\vstr{67}\lu
  \end{tangles}
  \quad\text{where}\quad
  \begin{tangles}{c}
    \ld\rd
  \end{tangles}\ = \
  \begin{tangles}{l}
    \vstr{67}\ld\\
    \vstr{67}\id\step\rd
  \end{tangles}
  \quad\text{and}\quad
  \begin{tangles}{c}
    \lu\ru
  \end{tangles}\ = \
  \begin{tangles}{r}
    \vstr{67}\lu\step\id\\
    \vstr{67}\ru
  \end{tangles}
\end{equation*}

\noindent
We can therefore define $\mapX(\cX\tensor \cY\tensor \cZ\tensor\dots)$
by nesting.

\subsection{The $\protect\mapI$ (\texttt{\normalfont[\i:]}) map of right
  coinvariants}\label{sec:mapI}
For right coinvariants in Hopf bimodules, $\mapX$ in~\eqref{chi-map}
reduces to the map $\mapI:\cX\tensor\cY\to\cX\tensor\cY$ given by the
diagram
\begin{equation*}
  \mapI:\ \ 
  \begin{tangles}{l}
    \id\step[1]\ld\\
    \ru\step[1]\id
  \end{tangles}
\end{equation*}
which yields right coinvariants and which inherits the property to
intertwine $\delta'_\Left$ and $\delta_\Left$, i.e.,
\begin{equation}\label{comodule-morphism}
  \begin{tangles}{l}
    \vstr{75}\ld\step[1]\ld[1]\\
    \vstr{75}\id\step[1]\hx\step[1]\d\\
    \hcu\step[1]\id\step\ld\\
    \vstr{75}\step[.5]\id\step[1.5]\ru\step\id
  \end{tangles}
  \ \ = \ \
  \begin{tangles}{l}
    \step\id\step\ld\\
    \step\ru\step\id\\
    \ld\step[2]\id
  \end{tangles}
\end{equation}

\noindent
as an immediate application of a Hopf-module axiom, with both sides
being equal to
\begin{equation*}
  \begin{tangles}{l}
    \vstr{67}\ld\step[1]\step[.5]\ld\\
    \vstr{90}\id\step\id\step[1]\hcd\step[.5]\id\\
    \vstr{67}\id\step[1]\hx\step[1]\id\step[.5]\id\\
    \vstr{67}\hcu\step[1]\ru\step[.5]\id
  \end{tangles}\ . 
\end{equation*}
We also have
$\mapI\ccirc\ast_\Right=\cdot_\Right\ccirc(\mapI\tensor\,\id)$.  In
verifying the associativity, we then find that
\begin{equation}\label{eq:II}
  \mapI(\mapI(\cX\tensor \cY)\tensor
  \cZ)=\mapI(\cX\tensor\mapI(\cY\tensor \cZ))
  =\ \
  \begin{tangles}{l}
    \id\step[.5]\ld\step[2]\hld\\
    \vstr{67}\hru\step[1]\id\step[1]\dd\step[.5]\id\\
    \vstr{67}\id\step[1.5]\hx\step[1]\hld\\
    \vstr{67}\id\step[.5]\dd\step[1]\id\step[1]\id\step[.5]\id\\
    \hru\step[2]\id\ru\step[.5]\id
  \end{tangles}
\end{equation}

This allows defining the action of $\mapI$ on multiple tensor products
of right coinvariants in Hopf bimodules by nesting.

Right coinvariants in Hopf bimodules can be viewed as Yetter--Drinfeld
modules under the adjoint action.  The $\mapI$ map intertwines the
diagonal action on tensor products, Eq.~~\eqref{adj-diag}, as follows:
\begin{equation}\label{adja-intertwine}
  \begin{tangles}{l}
    \hcd\step[1]\id\step\id\\
    \id\step[1]\hx\step\id\\
    \lu\object{\raisebox{8pt}{\kern-4pt\tiny$\blacktriangleright$}}
    \step[1]\lu\object{\raisebox{8pt}{\kern-4pt\tiny$\blacktriangleright$}}\\
    \step[1]\id\step\ld\\
    \step[1]\ru\step\id
  \end{tangles}
  \ \ = \ \
  \begin{tangles}{l}
    \vstr{80}\step[.5]\id\step[2.5]\id\step[1]\ld\\
    \hcd\step[2]\ru\step[1]\id\\
    \vstr{70}\id\step[1]\x\step[2]\dh\\
    \lu\step[1.5]\hcd\step[2]\id\\
    \vstr{70}\step[1]\id\step[1.5]\id\step[1]\x\\
    \vstr{70}\step[1]\id\step[1.5]\lu\step[1]\cd\\
    \vstr{70}\step[1]\id\step[2.5]\id\step[1]\O{\medA}\step[2]\id\\
    \vstr{50}\step[1]\id\step[2.5]\ru\step[1]\dd\\
    \vstr{70}\step[1]\d\step[1.5]\x\\
    \vstr{70}\step[2]\hd\step[1]\O{\medA}\step[2]\id\\
    \vstr{60}\step[2.5]\ru\step[2]\id\\
  \end{tangles}
\end{equation}

\smallskip

\begin{figure}[tbhp]
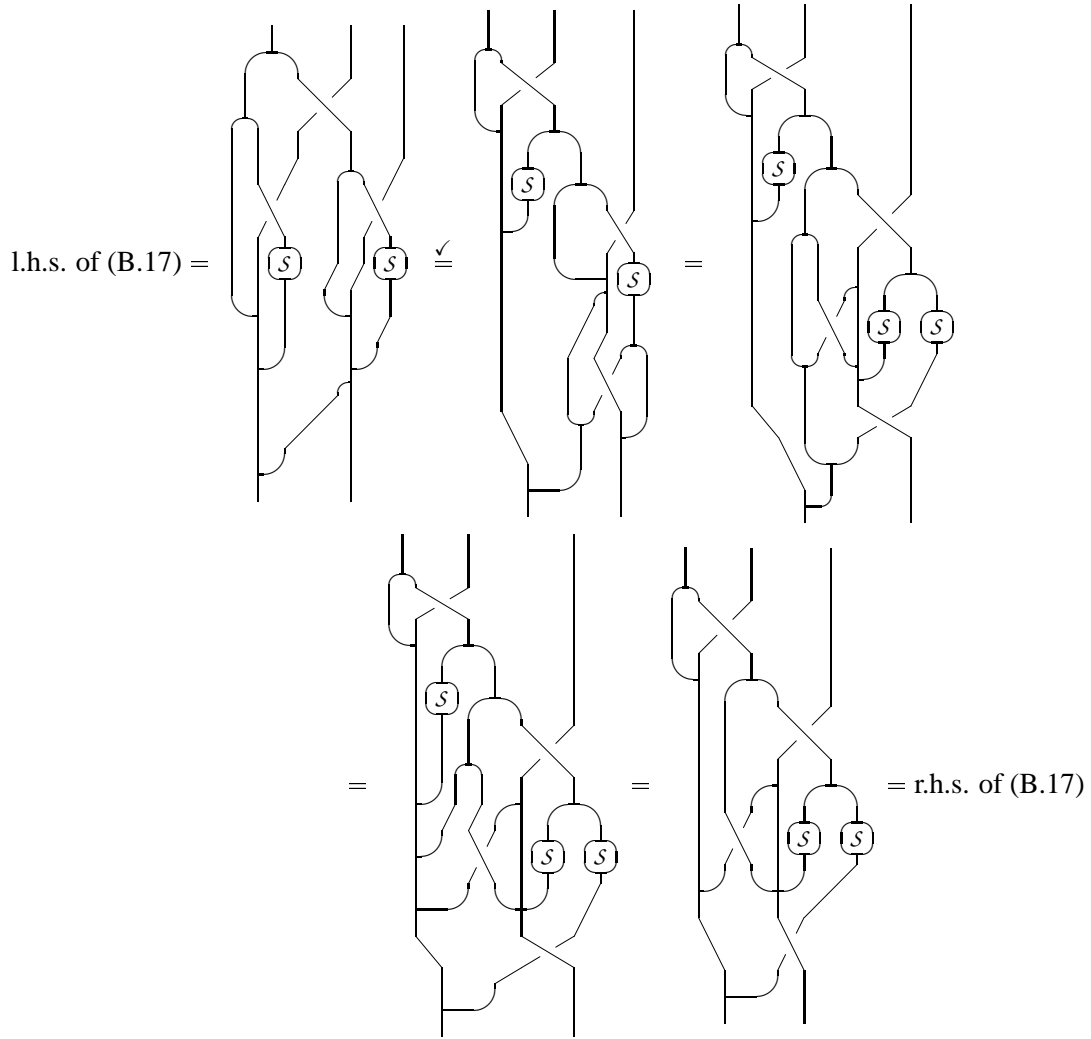
\small
  \centering
  \begin{multline*}
    \text{l.h.s. of~\eqref{adja-intertwine}} = \ \
    \begin{tangles}{l}
      \step[.5]\cd\step[2]\id\step[2]\id\\
      \hcd\step[1.5]\x\step[2]\id\\
      \id\step\id\step\ddh\step[1.5]\hcd\step[1]\ddh\\
      \id\step\hx\step[2]\id\step\hx\\
      \id\step\id\step\O{\medA}\step[1.5]\ddh\step[.5]\ddh\step\O{\medA}\\
      \lu\step\id\step[1.5]\lu\step\ddh\\
      \step[1]\ru\step[2.5]\ru\step[-1.5]\hld\\
      \step[1]\id\step[2]\ne2\step[.5]\id\\
      \step[1]\ru\step[2.5]\id
    \end{tangles}
    \ \ \ \stackrel{\checkmark}{=} \ \ 
    \begin{tangles}{l}
      \hcd\step[2]\id\step[3]\id\\
      \vstr{80}\id\step[1]\x\step[3]\id\\
      \lu\step\cd\step[2]\id\\
      \step[1]\id\step[1]\O{\medA}\step\cd\step\id\\
      \vstr{80}\step[1]\ru\step\id\step[2]\hx\\
      \step[1]\id\step[2]\lu[2]\step[-.5]\hld\step\O{\medA}\\
      \step[1]\id\step[2.5]\dd\ddh\step[.5]\hcd\\
      \step[1]\id\step[2.5]\id\step\hx\step\id\\
      \step[1]\d\step[1.5]\hcu\step\ru\\
      \step[2]\ru[2]\step[1.5]\id
    \end{tangles}
    \ \ \ \  = \ \ 
    \begin{tangles}{l}
      \hcd\step[2]\id\step[4]\id\\
      \vstr{60}\id\step\x\step[4]\id\\
      \lu\step\cd\step[3]\id\\
      \step[1]\id\step[1]\O{\medA}\step\cd\step[2]\id\\
      \step[1]\ru\step[.5]\hcd\step[1.5]\x\\
      \step[1]\id\step[1.5]\id\step[1]\id\step[1]\hld\step[1]\cd\\
      \step[1]\id\step[1.5]\id\step[1]\hx\step[.5]\id\step\O{\medA}\step[2]\O{\medA}\\
      \step[1]\id\step[1.5]\hcu\step\hlu\ru\step\dd\\
      \vstr{60}\step[1]\d\step[1]\id\step[2]\x\\
      \step[2]\d\cu\step[2]\id\\
      \vstr{60}\step[3]\ru\step[3]\id
    \end{tangles}
    \\
    \ \ = \ \ 
    \begin{tangles}{l}
      \hcd\step[2]\id\step[4]\id\\
      \vstr{60}\id\step[1]\x\step[4]\id\\
      \lu[1]\step\cd\step[3]\id\\
      \step[1]\id\step[1]\O{\medA}\step\cd\step[2]\id\\
      \step[1]\id\step[1]\id\step[.5]\hcd\step[1.5]\x\\
      \step[1]\ru\ddh\step[.5]\ddh\step[.5]\ld\step[1]\cd\\
      \step[1]\ru\step[1]\hx\step[1]\id\step[1]\O{\medA}\step[2]\O{\medA}\\
      \step[1]\ru[2]\step[1]\lu\ru\step\dd\\
      \vstr{60}\step[1]\d\step[3]\x\\
      \vstr{30}\step[2]\id\step[2]\dd\step[2]\id\\
      \step[2]\ru[2]\step[3]\id
    \end{tangles}
    \ \ \ = \ \
    \begin{tangles}{l}
      \hcd\step[2]\id\step[3]\id\\
      \id\step[1]\x\step[3]\id\\
      \lu[1]\step\cd\step[2]\id\\
      \step[1]\id\step[1]\id\step[2]\x\\
      \step[1]\id\step[1]\id\step[1]\ld\step[1]\cd\\
      \step[1]\id\step[1]\hx\step\id\step[1]\O{\medA}\step[2]\O{\medA}\\
      \step[1]\ru\step[1]\lu\ru\step\ne{2}\\
      \step[1]\d\step[2]\hx\\
      \step[2]\ru[2]\step\id
    \end{tangles}
    \ \ \ = \text{r.h.s. of~\eqref{adja-intertwine}}
  \end{multline*}
  \caption{\small Proof of~\eqref{adja-intertwine}.  The
    $\stackrel{\checkmark}{=}$ equality and the next one use two
    Hopf-module axioms and the braided anti-automorphism property of
    the antipode.  The\ equality connecting the two lines of the
    formula involves the right action associativity property, and the
    next equality uses coassociativity to isolate a ``bubble'' as
    in~\eqref{antipode-axiom}, which is then
    eliminated.}\label{fig:weary}
\end{figure}

\noindent
This is shown by direct calculation in Fig.~\ref{fig:weary}.
Everything below the top $\mapI$-diagram in the right-hand side
of~\eqref{adja-intertwine} shows how the Hopf algebra acts on
$\mapI(\cX\tensor\cY)$.  

The pattern extends to multiple tensor
products in the manner that is entirely clear from the example of a
triple product: the $\mapI$ map ``squared,'' i.e., \textit{the diagram
  in~\eqref{eq:II}, intertwines the diagonal adjoint action on a
  triple tensor product with the action given by}
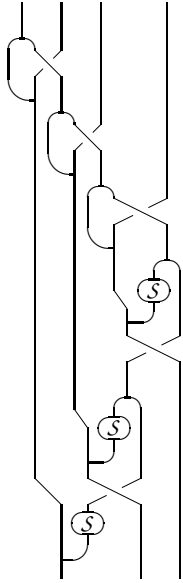
\begin{equation*}
  \begin{tangles}{l}
    \vstr{90}\step[1]\hcd\step[1]\id\step[1.5]\id\step[2.5]\id\\
    \vstr{50}\step[1]\id\step[1]\hx\step[1.5]\id\step[2.5]\id\\
    \vstr{90}\step[1]\lu\step[.5]\hcd\step[1]\id\step[2.5]\id\\
    \vstr{50}\step[2]\id\step[.5]\id\step[1]\hx\step[2.5]\id\\
    \vstr{90}\step[2]\id\step[.5]\lu\step[.5]\hcd\step[2]\id\\
    \vstr{50}\step[2]\id\step[1.5]\id\step[.5]\id\step[1]\x\\
    \vstr{90}\step[2]\id\step[1.5]\id\step[.5]\lu\step[1.5]\hcd\\
    \vstr{70}\step[2]\id\step[1.5]\id\step[1.5]\dh\step[1]\O{\medA}\step[1]\id\\
    \vstr{50}\step[2]\id\step[1.5]\id\step[2]\ru\step[1]\id\\
    \vstr{50}\step[2]\id\step[1.5]\id\step[2]\x\\
    \vstr{90}\step[2]\id\step[1.5]\id\step[1.5]\hcd\step[1.5]\id\\
    \vstr{70}\step[2]\id\step[1.5]\hd\step[1]\O{\medA}\step[1]\id\step[1.5]\id\\
    \vstr{60}\step[2]\id\step[2]\ru\step\id\step[1.5]\id\\
    \vstr{50}\step[2]\d\step[1]\x\step[1.5]\id\\
    \vstr{70}\step[3]\id\step[1]\O{\medA}\step[2]\id\step[1.5]\id\\
    \vstr{70}\step[3]\ru[1]\step[2]\id\step[1.5]\id
  \end{tangles}
\end{equation*}
The extension to $n$-fold tensor products is now entirely obvious (and
gives a fully braided version of the ``unexpected'' action
in~\cite{[Mj-hopf-in-braided]}).  The action on the \textit{last},
rightmost tensor factor is the adjoint in all cases.

\section{Yetter--Drinfeld braiding and diagonal
  braiding}\label{sec:YD-and-diag}
An important class of braidings is provided by Yetter--Drinfeld
categories.

Let $H$ be an ordinary Hopf algebra (with bijective antipode~$\hA$).
A left--left Yetter--Drinfeld $H$-module $R$ is by definition a left
$H$-module (with an action $H\tensor R\to R: h\tensor r\mapsto
h\leftii r$) and left $H$-comodule (with coaction $R\to H\tensor R:
r\mapsto r\mone\tensor r\zero$) satisfying the defining axiom
\begin{equation*}
  (h'\leftii r)\mone h''\tensor(h'\leftii r)\zero
  =h' r\mone\tensor h''\leftii r\zero
\end{equation*}
(which is just~\eqref{yd-axiom} with the braiding given by
transposition), or equivalently
\begin{equation}\label{YD-H}
  (h\leftii r)\mone\tensor(h\leftii r)\zero
  =h' r\mone \hA(h''')\tensor h''\leftii r\zero.
\end{equation}
The category $\HHyd$ of Yetter--Drinfeld $H$-modules is monoidal and
braided, with the braiding and the inverse given by
(cf.~\eqref{eq:YDbraiding})
\begin{align} \label{YD-braiding}
  \Psi: \cU\tensor \cV&\to \cV\tensor \cU:
  u\tensor v\mapsto u\mone\leftii v\tensor u\zero,\\
  \Psi^{-1}: \cV\tensor \cU&\to \cU\tensor \cV: v\tensor u\mapsto
  v\zero\tensor \hA^{-1}(v\mone)\leftii u.
\end{align}

A Hopf algebra $\cR$ in $\HHyd$ is a Yetter--Drinfeld $H$-module and a
braided Hopf algebra whose braiding is that of $\HHyd$ and whose
operations are $\HHyd$ morphisms (this notion can be defined purely
categorically).

Any braided Hopf algebra whose braiding is rigid can be realized as a
Hopf algebra in $\HHyd$ for some ordinary Hopf algebra $H$ (which is
by far not unique)~\cite{[T-survey]}.

For a Hopf algebra $\cR$ in $\HHyd$, we consider left $\cR$ modules
and comodules $\cY\in\HHyd$.  Because any such $\cY$ is also an $H$
module and comodule, we have to clearly distinguish the different
actions and coactions.  We write $r\leftact y$ and $h.y$ for the $\cR$
and $H$ actions, and $y\mapsto y\bmone\tensor y\bzero\in\cR\tensor
\cY$ and $y\mapsto y\mone\tensor y\zero\in H\tensor\cY$ for the
coactions.  That the $\cR$-coaction is an $H$-comodule morphism is
then expressed as
\begin{equation}\label{eq:Hcom-morph}
  y\mone\!\tensor
  y\zero{}\bmone\tensor y\zero{}\bzero=y\bmone{}\mone
  y\bzero{}\mone\tensor y\bmone{}\zero \tensor y\bzero{}\zero\in
  H\tensor\cR\tensor\cY.
\end{equation}
It is also an $H$-module morphism: $(h.y)\bmone\tensor(h.y)\bzero
=(h'\leftii y\bmone)\tensor h''.y$.

If we assume that the braided linear spaces $X$ and $Y$ in
Secs.~\ref{sec:H} and~\ref{sec:FUSION} are objects in $\HHyd$, then
$\Psi$ is given by~\eqref{YD-braiding} and there is an $H$-coaction
$x\mapsto x\mone\tensor x\zero$ on $X$, and similarly on $Y$, and then
\begin{align*}
  (x_1\tensor\dots\tensor x_n)\mone \tensor(x_1\tensor\dots\tensor
  x_n)\zero &=x_1{}\mone\dots x_n{}\mone \tensor
  (x_1{}\zero\tensor\dots\tensor x_n{}\zero).
  \\[-4pt]
  \intertext{Also, $(~)\bup{1}\tensor(~)\bup{2}$ and
    $(~)\bmone\tensor(~)\bzero$ are now the deconcatenation coproduct
    and coaction:}
  (x_1\tensor\dots\tensor x_n)\bup{1}
  {}\tensor(x_1\tensor\dots\tensor x_n)\bup{2} =&\sum_{i=0}^n
  (x_1\tensor\dots\tensor x_i)\tensor (x_{i+1}\tensor\dots\tensor x_n)
  \\[-6pt]
  \intertext{and} (x_1\tensor\dots\tensor x_n\tensor y)\bmone
  {}\tensor(x_1\tensor\dots\tensor x_n\tensor y)\bzero =&\sum_{i=0}^n
  (x_1\tensor\dots\tensor x_i)\tensor (x_{i+1}\tensor\dots\tensor
  x_n\tensor y).
\end{align*}
The adjoint action
\begin{equation*}
  \Adja_{1,s}:X\tensor X^{\otimes s}\tensor Y\to X^{\otimes(s+1)}\tensor Y
\end{equation*}
then evaluates as (for $x\in X$ and $b\in X^{\tensor s}\tensor Y$)
\begin{equation*}
  \Adja_{1,s}\bigl(x\tensor b\bigr)
  =x\mone\leftii b\bmone\tensor x\zero\tensor b\bzero
  - x\mone{}'\leftii\Bigl(b\bmone\tensor b\bzero{}\mone\hA(x\mone{}'')
  \leftii x\zero\tensor b\bzero{}\zero\Bigr).
\end{equation*}

Diagonal braiding in Yetter--Drinfeld categories is associated with
commutative cocommutative $H$, i.e., group algebras of (finite)
Abelian groups. This case has been considered in numerous papers, and
in particular in the ``Nichols'' context---in the framework of
Andruskiew\-itsch and Schneider's classification program
(\cite{[AS-pointed],[AS-onthe]}, and the references therein).  We
follow these papers in briefly recalling the relevant points about
diagonal braiding.

For a finite group $\Gamma$, left Yetter--Drinfeld
$k\Gamma$-modules are $\Gamma$-graded vector spaces
$X=\bigoplus_{g\in\Gamma}X_g$ with a $\Gamma$ action such that $g
X_{h} \subset X_{g h g^{-1}}$; the left comodule structure is here
given by $\delta x= 
g\tensor x$ for all $x\in X_g$.  If, moreover, $\Gamma$ is Abelian,
then Yetter--Drinfeld $k\Gamma$-modules are just $\Gamma$-graded
vector spaces with an action of $\Gamma$ on each~$X_g$.  The action
is then diagonalizable, and hence
$X=\bigoplus_{\chi\in\widehat\Gamma} X^{\chi}$, where $X^{\chi} =
\{x\in X\mid g x = \chi(g) x \text{ for all } g\in\Gamma\}$ and
$\widehat\Gamma$ is the group of characters of $\Gamma$. \ Then
\begin{equation*}
  X=\bigoplus_{g\in\Gamma,\;\chi\in\widehat\Gamma} X^{\chi}_g,
\end{equation*}
where $X^{\chi}_g=X^{\chi}\cap X_g$. \ Hence, each Yetter--Drinfeld
$k\Gamma$-module $X$ has a basis $(x_i)$ such that, for some
$g_i\in\Gamma$ and $\chi_i\in\widehat\Gamma$, \ $\delta x_i =
g_i\tensor x_i$ and $g x_i = \chi_i(g) x_i$ for all $g$. \ The
braiding~\eqref{YD-braiding} then takes the form
\eqref{diag-braiding} with $q_{ij}=\chi_j(g_i)$.\footnote{Diagonal
  braiding is also said to be of the Abelian group type in the
  terminology in~\cite{[AS-pointed]}.  Whenever the $q_{ij}$ are roots
  of unity, the braiding is of the \textit{finite group type}, which
  means that the group generated by the~$g_i$ in $\End X$ is finite.}

\parindent=0pt

\end{document}

